\documentclass[11pt,reqno]{amsart}

\usepackage[utf8]{inputenc}
\usepackage[T1]{fontenc}
\usepackage[french, english]{babel}
\usepackage{enumitem}

\usepackage{amsthm}
\usepackage{amsmath}
\usepackage{amssymb}
\usepackage{mathrsfs}
\usepackage{bbm}
\usepackage{graphicx}
\usepackage{listings}
\usepackage{mathrsfs}
\usepackage{pict2e}
\usepackage{stmaryrd}
\usepackage{mathtools}
\usepackage{enumitem}

\renewcommand{\Re}{\mathrm{Re}}
\renewcommand{\Im}{\mathrm{Im}}

\usepackage[all,cmtip]{xy}

\usepackage{multirow}

\usepackage{color}
\definecolor{marin}{rgb}   {0.,   0.1,   0.5} 
\definecolor{rouge}{rgb}   {0.8,   0.,   0.} 
\definecolor{sepia}{rgb}   {0.4,   0.25,   0.} 
\definecolor{mag}{rgb}   {0.3,   0,   0.3} 

\usepackage[colorlinks,citecolor=sepia,linkcolor=marin,urlcolor=mag ,bookmarksopen,bookmarksnumbered]{hyperref}

\newtheorem{theorem}{Theorem}[section]
\newtheorem{corollary}[theorem]{Corollary}
\newtheorem{lemma}[theorem]{Lemma}
\newtheorem{proposition}[theorem]{Proposition}
\newtheorem{definition}[theorem]{Definition}
\newtheorem{remark}[theorem]{Remark}

\newtheorem{Assumptions}{Assumptions}

\newcommand{\N}{\mathbb{N}}

\newcommand{\Gc}{\mathcal{G}}

\newcommand{\Lc}{\mathcal{L}}

\newcommand{\Dc}{\mathcal{D}}

\newcommand{\R}{\mathbb{R}}

\newcommand{\C}{\mathbb{C}}

\newcommand{\T}{\mathbb{T}}
\newcommand{\Z}{\mathbb{Z}}
\newcommand{\Tc}{\mathcal{T}}

\newcommand{\Sc}{\mathcal{S}}

\newcommand{\K}{\mathcal{K}}
\renewcommand{\div}{\mathrm{div}}

\newcommand{\Dir}{\mathsf D}
\newcommand{\Neu}{\mathsf N}

\newcommand{\Ac}{\mathcal{A}}

\newcommand{\dt}{\partial_t}
\newcommand{\dx}{\partial_x}
\newcommand{\jb}[1]{\langle #1\rangle}

\newcommand{\les}{\lesssim}

\newcommand{\D}{\mathcal{D}}
\newcommand{\ft}{\widehat}
\newcommand{\F}{\mathcal{F}}
\newcommand{\M}{\mathcal{M}}

\newcommand{\cj}[1]{\overline{#1}}

\renewcommand{\P}{\mathbf{P}}

\newcommand{\scal}{\mathrm{scal}}

\newcommand{\reff}{\mathrm{ref}}
\newcommand{\Q}{\mathbf{Q}}

\DeclareMathOperator{\Diff}{\textsf{Diff~}}

\renewcommand{\vec}{\textsf{Vec~}}

\usepackage{scalerel} % to rescale the paraproduct symbols
\usepackage{tikz}

\begin{document}
\title[Low-regularity nonlinear Schr\"odinger equation with bilinear interactions]{Low-regularity evolution and control for the nonlinear Schr\"odinger equation with singular bilinear interactions}

%\title[]{Small time exact bilinear control of periodic NLS at low regularity}

\author{Rémi Buffe}

\address{\small{Université de Lorraine, CNRS, IECL, F-54000 Nancy, France}}

\email{remi.buffe@univ-lorraine.fr}

\author{Alessandro Duca}

\address{\small{Université de Lorraine, CNRS, Inria, IECL, F-54000 Nancy, France}}

\email{alessandro.duca@inria.fr}

\author{Tristan Robert}

\address{\small{Université de Lorraine, CNRS, IECL, F-54000 Nancy, France}}

\email{tristan.robert@univ-lorraine.fr}

%\author{Laurent Thomann}
%
%\address{\small{Université de Lorraine, CNRS, IECL, F-54000 Nancy, France}}
%
%\email{laurent.thomann@univ-lorraine.fr}

%
%\subjclass[2020]{35Q53,35Q55,60L50}

\begin{abstract}
 
We consider the one-dimensional nonlinear Schrödinger equation (NLS) $i\dt \psi -\dx^2 \psi =  (u(t)\cdot  \mu) \psi +\kappa |\psi|^{2m}\psi,$ on the torus under the presence of singular, time-dependent bilinear external fields $(u(t)\cdot  \mu)$. We prove well-posedness in $H^s(\T)$ for low regularity indices $s>\max\big\{\tfrac12-\tfrac1m,0\big\}$ ($s>0$ when $\kappa=0$), under sharp assumptions on the singular potentials, covering both the linear and nonlinear cases.
We establish small-time exact controllability in $H^{\frac12+}(\T)$ for a large class of states, by combining saturation techniques with the moment method. The result extends previous bilinear exact controllability results for linear Schrödinger equations, which are typically obtained in large time and under high regularity assumptions.
Crucial to our approach is the role of the singular time-dependent potentials, which enable the recovery of controllability properties in small time and at low regularity, overcoming obstructions present in the smooth setting.
\end{abstract}

%\date{\today}
%%
%
\maketitle

\vspace{-5mm}

\tableofcontents

\section{Introduction}

While much of the literature on partial differential equations concerns additive or boundary forcing terms, many evolution models arising in quantum mechanics involve multiplicative interactions. A representative example is a quantum particle on the torus interacting with an external electromagnetic field, whose spatial profile is fixed, whereas its time-dependent intensity acts as a scalar modulation of the dynamics.
The dynamics of the particle can be modeled by a nonlinear equation of the form
\begin{equation}\label{NLS}\tag{NLS}
\begin{split}
\begin{cases}
i\dt \psi -\dx^2 \psi = (u(t)\cdot\mu(x))\psi +\kappa |\psi|^{2m}\psi,~~&(t,x)\in\R\times\T,\\
\psi(t=0,x)=\psi_0(x),~~&x\in\T,
\end{cases}
\end{split}
\end{equation}
where $\T=\R/(2\pi\Z)$ is the torus representing the medium, the operator $-\dx^2$ is the standard periodic Laplacian, and $u(t)\cdot\mu(x)$ is the external potential field. Here, the fixed vector-valued function $\mu:\T\rightarrow \R^d$ encodes the spatial action of the potential, while $u:[0,T]\rightarrow \R^d$ represents the field intensity and can be interpreted as a control.
The parameters $\kappa\in\R$ and $m\in\N^*$ determine a polynomial nonlinearity, which arises naturally from mean-field approximations of many-body interactions and reproduces phenomena such as self-focusing, defocusing, and solitary waves.

When the bilinear forcing term in \eqref{NLS} is interpreted as a control acting on the system, the analysis is considerably delicate. Although bilinear controls are physically natural in quantum mechanics, they are typically too weak to ensure exact controllability in infinite-dimensional settings. In particular, when the potential $\mu$ is too \emph{regular}, one encounters intrinsic obstructions to controllability in the underlying Hilbert space, as first observed in \cite{Ball}.

This motivates the use of \emph{singular} potential fields to achieve exact controllability while preserving the well-posedness properties. In this work, we pursue this strategy by considering a class of singular potentials for which the associated nonlinear Schrödinger dynamics remain well-posed in low regularity Sobolev spaces. At the same time, we show that the increased “strength” induced by the singular structure allows to recover global exact controllability in the spaces. 

In detail, the aim of this work is twofold.
\begin{itemize}[leftmargin=2em, rightmargin=2em]
    \item First, we prove the well-posedness of the periodic nonlinear Schrödinger equation \eqref{NLS}
    in low-regularity Sobolev spaces $H^s(\T;\C)$  for suitable $s>0$ and for certain  \emph{singular} potentials $\mu_j$ rough enough so that for $\psi\in H^s(\T;\C)$, $\mu_j\psi\not\in H^s(\T;\C)$.
    
    \medskip
    
    \item Second, we establish small-time global exact controllability of \eqref{NLS} for a class of states in $H^s(\T;\C)$ and small-time exact controllability between eigenmodes of the Laplacian.
\end{itemize}

\subsection{Main results}

%\subsubsection{Some notations}
Consider the one-dimensional nonlinear Schr\"odinger equation \eqref{NLS} on the one-dimensional torus $\T$ which we identify with the interval $(-\pi,\pi)$. %Theorem \ref{THM:control_global_intro} is consequence of few key results and the first is the following local exact controllability between eigenfunctions for the \eqref{NLS} in the space $H^s(\T)$. The results is ensured by 
Recall that in $L^2(\T;\C)$ the unbounded operator $-\dx^2$ with domain $H^2(\T;\C)$ is closed, self-adjoint and has compact resolvent. It also admits a family of eigenfunctions $(e_k)_{k\in\Z}$ forming an Hilbert basis for the space $L^2(\T;\C)$.  We denote $(\lambda_k)_{k\in\Z}$ the corresponding eigenvalues and we assume
$$\lambda_k=k^2,\ \ \ \ \ \ \ e_k(x)=\frac1{\sqrt{2\pi}}e^{ikx},\ \ \ \ \ \ k\in\Z,\ \  x\in\T.$$
We consider the Fourier series of $\psi\in\D'(\T;\C)$ defined by
$$\ft\psi_k = \langle \psi,e_k\rangle_{\D'\times\D},~~k\in\Z,$$
and we have in the sense of distributions
$$\psi =\sum_{k\in\Z}\ft\psi_k e_k.$$
Let $\jb{\nabla}^s$ denote the Fourier multiplier with symbol $\jb{k}^s$, where $\jb{k} = \sqrt{1+k^2}$.

\subsubsection{Well-posedness in low regularity Sobolev spaces with singular potentials}

\begin{Assumptions}\label{AssumptionLWP}
Assume that for some $d\in\N^*$, $u=(u_1,\dots,u_d)$ and $\mu=(\mu_1,\dots,\mu_d)$ satisfy $u\in L^2(\R;\R^d)$ and that the potentials $\mu_j\in\D'(\T;\R)$ are in $\F\ell^{s,\infty}(\T;\R)$, i.e. satisfy
\begin{equation}\label{condmu}\tag{$\mathbf{H^{\text{\tiny (s)}}_{\text{\tiny WP}}}$}
~\exists C>0,~\forall j = 1,\dots,d,~\forall k \in \Z,~|\ft{(\mu_j)}_k|\leq C\jb{k}^{-s},
\end{equation}
with finite norm $\|\mu_j\|_{\F\ell^{s,\infty}}=\sup_{k\in\Z}\jb{k}^s|\ft{(\mu_j)}_k|$.
\end{Assumptions}

\begin{theorem}\label{THM:LWP_intro}
Let $m\in\N^*$. For any $\kappa\in\R$, then the equation \eqref{NLS} is locally well-posed in $H^s(\T;\C)$ provided that:

\smallskip

\begin{enumerate}[label=\textup{(\roman*)}, leftmargin=2em, rightmargin=2em]

 \item (linear case) $s>0$ if $\kappa=0$, and that \eqref{condmu} holds. In this case the solution is global in time;

 \item (nonlinear case) $s>s_\scal(m):=\max\big\{\tfrac12-\tfrac1m,0\big\}$ if $\kappa\neq 0$, and that \eqref{condmu} holds. In this case, if $\psi_\reff$ is the global solution to \eqref{NLS} with $u\equiv 0$ starting from $e_0$, then for any $T>0$, provided that $u$ and $\psi(0)-\psi_\reff(0)$ are small enough, the solution $\psi$ extends to $[0,T]$.

\item In both cases, the mass $t\mapsto \|\psi(t)\|_{L^2}^2$ is conserved, and the flow map $\Phi : (u,\mu,\psi_0)\to \psi$ is smooth. 

\end{enumerate}
\end{theorem}
We refer to Section~\ref{SEC:bilinear} and Theorem~\ref{THM:LWP} below for a more detailed version of this result. In particular, Remark~\ref{RM:uniqueness} below discusses our notion of mild solution as well as the conditional or unconditional nature of our uniqueness statement.

\medskip
Let us point out that Theorem~\ref{THM:LWP_intro} provides an extension to lower regularities (also valid in the nonlinear case) of the well-posedness results of \cite{BL}\footnote{These are concerned with Dirichlet or Neumann boundary conditions on the interval $(0,\pi)$, instead of periodic boundardy conditions on $(0,2\pi)$ as in this work. At low regularity, the former setting can be embedded into the latter by restricting the dynamics to odd or even functions, which are invariant subspaces provided that $V(t)$ also satisfies a parity or imparity condition.} and of \cite[Theorem 1.2]{Alessandro}. Moreover, Theorem~\ref{THM:LWP_intro} also applies to larger non-integer values of the regularity $s$ (compared to $s\geq 3$ odd in the Dirichlet case and $s\geq 2$ even in the Neumann case of \cite{BL}, or $s=1$ in \cite{Alessandro}), which could not be treated using the techniques employed in the mentioned works. 

\medskip
Regarding the range of regularity $s$ where we prove well-posedness, we note that when $m=1,2$, Theorem~\ref{THM:LWP_intro} barely misses the case $s=0$ (well-posedness in $L^2(\T;\C)$) even though $s_\scal(m)=0$. Note that when $u\equiv 0$, well-posedness of the cubic (corresponding to $m=1$) NLS in $L^2(\T;\C)$ is known since the pioneering work of Bourgain \cite{Bourgain}. However, the condition \eqref{condmu} on $\mu$, motivated by the scaling for proving the exact controllability in $H^s$ of the linearized equation around $(u,\psi_0) = (0,e_n)$ (see Proposition~\ref{THM:control-periodic}), leads to a loss of regularity in the equation. More precisely, there is an arbitrarily small loss of regularity in the crucial endpoint $L^\infty_xL^4_t$ type estimate of Proposition~\ref{PROP:bilinear} below (see Lemma~\ref{LEM:counting}). This is not only technical but actually prevents the use of a fixed-point argument to show well-posedness in $L^2$, and therefore to have smoothness of the flow map as in Theorem~\ref{THM:LWP_intro}~(v) above, as shown by the last author with Thomann \cite{RT}. Note that even when $\kappa=0$ (linear case), the equation \eqref{NLS} is $L^2$-critical under \eqref{condmu}; see the discussion in Subsection~\ref{SUBS:WP} below. The mild ill-posedness proved in \cite{RT} is therefore similar to the situation encountered for the $L^2$-critical periodic quintic \eqref{NLS} (corresponding to $m=2$ and $u\equiv 0$) in \cite{Bourgain}, for which the loss of regularity in the $L^6_{t,x}$ periodic Strichartz estimate also prevents local well-posedness in $L^2$ via a fixed-point argument, and where counterexamples to smoothness of the flow map is showed in \cite{Kishimoto}.

\subsubsection{Small-time exact controllability}

\begin{Assumptions}\label{AssumptionControl}
Assume that $d=Q^{\mathsf D}+Q^{\mathsf N}+3$ for some $Q^{\mathsf D},Q^{\mathsf N}\in\N^*$, and that $$\mu = (\mu_1^{\mathsf D},...,\mu_{Q^{\mathsf D}}^{\mathsf D},\mu_1^{\mathsf N},,...,\mu_{Q^{\mathsf N}}^{\mathsf N}, 1, \sin(x), \cos(x)),$$ where $\mu_1^{\mathsf D},\dots,\mu_{Q^{\mathsf N}}^{\mathsf N}\in \D'(\T,\R)$ are real-valued periodic distributions satisfying
\begin{equation}\label{condmu2}\tag{$\mathbf{H^{\text{\tiny (s)}}_{\text{\tiny C}}}$}
\begin{split}
\exists C>0,~\begin{cases}
\mu_j^\Dir \text{ is odd, and } \forall k\in\N^*,~~|(\ft{\mu^{\mathsf D}_j})_k|\geq C\jb{k}^{-s},\ \ \ \ \ \ \  &\forall j\in\{1,\dots, Q^\Dir \},\\
\mu_l^\Neu  \text{ is even, and } \forall k\in\N,~~|(\ft{\mu^{\mathsf N}_l})_k|\geq C\jb{k}^{-s},\ \ \ \ \ \ \  &\forall l\in\{1,\dots, Q^\Neu \}.
\end{cases}\end{split}
\end{equation}
\end{Assumptions}

\smallskip

\begin{theorem}\label{THM:control_global_intro}
Let $m\in\N^*$ and $r>0$. Assume \begin{equation}\label{resonances}\kappa\not\in\Big\{\frac{2^{m-1}\pi^m}{m\, r^{2m}}(j^2+k^2):\ j,k\in\N^*,\ j\neq k\Big\}.\end{equation} Let $s_1>\tfrac{1}{2}$. Consider \eqref{NLS} in the presence of $d=5$ controls $$\mu = (\mu^{\mathsf D},\mu^{\mathsf N}, 1, \cos(x),\sin(x))$$  with $\mu^{\mathsf D}$ and $\mu^{\mathsf N}$ verifying the conditions \eqref{condmu}-\eqref{condmu2} for some $s\in (0,\min(s_1,1))$ such that $s>s_\scal(m):=\max\big\{\tfrac12-\tfrac1m,0\big\}$ if $\kappa\neq 0$.
Then, for any $T>0$ and $$\psi_1,\psi_2\in\{\psi= \rho e^{i\theta},\ \rho \in H^{s_1}(\T;\R_+^*),\ \theta\in H^{s_1}(\T;\R),\ \|\rho\|_{L^2}=r\},$$  there exists $u\in L^2((0,T);\R^{d})$ such that the  unique mild solution $\psi$ to the equation \eqref{NLS} (in the sense of Theorem \ref{THM:LWP_intro}) starting from $\psi(0)=\psi_1$ satisfies
$$\psi(T)=\psi_2.$$
\end{theorem}

Theorem \ref{THM:control_global_intro} allows us to identify a low-regularity Sobolev space that can be exactly controlled in arbitrarily small time. To the best of our knowledge, this is the first results establishing global exact controllability in arbitrarily small time for \eqref{NLS}, novel not only in the nonlinear case but also in the linear setting.   

\begin{remark}\label{rmk-res}\rm~
The non-resonance hypothesis \eqref{resonances} is required in the proof to identify a real control ensuring the local controllability which is part of the proof. These resonances happen when
$$\frac{m\, r^{2m}\kappa}{2^{m-1}\pi^m}\in\{j^2+k^2:\ j,k\in\mathbb{N}^*,\ j\neq k\}=\{5,10,13,17,20,25,\dots\}.$$
Note that one can still ensure controllability also in this case by including two more controls $\widetilde \mu^{\mathsf D}$ and $\widetilde \mu^{\mathsf N},$ verifying the conditions \eqref{condmu}-\eqref{condmu2}.

\end{remark}

\smallskip

\begin{remark}\rm~
The potentials $\mu_1^{\mathsf D},\dots,\mu_{Q^{\mathsf N}}^{\mathsf N}$ in Assumption~\ref{AssumptionControl} and $1$ if necessary, are used to prove local exact controllability around eigenmodes, in the spirit of \cite{BL,Alessandro}, while the smooth potentials $1,\sin(x),\cos(x)$ are used to prove global approximate controllability between suitable states, in the spirit of \cite{Vahagn2,BP04-2,BCP}. Note that our motivation for the well-posedness result of Theorem~\ref{THM:LWP_intro} above under Assumption~\ref{AssumptionLWP} is to have a well-defined dynamics for potentials satisfying the controllability condition \eqref{condmu2}. Indeed, \eqref{condmu} and \eqref{condmu2} are obviously compatible, and the class of potentials $\mu_1^{\mathsf D},\dots,\mu_{Q^{\mathsf N}}^{\mathsf N}$ admissible for the controllability is clearly generic in the class of potentials admissible for the well-posedness; compare with Remark 2 in \cite{BL}.
\end{remark}

\smallskip
 
\begin{remark}\rm
The controllability result in Theorem~\ref{THM:control_global_intro} allows to stabilize the usual ($u\equiv 0$) nonlinear Schrödinger dynamics by chosing carefully the control $u$ depending on the problematic initial data. For example, in the mass-critical focusing case $r=1$, $m=2$, $\kappa=1$, without potential ($u\equiv 0$), the construction of Ogawa-Tsutsumi \cite{OT} provides smooth, real-valued and positive initial data $\psi_1$ with mass at or above the mass-threshold of the ground state, and which lead to a self-similar type finite-time blow-up. However, in view of Theorem~\ref{THM:control_global_intro}, the controls $u_j$ can be chosen so as to perturb the solution starting from the same data above the mass-threshold in order for it to merge with the constant in space global trajectory $\|\psi_1\|_{L^2}e^{-it\kappa\|\psi_1\|_{L^2}^{2m}}$ before the blow-up time.
\end{remark}

We now observe that the small-time global exact controllability \eqref{NLS} can also be achieved between eigenmodes. This case is not covered by Theorem~\ref{THM:control_global_intro} and is stated in the following theorem.

\begin{theorem}\label{THM:control_global_intro_eigen}
Let $m\in\N^*$. Assume \begin{equation*}\kappa\not\in\Big\{\frac{2^{m-1}\pi^m}{m}(j^2+k^2):\ j,k\in\N^*,\ j\neq k\Big\}.\end{equation*} 
Consider the nonlinear Schr\"odinger equation \eqref{NLS} in the presence of $d=5$ controls $$\mu = (\mu^{\mathsf D},\mu^{\mathsf N}, 1, \cos(x),\sin(x))$$  with $\mu^{\mathsf D}$ and $\mu^{\mathsf N}$ verifying the conditions \eqref{condmu}-\eqref{condmu2} for some $s\in (0,1/2)$ such that $s>s_\scal(m):=\max\big\{\tfrac12-\tfrac1m,0\big\}$ if $\kappa\neq 0$. For any $T>0$ there exists $u\in L^2((0,T);\R^{5})$ such that there exists $\psi$ a unique mild solution to the equation \eqref{NLS} (in the sense of Theorem \ref{THM:LWP_intro}), satisfying
$$\psi(0)=e_{n_1},\qquad  \psi(T)=e_{n_2}.$$
\end{theorem}

\begin{remark}
\rm
To the best of our knowledge, Theorem \ref{THM:control_global_intro_eigen} is the first result of small-time exact controllability for the \eqref{NLS} between eigenmodes. It ensures that for suitable nonlineairites all eigenmodes can be exactly controlled in arbitrarily small time. Note that the result can also be obtained for general $\kappa \in \R$ (as for Theorem \ref{THM:control_global_intro_eigen} as explained in Remark \ref{rmk-res}) by adding 2 more controls to the equation.
\end{remark}

\subsection{State of the art on bilinear control}

Let $X$ be a Banach space, $A : D(A) \to X$ a generator of a $C^0$-semigroup, $B : X \to X$ a bounded operator, and $u : [0, +\infty) \to \mathbb{R}$. We consider the system
\begin{equation}\label{bil}
i\partial_t \psi = A\psi + u(t)B\psi.
\end{equation}

Controllability of bilinear systems such as \eqref{bil} is delicate due to the multiplicative nature of the control. A fundamental obstruction was identified in~\cite{Ball}, where it is shown that the reachable set is contained in a countable union of compact subsets of $X$. In infinite dimensions, this implies that the reachable set has empty interior and thus excludes exact controllability.

This obstruction to exact controllability for bilinear systems applies in particular to \eqref{NLS} with $\kappa = 0$ in $L^2$, as proved in~\cite{Turinici}, and extends to Sobolev spaces and to unbounded controls in~\cite{UP}, and to \eqref{NLS} with $\kappa\neq 0$ at low regularity in \cite{RT}. One must therefore consider weaker notions of controllability, such as approximate controllability or exact controllability in higher-regularity spaces not preserved by the operator $B$.

\subsubsection{The controllability issue}

Approximate controllability for \eqref{NLS} when $\kappa=0$ has been extensively studied in the literature in large time. Global results were proved via Lyapunov methods in~\cite{Mir09, Ner10}, adiabatic arguments in~\cite{BCMS12, BGRS15}, and Lie--Galerkin techniques in~\cite{BdCC13, BCMS09}. Nevertheless until recently, it remained unclear whether approximate controllability could hold in small time.

A breakthrough on this matter was initiated in~\cite{Vahagn}, where small-time global approximate controllability between eigenmodes was established, opening the way to several subsequent works~\cite{small-time-molecule, coron-2023, small-time-momentum, small-time-wave}. This breakthrough was ultimately concretized in the works~\cite{BP04-1,BP04-2}, which provided the first results of small-time global approximate controllability in $L^2$ for \eqref{NLS} with $\kappa=0$ on $\R^d$ and $\T^d$. Nevertheless, to the best of our knowledge, the extension of these results to higher regularity spaces $H^s$ and to polynomial nonlinearities remains an open problem (we refer to~\cite{BCP} for related results in the logarithmic nonlinear setting).

Concerning the exact controllability of \eqref{NLS}, a major breakthrough was achieved in~\cite{Bea05} and later extended in~\cite{BL}. In \cite{BL}, the authors study \eqref{NLS} on an interval, driven by the Dirichlet Laplacian $\Delta_{\mathsf D}$ (resp. the Neumann Laplacian $\Delta_{\mathsf N}$) in the case $\kappa = 0$ (resp. $\kappa=-1$). The control potential $\mu$ is taken in $H^3(0,1)$ (resp. in $H^2(0,1)$) without imposing any compatibility conditions at the boundary; as a consequence, $\mu$ exhibits a singular behavior at the boundary.

In this setting, the multiplication operator with respect to $\mu$ is continuous on $D(\Delta_{\mathsf D})$ (resp. in $H^1(0,1)$). Thus, the Ball--Marsden--Slemrod obstruction prevents the exact controllability in such space. Nevertheless, this obstruction disappears in higher regularity spaces, in particular in $D(|\Delta_{\mathsf D}|^{3/2})$ (resp. in $D(\Delta_{\mathsf N})$). In \cite{BL}, well-posedness is first established and then the authors prove the local exact controllability around the ground state under the assumption
\begin{gather}\label{cond_kar}\forall k\in\N^*,\; \big|\langle \sin(k\pi x),\mu\sin(\pi x)\rangle_{L^2}\big|\geq  {C}{k^{-3}},\\
\text{\Big(resp. $\forall k\in\N,\; \big|\langle \cos(k\pi x),\mu\rangle_{L^2}\big|\geq  {C}\,{\min\{1,k^{-2}\}}$\Big)}.\nonumber
\end{gather}

This approach has been extended to simultaneous control problems in~\cite{Duc20b, Mor14, MN15}, and to nonlinear equations in~\cite{Karine2, Vahagn2}. More recently, further obstructions in the Neumann case were identified in~\cite{Karine3}. Despite these developments, the method of~\cite{BL} has remained essentially restricted to Dirichlet settings in $D(|\Delta_{\mathsf D}|^{3/2})$  (resp. Neumann settings in $D(\Delta_{\mathsf N})$) or 
high regularity spaces.

A recent development was achieved in~\cite{Alessandro}, where the local exact controllability is proved in $H^1_0(0,1)$ for \eqref{NLS} with $\kappa=0$ in the Dirichlet case. The key idea is to replace boundary singularities in $\mu$ by interior ones via piecewise constant potentials. This choice of $\mu$ bypasses the Ball--Marsden--Slemrod obstruction in $H^1_0(0,1)$ and yields the controllability under
\begin{align}\label{cond_Ale}\forall k\in\N^*,\; \big|\langle \sin(k\pi x),\mu\sin(\pi x)\rangle_{L^2}\big|\geq  {C}{k^{-1}}.\end{align}

In~\cite{Alessandro}, the periodic boundary conditions are also treated, and the main difficulty is the presence of spectral resonances for the problem. The authors overcome this issue via a Zeeman-type splitting induced by an additional magnetic field, restoring simplicity of the spectrum and extending the Dirichlet strategy to $H^1(\T;\C)$.
 
Despite these advances, the existing results were only sufficient to establish global exact controllability for \eqref{NLS} in large time, and only in the linear setting. This could be achieved by combining global approximate controllability results in large time in sufficiently regular Sobolev spaces with local exact controllability in the same topology. 

Nevertheless, it remained unclear whether small-time global exact controllability, as in Theorem~\ref{THM:control_global_intro}, or small-time exact controllability between eigenmodes, as in Theorem~\ref{THM:control_global_intro_eigen}, could be achieved. In this work, we establish these properties also in the nonlinear setting. 

\subsubsection{Well-posedness of the free evolution}\label{SUBS:WP}

When $u\equiv 0$, the periodic nonlinear Schrödinger equation \eqref{NLS} has been extensively studied in the literature. In the Eulidean setting, the equation enjoys both scaling and Galilean invariance, namely
$$\psi(t,x)\mapsto \psi_\lambda(t,x):=\lambda^{-\frac1m}\psi(\lambda^{-2}t,\lambda^{-1}x),\qquad \lambda>0,$$
and
\begin{equation*}
\psi(t,x)\mapsto \psi_{v}(t,x):=e^{i(xv+t\frac{v^2}{2})}\psi(t,x-vt),\qquad v\in\R.
\end{equation*}
These invariances suggest that, also for the periodic problem, well-posedness should naturally hold in Sobolev spaces at regularities no lower than the threshold
\begin{equation*}
s \geq s_{\scal}(m):=\max\big(\frac12-\frac1m,0\big),
\end{equation*}
where each of the critical regularities is linked to the corresponding invariance.

The study of low-regularity well-posedness of the periodic problem was initiated by Bourgain \cite{Bourgain}, and subsequently developped in a series of works exploiting multilinear estimates, periodic Strichartz-type bounds, and refined Fourier restriction methods; see for instance \cite{BD,HTT,KV,Wang}. As a consequence of these contributions, local well-posedness for \eqref{NLS} with $u\equiv 0$ is now understood in $H^s(\T;\C)$ in the following regimes:
\begin{align*}
\begin{cases}
s \geq 0, & m=1,\\
s>0, & m=2,\\
s \geq s_{\scal}(m), & m\geq 3.
\end{cases}
\end{align*}
In particular, the mass-critical case $m=2$, corresponding to the scaling threshold $s_{\scal}(2)=0$, remains especially delicate in the periodic setting, and global well-posedness in $L^2(\T;\C)$ is still a challenging open problem; see \cite{Kishimoto}, and \cite{HerrKwak1,HerrKwak2,Guo} for recent progress in this direction.

In the case of a non-trivial external multiplicative interaction given by $u\neq 0$ and $\mu$ of low regularity, the structure of the equation \eqref{NLS} changes, since the dynamics is driven by a bilinear coupling between a time-dependent control of limited regularity and a rough spatial potential. In this regime, a scaling argument suggests a second critical threshold associated with the regularity of the potential. Indeed, the rescaling
$$\big(u,\mu,\psi\big) \mapsto \big(\lambda^{-1}u(\lambda^{-2}t),\lambda^{-1}\mu(\lambda^{-1} x),\lambda^{s-\frac12}\psi(\lambda^{-2}t,\lambda^{-1}x)\big),~~\lambda>0,$$ 
leaves both the equation invariant, as well as the $L^2$ norm of $u$ and the $H^s$ norm of $\psi$ invariant, and thus suggests from the regularity of $\mu$ that 
\begin{equation}\label{scalingmu} s_\mu = 0 \end{equation} is critical for this term. 
This indicates that the regime $s>0$ should be interpreted as subcritical with respect to the multiplicative interaction. In the purely bilinear case $\kappa=0$ or in the cubic case $\kappa\neq 0$, $m=1$, Theorem~\ref{THM:LWP_intro} therefore covers the whole sub-critical regime $s>0=s_\mu$ where it is possible to control the bilinear coupling. For higher values of $m$, we have $s_{\scal}(m)= \frac12-\frac1m\geq 0=s_\mu$, and therefore the bilinear interaction becomes sub-critical compared to the nonlinear contribution.

\subsection{The developed techniques and scheme of the work}

As mentioned above, bilinear control systems associated with \eqref{NLS} are in general weak, since, as shown in~\cite{Ball}, even regular potentials already induce obstructions to the exact controllability.
For this reason, Theorems~\ref{THM:control_global_intro} and~\ref{THM:control_global_intro_eigen} rely on singular potentials and the well-posedness result of Theorem~\ref{THM:LWP_intro} (proved in Section~\ref{Sec:WP}) is a key ingredient of the analysis. 

The proofs of Theorems~\ref{THM:control_global_intro} and~\ref{THM:control_global_intro_eigen} (Section~\ref{Sec:global_proof}) combine the following results. 

\begin{enumerate}[label=\textup{(\roman*)}, leftmargin=2em, rightmargin=2em]
\item Small-time global approximate controllability for states in $H^{s_1}(\T)$ with respect to the $H^s$ topology ($s_1=1$ in the proof of Theorem~\ref{THM:control_global_intro_eigen}). The result is ensured by Proposition~\ref{THM:control_app_1} in Section~\ref{Sec:global-approx}.

\item Local exact controllability in $H^s(\T;\C)$ around the ground state proved in Proposition~\ref{THM:contr_ground} in Section~\ref{Sec:control_ground_proof}.
\end{enumerate}

A central point is the alignment of the functional settings. Small-time global approximate controllability has been established in the literature in $L^2$ in~\cite{BCP, BP04-1, BP04-2}, while local exact controllability has naturally been formulated in higher-regularity Sobolev spaces~\cite{BL, Alessandro}. Here, we show that approximate controllability also extends to low-regularity Sobolev spaces, while exact controllability can be achieved under weaker regularity assumptions. This makes it possible for both results to hold within a common functional framework, and thus to be coherently combined into a unified dynamical setting.

\subsubsection{Global approximate controllability: the techniques}

As mentioned before, the recent work~\cite{Vahagn} introduced a saturation method yielding small-time approximate controllability of phases for \eqref{NLS} on the $N$-torus, implying controllability between eigenmodes. This approach was refined in~\cite{BCP, BP04-1, BP04-2}, leading to small-time global approximate controllability in $L^2$.

We extend these techniques to Sobolev spaces $H^s$ by adapting the strategy of~\cite{BP04-2}. The proof of Proposition~\ref{THM:control_app} relies on a decoupling of angular and radial dynamics. The angular component is controlled using~\cite{Vahagn}, while the radial component is approximated in $H^s$ by unitary actions of smooth diffeomorphisms.

These actions are approximated via flows generated by suitable vector fields as \cite{BP04-2} but this time in $H^s$. Unlike the $L^2$ setting, where semigroup methods and the Trotter–Kato theorem apply, the $H^s$ framework requires direct energy estimates that we develop in Section~\ref{Sec:global-approx}.

\subsubsection{Local exact controllability: the techniques}
We lower the regularity required to ensure local exact controllability in Propositions~\ref{THM:contr_ground} compared with the work in~\cite{Alessandro}, which dealt with the $H^1$ framework. To this end, we introduce classes of controls satisfying assumptions \eqref{condmu}–\eqref{condmu2}, inspired by \eqref{cond_kar} in \cite{BL} and \eqref{cond_Ale} in \cite{Alessandro}. This choice allows us to bypass the well-known Ball–Marsden–Slemrod~\cite{Ball} obstruction to the exact controllability.

%These conditions ensure both well-posedness in low-regularity Sobolev spaces and solvability of the associated moment problems.

The local controllability result is reduced to the exact controllability of suitable linearized systems  obtained via moment problems involving resonant frequencies associated with the double eigenvalues of the Laplacian.
To handle these resonances, we introduce independent controls $\mu^{\mathsf D}$ and $\mu^{\mathsf N}$ and decompose the dynamics on $\T\simeq(-\pi,\pi)$ into Dirichlet and Neumann subproblems on $(0,\pi)$. Each component exhibits resonances, requiring a decoupled control structure.

This idea was already suggested in~\cite[Remark 4.3]{Alessandro}, but could not be implemented there due to the lack of a functional framework ensuring both well-posedness and controllability in $H^1(\T;\C)$. In~\cite{Alessandro}, only the Dirichlet component is treated successfully, while the Neumann case cannot be handled within the same framework.

Thanks to Theorem~\ref{THM:LWP_intro}, we overcome this limitation and construct controls simultaneously acting on both components. For the Neumann part, one may for instance take $\mu^{\mathsf N} = D^{-s}\delta$.
Note that when $s=1$, these controls are not piecewise $H^1$ functions, and the obstruction of~\cite[Proposition 4.1]{Alessandro} does not apply, which is crucial to close the argument.

\subsubsection{Well-posedness: the techniques}

The well-posedness of \eqref{NLS} at low regularity is delicate due to potential derivative loss coming from the roughness of the potential, as $\F\ell^{s,\infty}$ scales like $H^{s-\frac12}$, leading to a potential $\frac12$ derivative loss. Actually, this also scales like $W^{s,1}$, which shows instead an integrability loss to control the product $\mu\psi$ in $H^s$. Therefore the proof of Theorem~\ref{THM:LWP_intro} uses harmonic analysis techniques and Bourgain-type $X^{s,b}$ spaces to compensate for this loss, by relying strongly on a Fourier-side bilinear refinement of an endpoint $L^\infty_xL^4_t$ Strichartz-type estimate; see Lemma~\ref{LEM:bilinear} and the discussion in Subsection~\ref{SEC:bilinear} below.

This framework extends previous results, which either required high regularity in the presence of singular potentials or excluded such perturbations in low-regularity regimes.

\medskip 

\noindent
{\bf Acknowledgments} The first two authors were partially supported by the French National Research Agency (ANR), Project QuBiCCS (ANR-24-
CE40-3008). R.B was partially supported by the French National Research Agency (ANR), Project COINS (ANR-26-CE40-7933-01). T.R. was partially supported by the French National Research Agency (ANR), Project Smooth (ANR-22-CE40-0017).

\medskip 

\noindent
{\bf AI disclaimer}. The authors used AI tools to improve the presentation of this article. All ideas and proofs are due to the authors, who take full responsibility for the content of this article.

\section{Preliminaries}
\subsection{Notations and counting estimates}

For complex numbers $a\in\C$, we write $a^-=\cj a$ and $a^+=a$. For positive numbers $A$ and $B$, we write $A\les B$ if there is a constant $C>0$, independent of the various parameters, such that $A\le CB$; $C=10^{100}$ should work throughout this text. We also write $A\sim B$ if both $A\les B$ and $B\les A$. We also write 
$$A\vee B = \max\{A,B\},\ \ \ \ A\wedge B=\min\{A,B\}.$$

In the proof of bilinear estimates, we will make use of the following basic counting estimates; see e.g. \cite[Lemmata 5.2~\&~5.4]{KP1} for the first two, and \cite{ET} for the last one.
\begin{lemma}\label{LEM:basic}
The following estimates hold.
\begin{itemize}[leftmargin=2em, rightmargin=2em]
    \item[\textup{(i)}] Let $\Lambda \subset \Z^2$. We assume that the projection of $\Lambda$ on the $x$ axis is contained in an interval $I\subset \R$. Moreover, we assume that the measure of the $k_2$-sections of $\Lambda$ (that is the sets $\left\{k_2 \in \Z, (k_1,k_2)\in \Lambda\right\}$ for a fixed $k_1$) is uniformly (in $k_1$) bounded by a constant $C$. Then we have
\begin{equation}\label{mesure section}
\#\Lambda \leq C |I|.
\end{equation}
% \item[\textup{(ii)}] Let $I$, $J$ be two intervals in $\R$, and let $f : I \rightarrow \R$ be a $\mathcal{C}^1$ function with $\inf_{\xi\in I}\left|f'(\xi)\right|>0$. Then
%\begin{equation}\label{estimation mesure cas Z}
%\#\left\{k \in I\cap \Z,~f(k) \in J\right\} \lesssim 1+\frac{|J|}{\inf_{\xi\in I}|f'(\xi)|}.
%\end{equation}
 \item[\textup{(ii)}] Let $a\neq 0$ ,$b$, $c$ be real numbers and $J\subset \R$ a bounded interval. Then
\begin{equation}\label{estimation mesure parabole Z}
\#\left\{k \in \Z,~ak^2+bk+c \in J\right\} \lesssim 1+\frac{|J|^{1/2}}{|a|^{1/2}}.
\end{equation}
\item[\textup{(iii)}] For any $\epsilon>0$, there is $C>0$ such that for any $m\in\Z^*$, 
\begin{equation}\label{estimation diviseurs}
\#\left\{k \in \Z,~k|m\right\}\leq C|m|^\epsilon.
\end{equation}
\end{itemize}
\end{lemma}

We will also make crucial use of a refined version of the classical divisor estimate \eqref{estimation diviseurs} due to Deng-Nahmod-Yue \cite[Lemma 3.4 (1)]{DNY}.
\begin{lemma}[\cite{DNY}]\label{LEM:DNY}
For any $\epsilon>0$, there is $C>0$ such that for any $m,q\in\Z^*$ and $\rho>0$ such that $|q|\geq |m|^\epsilon>0$, it holds 
\begin{equation}\label{estimation diviseurs DNY}
\#\left\{k \in \Z,~k|m,~|k-q|\leq \rho\right\}\leq C\rho^\epsilon.
\end{equation}
\end{lemma}

In turn, the previous estimates will be used to establish the following counting estimates, used below in the proof of bilinear periodic Strichartz type estimates. 
\begin{lemma}\label{LEM:counting}
The following estimates hold.
\begin{itemize}[leftmargin=2em, rightmargin=2em]
\item[\textup{(i)}] Let $J\subset \R$ be an interval of size $|J|\gtrsim 1$, and $N_1,N_2$ be dyadic integers such that $N_1\not\sim N_2$. Then for any $\epsilon>0$, 
$$\#\big\{(k_1,k_2)\in\Z^2,~|k_1|\sim N_1,~|k_2|\sim N_2,~k_1^2-k_2^2\in J\big\}\les \min\big(|J|(N_1\wedge N_2)^\epsilon,|J|^\frac12(N_1\wedge N_2)\big).$$
\item[\textup{(ii)}] Let $J\subset \R$ be an interval of size $|J|\gtrsim 1$, and $N$ be a dyadic integer. Then for any $\epsilon>0$, 
$$\#\{(k_1,k_2)\in\Z^2,~|k_1|\sim N,~k_1\neq \pm k_2,~k_1^2-k_2^2\in J\}\les \min\big(|J|N^\epsilon,|J|^\frac12N\big).$$
\item[\textup{(iii)}] Let $J\subset \R$ be an interval of size $|J|\gtrsim 1$, and  $N,N_1$ be a dyadic integers satisfying $N\ll N_1$. Then for any intervals $I_1,I_2\subset \R$ of size $|I_1|+|I_2|\les N$, and any $0<\epsilon<\frac12$,
$$\#\{(k_1,k_2)\in\Z^2,~k_1\in I_1,~k_2\in I_2,~|k_1|\sim|k_2|\sim N_1,~0<|k_1-k_2|\les N,~k_1^2-k_2^2\in J\}\les \min\big(|J|N^\epsilon, N^2).$$
\end{itemize}

\end{lemma}

\begin{proof}
We start with the first bound of estimate (i). We have
\begin{align*}
&\#\big\{(k_1,k_2)\in\Z^2,~|k_1|\sim N_1,~|k_2|\sim N_2,~k_1^2-k_2^2\in J\big\}\\
&\qquad\qquad\les |J|\sup_{m\in\Z^*}\#\big\{(k_1,k_2)\in\Z^2,~|k_1|\sim N_1,~|k_2|\sim N_2,~k_1^2- k_2^2=m\big\}.
\end{align*}
By symmetry it suffices to consider the case $N_1\gg N_2$. Then $|m|=|k_1^2-k_2^2|\sim N_1^2\neq 0$ and $|k_1-k_2|\sim|k_1+k_2|\sim N_1$. Following e.g. \cite[Lemma 3.2]{BGT}, either $N_1\les N_2^2$, then using the divisor bound \eqref{estimation diviseurs} (which is possible since $m\neq 0$) we get for arbitrary $\delta>0$
\begin{align*}
&\sup_{m\in\Z}\#\big\{(k_1,k_2)\in\Z^2,~|k_1|\sim N_1,~|k_2|\sim N_2,~k_1^2- k_2^2=m\big\}\\
&\les \sup_{|m|\sim N_1^2}\#\big\{(n_1,n_2)\in\Z^2,~|n_1|\sim |n_2|\sim N_1,~n_1n_2=m\big\}\\
&\les N_1^{2\delta}\les N_2^{4\delta},
\end{align*}
which, by taking $\delta\ll\epsilon$, is enough for the first bound in the right-hand side of (i). Or else $N_1\gg N_2^2$, but then observe that the constraints $k_1^2-k_2^2 = m$ with $\frac12 N_j\leq |k_j|\leq 2N_j$ imply that $N_1^2\leq m\leq 10N_1^2$. Moreover, $k_1$  takes value in the union $[-\sqrt{m+4N_2^2},-\sqrt{m+\frac14 N_2^2}]\cup[\sqrt{m+\frac14 N_2^2},\sqrt{m+4N_2^2}]$, where each interval is of size $O(N_1^{-1}N_2^2)\ll 1$ since $N_1 \gg N_2^2$. Therefore the number of possible $k_1$ is bounded uniformly in $N_1,N_2$, while clearly for each such $k_1$ there are also at most 2 values of $k_2^2=k_1^2-m$. Therefore we get
\begin{align*}
\sup_{m\neq 0}\#\{(k_1,k_2)\in\Z^2,~|k_j|\sim N_j,~k_1^2-k_2^2=m\} \les 1
\end{align*}
in this case. All in all, this shows the first bound in (i).

As for the second bound, by symmetry of $k_1$ and $k_2$ we can always assume again $N_1\geq N_2$. Then using Lemma~\ref{LEM:basic}, more precisely using together \eqref{mesure section} with $k_2$ lying in the interval $[\frac12 N_2,2N_2]$, and \eqref{estimation mesure parabole Z} to measure the size of the $k_1$ sections with $k_1^2-k_2^2=ak_1^2 +b k_1+c$ for fixed $k_2$, we get
\begin{align*}
\#\big\{(k_1,k_2)\in\Z^2,~|k_1|\sim N_1,~|k_2|\sim N_2,~k_1^2-k_2^2\in J\big\}\les (N_1\wedge N_2)|J|^\frac12.
\end{align*}
This gives the second bound in (i).

The proof of (ii) is similar: now we need to estimate
$\sup_{m\in\Z^*}\#\big\{(k_1,k_2)\in\Z^2,~|k_1|\sim N,~k_1^2-k_2^2=m\big\}$. Note that the condition $m\neq 0$ is ensured by the assumption $k_1\neq \pm k_2$.
We repeat the previous counting argument but reasoning on $0<|m|\les N^4$ or $|m|\gg N^4$ instead: in the former case, the divisor bound gives again
$$\#\{(k_1,k_2)\in\Z^2,~|k_1|\sim N,~k_1^2-k_2^2=m\}\les |m|^{\frac\epsilon4}\les N^\epsilon.
$$
In the latter case, we must have $|k_2|\gtrsim N^2$ and $m<0$ to ensure $|m|\gg N^4$ given that $|k_1|\sim N$. Therefore, for each fixed $k_1$, $k_2$ takes value in $[-\sqrt{4N^2-m},-\sqrt{\frac14 N_2^2-m}]\cup[\sqrt{\frac14 N^2-m},\sqrt{4N^2-m}]$ which is of size $O(|m|^{-\frac12}N^2)\ll 1$ since $|m|\gg N^4$. This ensures as before that the number of points is bounded by 1 is this case.

The other required bound
$$\#\big\{(k_1,k_2)\in\Z^2,~|k_1|\sim N,~k_1^2-k_2^2\in J\big\}\les N|J|^\frac12$$
is left unchanged by measuring the $k_2$-sections as before (with the role of $k_1$ and $k_2$ exchanged). This proves (ii).

At last, we turn to the proof of (iii). Let $I_1=[k_1^*-cN,k_1^*+cN]$ and $I_2=[k_2^*-cN,k_2^*+cN]$. Then from the assumptions, $|k_1^*-k_2^*|\les N$ and $|k_1^*|\sim|k_2^*|\sim N_1\gg N$, therefore $|k_1^*+k_2^*|\sim N_1 \geq |m|^\epsilon \sim N^\epsilon N_1^\epsilon$ for $0<\epsilon<\frac12$. We can thus invoke Lemma~\ref{LEM:DNY} with $q=k_1^*+k_2^*$ and $\rho=N$ to bound
\begin{align*}
&\#\{(k_1,k_2)\in\Z^2,~k_1\in I_1,~k_2\in I_2,~|k_1|\sim|k_2|\sim N_1,~|k_1-k_2|\les N,~k_1^2-k_2^2=m\}\\
&\leq \#\{(n_1,n_2)\in\Z^2,~|n_1|\les N,~|n_2-q|\les N,~n_1n_2=m\}\\
&\les N^\epsilon,
\end{align*}
uniformly in $m\neq 0$.

The second required bound is the trivial one
$$\#\big\{(k_1,k_2)\in\Z^2,~k_1\in I_1,~k_2\in I_2,~k_1^2-k_2^2\in J\}\les N^2$$
following from $|I_1|+|I_2|\les N$. This proves (iii).
\end{proof}

\begin{remark}
\rm
Note that in case $N_1\sim N_2$, if $0\in J$, then
\begin{multline*}
\#\{(k_1,k_2)\in\Z^2,~|k_1|\sim N_1,~|k_2|\sim N_2,~k_1^2-k_2^2\in J\}\\\geq \#\{(k_1,k_2)\in\Z^2,~|k_1|\sim N_1,~|k_2|\sim N_2,~k_1^2-k_2^2=0\}\gtrsim N_1,
\end{multline*}
which would lead to a bigger loss of derivative in Lemma~\ref{LEM:bilinear} below, and therefore limit the well-posedness in $H^s$ to the range $s>1$.
\end{remark}

\subsection{Littlewood-Paley theory and function spaces}\label{SUBS:LP}
We now define the phase spaces for the dynamics of \eqref{NLS}. \subsubsection{Fourier-Lebesgue spaces}
Recall that for $s\in\R$ and $p\in[1,\infty]$, the Fourier-Lebesgue space $\F\ell^{s,p}(\T;\C)$ is defined as the space of complex-valued periodic distributions with finite norm
\begin{equation}\label{FLnorm}
\|\psi\|_{\F\ell^{s,p}}:=\|\jb{\nabla}^s\psi\|_{\F\ell^p}=\|\jb{k}^s\ft\psi_k\|_{\ell^p_k},
\end{equation}
where $\jb{k}=\sqrt{1+k^2}$ and $\jb{\nabla}^s$ is the Fourier multiplier with symbol $\jb{k}^s$. In particular, the standard $L^2$-based Sobolev spaces are $H^s(\T)=W^{s,2}(\T)=\F L^{s,2}(\T)$.

\subsubsection{Littlewood-Paley theory}
In the following, we use $N,N_j,M\in 2^{\N}$ to denote dyadic integers in $2^{\N}=\{2^m,~~m\in \N\}$. Given dyadic integers $N_0,\dots,N_k$, we write $N_0^*\geq \dots \geq N_k^*$ for their decreasing rearrangement. 

For $N\in 2^\N$ a dyadic integer and $\chi\in C^\infty_0(\R)$ even with $\chi\equiv 1$ on $[-1,1]$ and $\mathrm{supp}\chi \subset [-\frac54,\frac54]$, we first define $\P_{\le N}$ to be the smoothed out Dirichlet projection onto the frequencies $\{|k|\leq N\}\subset \Z$, i.e. 
\begin{equation}\label{PN}
\P_{\leq N} \psi= \sum_{k\in\Z}\chi(N^{-1}k)\ft \psi_k e_k.
\end{equation}
Then we define the projection on annuli by $\P_N=\P_{\le N}-\P_{\le \frac{N}2}$, with the convention that $\P_{\le \frac12}\psi=\P_0 \psi=\ft \psi_0$ and $\P_0^\perp \psi = \psi-\P_0 \psi$.
More generally, for any subset $S\subset \R$, we denote by $\P_S$ the frequency projection on $S$, i.e.
\begin{equation}\label{PS}
\P_S \psi = \sum_{k\in S\cap \Z}\ft\psi_k e_k.
\end{equation}

%We record the following basic Bernstein type inequalities.
%\begin{lemma}\label{LEM:Bernstein}
%For any $s\in\R$ and $1\leq p\leq q\leq \infty$, there is $C>0$ such that for any dyadic $N$ and $f\in \S(\T;\C)$, it holds
%\begin{align*}
%\|\P_N f\|_{L^p}\le CN^{\frac1q-\frac1p}\|\P_N f\|_{L^q}.
%\end{align*}
%\end{lemma}
%\begin{proof}
%Write 
%\begin{align*}
%\P_N f = \sum_{k\in\Z}\eta(N^{-1}k)\ft f_ke_k = \sum_{k\in\Z}\widetilde{\eta}(N^{-1}k)\widehat{(\P_N f)}_ke_k = \varphi_N * \P_N f,
%\end{align*}
%where $\eta(x)=\chi(x)-\chi(2x)$ as in the definition of $\P_N$, and $\widetilde{\eta}(x)=\chi(\frac{x}{2})-\chi(4x)$ satisfies $\eta\in C^{\infty}_0(\R;\R)$, $\widetilde{\eta}\equiv 1$ on $\supp~\eta$, and ${\displaystyle\varphi_N = \sum_{k\in\Z}\widetilde{\eta}(N^{-1}k)e_k}$. By Poisson summation formula (see e.g. \cite{Grafakos}), we have for all $x\in\T$ that
%\begin{align*}
%\varphi_N(x)=C\sum_{k\in\Z}\widehat{\wt\eta(N^{-1}\cdot)}(x+2\pi k).
%\end{align*}
%Therefore,
%\begin{align*}
%\|\varphi_N\|_{L^1(\T)} &\les \sum_{k\in\Z}\int_{2\pi k}^{2\pi (k+1)}|\widehat{\wt\eta(N^{-1}\cdot)}(y)|dy = \|\widehat{\wt\eta(N^{-1}\cdot)}\|_{L^1(\R)} = \|\widehat{\wt\eta}\|_{L^1(\R)},
%\end{align*}
%and by the triangle inequality,
%\begin{align*}
%\|\varphi_N\|_{L^\infty(\T)}&\les \sum_{|k|\les N}1 \sim  N.
%\end{align*}
%Interpolating using H\"older's inequality gives
%\begin{align*}
%\|\varphi_N\|_{L^p}\les N^{\frac1{p'}},
%\end{align*}
%from which the proof of Lemma~\ref{LEM:Bernstein} follows through Young's convolution inequality.
%\end{proof}

\subsubsection{Solution space}
We will build the solution $\psi$ to \eqref{NLS} in the Bourgain space $X^{s,b}$ adapted to the linear Schr\"odinger equation, which is the Banach space defined through the norm
\begin{equation}\label{Xsb}
\|\psi\|_{X^{s,b}} := \big\|\langle k\rangle^{s}\langle \tau - k^2\rangle^b \hat{\psi}_k(\tau)\big\|_{L^2_\tau\ell^2_k},
\end{equation}
for any $s,b\in\mathbb{R}$, where $\hat{\psi}_k(\tau) = \int_{\mathbb{R}\times\mathbb{T}}\psi(t,x)e^{-i(t\tau+kx)}dxdt$ is the space-time Fourier transform of $\psi : (t,x)\in\mathbb{R}\times\mathbb{T}\mapsto \psi(t,x)\in \mathbb{C}$.

As in Subsection~\ref{SUBS:LP}, we define a Littlewood-Paley characterization of these spaces using smooth projections also in modulation, namely for $\psi\in X^{s,b}$ and a dyadic integer $L$ we define
\begin{equation}\label{QleL}
\Q_{\leq L} \psi= \int_{\R}\sum_{k\in\Z}\chi(L^{-1}(\tau-k^2))\ft \psi_k(\tau)e^{it\tau}d\tau e_k,
\end{equation}
with $\chi$ as in \eqref{PN}, and $\Q_L\psi := \Q_{\leq L}\psi-\Q_{\leq L/2}\psi$, with the convention $\Q_{\leq \frac12} = 0$.
Then we have
\begin{align*}
\|\psi\|_{X^{s,b}}^2\sim \sum_{N,L}\jb{N}^{2s}\jb{L}^{2b}\|\P_N\Q_L\psi\|_{L^2_{t,x}}^2. 
\end{align*}

For $I\subset \R$ a time interval, we also define its time-localized version
\begin{equation}\label{XsbT}
\|\psi\|_{X^{s,b}_I}=\inf\{\|\Psi\|_{X^{s,b}},~~\Psi\equiv \psi\text{ on }I\},
\end{equation}
with the short-hand notation $X^{s,b}_T$ when $I=[0,T]$ for some $T>0$.

Next we collect some linear estimates in the space $X^{s,b}_T$.
\begin{lemma}\label{LEM:linear}
The following properties hold:
\begin{itemize}[leftmargin=2em, rightmargin=2em]
    \item[\textup{(i)}] \textup{($X^{s,b}$ as a resolution space)} Let $s\in\mathbb{R}$ and $b>\frac12$. Then there exists $C>0$ such that for any $\psi\in X^{s,b}$ it holds $\psi\in C(\mathbb{R};H^s(\mathbb{T}))$ and $\|\psi\|_{L^\infty_tH^s}\le C\|\psi\|_{X^{s,b}}$.\\
\item[\textup{(ii)}]\textup{ (Time localization)} Let $s\in\mathbb{R}$ and $-\frac12 <b'\le b <\frac12$. Then there exists $C>0$ such that for any $T\in(0,1]$ and $\psi\in X^{s,b}_T$, it holds $\|\psi\|_{X^{s,b'}_T}\le C T^{b-b'}\|\psi\|_{X^{s,b}_T}$.\\
\item[\textup{(iii)}]\textup{ (Linear estimate)} For any $s\in\mathbb{R}$, $b>\frac12$, there is $C>0$ such that for any $T>0$ and $\psi_0\in H^s(\mathbb{T};\C)$ it holds $$\|e^{-it\dx^2}\psi_0\|_{X^{s,b}_T}\leq C \jb{T}^\frac12\|\psi_0\|_{H^s}.$$
\item[\textup{(iv)}]\textup{ (Energy estimate)} For any $s\in\mathbb{R}$, $b>\frac12$, there is $C>0$ such that for any $T>0$ and $F\in X^{s,b-1}_T$ it holds $$\Big\|\int_0^te^{-i(t-t')\dx^2}F(t')dt'\Big\|_{X^{s,b}_T}\leq C\jb{T}^2\|F\|_{X^{s,b-1}_T}.$$
\end{itemize}
\end{lemma}
These properties are standard, and we refer to \cite{Tao06} and \cite{ET}; see also \cite{BGT_TAMS} for the long-time version of properties (iii) and (iv).

%We also record the following version, due to Bourgain \cite{Bourgain}, of Zygmund's $L^4$ Strichartz type estimate; see also \eqref{Zygmund} below.
%
%\begin{lemma}[\cite{Bourgain,Tao06}]\label{LEM:L4}
%There exists $C>0$ such that for any $\psi\in X^{0,\frac38}$, it holds
%\begin{align*}
%\|\psi\|_{L^4_{t,x}}\leq C\|\psi\|_{X^{0,\frac38}}.
%\end{align*}
%\end{lemma}

\section{Existence and uniqueness of solutions}\label{Sec:WP}

The main goal of this section is to prove the following detailed version of Theorem~\ref{THM:LWP_intro}.
\begin{theorem}\label{THM:LWP}
For any $m\in\N^*$, and any $\kappa\in\R$, then the equation \eqref{NLS} is well-posed in $H^s(\T;\C)$ for suitable $s>0$ in the following sense.

\smallskip

\begin{enumerate}[label=\textup{(\roman*)}, leftmargin=2em, rightmargin=2em]

 \item (Linear case, high regularity). Assume $\kappa=0$, $s>\frac13$ and that Assumption~\ref{AssumptionLWP} holds. For every $\psi_0\in H^s(\T;\C)$, there exists a unique global mild solution $\psi\in C^0(\R; H^s(\T;\C))$  that is a solution to the Duhamel formula:
\begin{align}
    \label{duhamel_lin}\psi(t)=e^{- it\dx^2}\psi_0-i\int_0^te^{ -i(t-t')\dx^2}\big\{(u(t')\cdot\mu)\psi(t')\big\}dt',
\end{align}
for all $t\in \R$, and that satisfies conservation of mass: $\|\psi(t)\|_{L^2}=\|\psi_0\|_{L^2}$ for any $t\in \R$. Moreover, there is $\delta(s)>0$ such that for any $r_1,r_2>0$, there is $\tau=\tau(r_1,r_2)>0$ such that for any $u_j\in L^2(\R;\R)$ and $\mu_j\in\F\ell^{s,\infty}(\T;\R)$ satisfying $\|u_j\|_{L^2}\leq r_1$ and $\|\mu_j\|_{\F\ell^{s,\infty}}\leq r_2$, $j=1,...,d$, and any time interval $I\subset \R$ of length at most $\tau$, the solution $\psi$ belongs to the Bourgain-Klainerman-Machedon space $X^{s,\frac12+\delta}_I$; see Definition~\ref{XsbT} below.

\smallskip

\item (Linear case, low regularity). Assume $\kappa=0$, $\tfrac13\geq s>0$, and that Assumption~\ref{AssumptionLWP} holds. For every $\psi_0\in H^s(\T;\C)$, there exists a global mild solution $\psi\in C^0(\R; H^s(\T;\C))$ satisfying Duhamel's formula \eqref{duhamel_lin} and mass conservation, which is unique in $X^{s,\tfrac12+\delta}_I$ for all $I$ as in \textup{(i)}.

\smallskip

\item (Nonlinear case, high regularity). Assume $\kappa\neq 0$, $s>\tfrac12$, and that Assumption~\ref{AssumptionLWP} holds. For every $r_1,r_2,R>0$, there exists $\tau(r_1,r_2,R)>0$ such that for any $u\in L^2(\R;\R^d)$, $\mu\in\F\ell^{s,\infty}(\T;\R^d)$ and $\psi_0\in H^s(\T;\C)$ with $\|u_j\|_{L^2}\leq r_1$, $\|\mu_j\|_{\F\ell^{s,\infty}}\leq r_2$, $j=1,...,d$, and $\|\psi_0\|_{H^s}\leq R$, there exists a unique local mild solution $\psi\in C^0([0,\tau]; H^s(\T;\C))$ to \eqref{NLS}, that is a solution to the Duhamel formula:
\begin{align}
    \label{duhamel_nonlin}\psi(t)=e^{- it\dx^2}\psi_0-i\int_0^te^{ -i(t-t')\dx^2}\big\{(u(t')\cdot\mu)\psi(t') +\kappa|\psi(t')|^{2m}\psi(t')\big\}dt',
\end{align}
for any $t\in [0,\tau]$, and that satisfies $\|\psi(t)\|_{L^2}=\|\psi_0\|_{L^2}$ for any $t\in [0,\tau]$. 

\smallskip Moreover, there is $\delta(s)>0$ such that $\psi\in X^{s,\frac12+\delta}_\tau$.

 \smallskip
 
 \item (Nonlinear case, low regularity). Assume $\kappa\neq 0$, $\tfrac12\geq s>s_c(m):=\max\big\{\tfrac12-\tfrac1m,0\big\} $, and that Assumption~\ref{AssumptionLWP} holds. For every $r_1,r_2,R>0$, there exists $\tau(r_1,r_2,R)>0$ such that for any $u\in L^2(\R;\R^d)$, $\mu\in\F\ell^{s,\infty}(\T;\R^d)$ and $\psi_0\in H^s(\T;\C)$ with $\|u_j\|_{L^2}\leq r_1$, $\|\mu_j\|_{\F\ell^{s,\infty}}\leq r_2$, $j=1,...,d$, and $\|\psi_0\|_{H^s}\leq R$, there exists a local mild solution $\psi\in C^0([0,T]; H^s(\T;\C))$ to \eqref{NLS} solving Duhamel's formula \eqref{duhamel_nonlin} and satisfying mass conservation, which is unique in $X^{s,\frac12+\delta}_T$ with $\delta(s)$ as in (iii).

\smallskip 
\item (Globalization of perturbation solutions) Let any $T>0$. Then there exists $\eta=\eta(T,\psi_\reff,\mu)>0$ such that if $\|\psi_0-\psi_\reff(0)\|_{H^s}\leq \eta$, then the unique local solution $\psi$ extends to $[0,T]$.

\smallskip
\item (Regularity of the flow map). In addition, for any  $r_1,r_2,R>0$, with $\tau(r_1,r_2,R)>0$ as above, the flow map
\begin{equation*}
\Phi : \begin{cases}
B_{L^2}(r_1)^d\times B_{\F\ell^{s,\infty}}(r_2)^d\times B_{H^s}(R)&\to  C^0\big([0,\tau];H^s(\T,\C)\big)\\
\hfill (u,\mu,\psi_0)&\mapsto \psi.
 \end{cases}
\end{equation*}
is smooth, where $B_{X}(R)$ denotes the ball of radius $R$ centered at the origin in the space $X$.
\end{enumerate}
\end{theorem}
 Recall that in (v) above, 
 \begin{equation}\label{psireff}
\psi_\reff(t)= \frac1{\sqrt{2\pi}} e^{-i \frac{\kappa}{(2\pi)^m} t}
\end{equation}
is the unique global solution to \eqref{NLS} with $u\equiv 0$ and initial data $\psi_\reff(0)=e_0$.

\begin{remark}[On the notion of solution and uniqueness]\label{RM:uniqueness}\rm~
\begin{enumerate}[label=\alph*)]
\item Note that in Theorem \ref{THM:LWP} the solution $\psi$ is defined in the mild sense, satisfying Duhamel's formula \eqref{duhamel_lin}, and \eqref{duhamel_nonlin}, respectively. Consequently, we implicitly assume that the right-hand side is well-defined as an $H^s$-valued distribution for each $t\in[0,T]$. To justify this, we first make a fixed-point argument directly in the space $X^{s,\frac12+\delta}_I$, which is a closed subspace of $C(I;H^s(\T;\C))$. For a function $\psi\in X^{s,\frac12+\delta}_I$ and a potential $u\cdot \mu$ with $u\in L^2$ and $\mu\in\F\ell^{s,\infty}$, Proposition~\ref{PROP:bilinear} below shows that the product $(u\cdot\mu)\psi$ is indeed well-defined as a space-time distribution, belonging to the space $X^{s,-\frac12+\delta}_I$. This shows existence of a solution in $X^{s,\frac12+\delta}_I$, a priori only unique in this closed subspace of $C(I;H^s(\T;\C))$. In the ``high'' regularity cases (i) and (iii) of Theorem~\ref{THM:LWP}, the uniqueness statement can be upgraded to all mild solutions in the class $C([0,T];H^s(\T;\C))$, regardless of how these solutions have been constructed (i.e. how one has been able to make sense of the right-hand side in Duhamel's formula), and the fact that this unique solution actually belongs to the smaller space $X^{s,\frac12+\delta}_I$ is given as an extra property. On the other hand, in the low regularity cases (ii) and (iv), we are only able to show uniqueness in the smaller class $X^{s,\frac12+\delta}_I$ since $\mu$ and $\psi$ lack of sufficient spatial regularity to be able to make sense of $(u\cdot\mu)\psi$ without knowing that $\psi\in X^{s,\frac12+\delta}_I$. Moreover, due to the poor temporal regularity of $u$, usual techniques to get unconditional uniqueness in $C([0,T];H^s(\T;\C))$ such as Poincaré-Dulac normal form reductions do not seem to be implementable.
\item From the perspective of control theory, the local exact controllability result of Theorem~\ref{THM:contr_ground} at low regularity, building on the well-posedness provided by Theorem~\ref{THM:LWP}, shows that for suitable electric potentials $\mu$, it is possible to act on the dynamics of \eqref{NLS} by tuning the amplitudes $u$, but a priori only for the solution we constructed in $X^{s,\frac12+\delta}_I$; it is left open wether other solutions could exist, and therefore how these other trajectories would react when using the controls $u$ suited for the unique $X^{s,\frac12+\delta}_I$ trajectory.
\end{enumerate}
\end{remark}

\begin{remark}[On the time of existence]~\\
\rm
We remark that in the nonlinear case $\kappa\neq 0$, Theorem~\ref{THM:LWP} only claims existence and uniqueness of a local solution. The easiest way to globalize local solutions to a nonlinear dispersive PDE is to work at the regularity level of a conservation law. In the case of \eqref{NLS}, due to the time dependence of the control $u$, the energy associated with the (time-dependent) Hamiltonian of the equation is no longer a conserved quantity. The mass is therefore the only available conserved quantity, but we cannot reach well-posedness in the critical space $L^2$; recall the discussion after Theorem~\ref{THM:LWP_intro}. In \cite{BL}, the authors show global well-posedness of \eqref{NLS} for small control and data close to a constant. In our case however, our estimates used to reach low regularity spaces rely crucially on the fact that the eigenvalues of $-\dx^2$ are squares. Therefore the strategy of \cite[Proposition 9]{BL} to write Duhamel's formula based on the linear operator obtained by linearizing the equation round the constant trajectory does not seem to be implementable. But this is not an issue since, provided that we work closeley enough to the constant trajectory, the control is small and the exact controllability can be obtained in arbitrarily small time, in particular smaller than the maximal time of existence of the solution.
\end{remark}

\subsection{Bilinear and multilinear estimates}
In order to solve \eqref{NLS} through a fixed-point argument, we need several multilinear estimates to handle both terms $(u\cdot\mu)\psi$ and $|\psi|^{2m}\psi$ in \eqref{NLS}.

\subsubsection{Bilinear estimates}\label{SEC:bilinear}
We start with the estimate used to handle the bilinear term in \eqref{NLS}. In \cite{BL,Alessandro}, Ingham's inequality is used  to prove the well-posedness of the linearized equation with control in $L^2$. The upper bound in Ingham's inequality can be seen as an $L^\infty_xL^2_t$ linear estimate on the free solution $e^{-it\dx^2}\varphi$. To tackle the full bilinear term, we thus need a bilinear version of this estimate, which is therefore an $L^\infty_xL^4_t$ type estimate. Let us point out that while the $L^\infty_x L^2_t$ estimate scales like Zygmund's estimate
\begin{equation}\label{Zygmund}
\|e^{-it\dx^2}\varphi\|_{L^4_{t,x}(\T\times\T)}\les \|\varphi\|_{L^2};
\end{equation}
(see e.g. \cite{ET}), the $L^\infty_x L^4_t$ type estimate scales like Bourgain's $L^6_{t,x}$ estimate \cite{Bourgain}
\begin{equation}\label{BourgainL6}
\|e^{-it\dx^2}\varphi\|_{L^6_{t,x}(\T\times\T)}\les \|\varphi\|_{H^\varepsilon}.
\end{equation} 
Following the recent breakthrough of Bourgain-Demeter \cite{BD}, sharp $L^p_{t,x}$ estimates on tori $\T^d$ are known to hold as in the Euclidean case $\R^d$, up to endpoints, at which there is an unavoidable, arbitrarily small, loss of derivative. The full range of $L^p_tL^q_x$ Strichartz estimates with arbitrary loss on tori is still open though; in particular, we are not aware of a sharp (up to $\varepsilon$ loss) $L^4_tL^\infty_x$ estimate on $\T$. The best result seems to be a $\frac16+$ loss of derivatives by using Sobolev and H\"older inequalities with the scale invariant $L^6_{t,x}$ estimate of Bourgain with arbitrary loss. Under the milder condition $V\in L^2_t W^{s,1}_x$ (compared to Assumption~\ref{AssumptionLWP}), together with Galilean invariance and a similar bilinear analysis as the one below, this Strichartz estimate would yield well-posedness of \eqref{NLS} in $H^s$ with the extra condition $s>\frac13$ instead of the whole sub-critical range $s>\max(0,s_c(m))$ as in Theorem~\ref{THM:LWP_intro}.  Here we rely instead on (multilinear refinements of) $L^\infty_xL^4_t$ type estimates, which we show to hold with arbitrary loss. Actually, our estimate is slightly stronger than $L^\infty_x L^4_t$, as it allows to deal with a space-time potential $u\cdot\mu$ with $u\in L^2(\R;\R^d)$ and $\mu\in\F\ell^{s,\infty}(\T;\R^d)$ instead of just $\mu\in W^{s,1}(\T;\R^d)$ (or $\jb{\nabla}^s\mu\in \M(\T;\R^d)$, the space of $\R^d$-valued signed Radon measures).

\begin{lemma}\label{LEM:bilinear}
For any $\varepsilon>0$, there exists $C,\delta>0$ such that for any $\psi,\varphi\in X^{0,\frac12-\delta}$, and any dyadic $N,N_1,N_2$, the following holds:\\
\textup{(i)} If $N_1\not\sim N_2$, then
\begin{equation}\label{LinfL4bil}
\big\|(\ft{\P_{N_1}\psi})_{k_1}(\tau_1) \cj{(\ft{\P_{N_2}\varphi})_{k_1+k}}(\tau_1+\tau)\big\|_{L^2_{\tau}L^1_{\tau_1}\ell^1_{k,k_1}} \leq C(N_1\wedge N_2)^\varepsilon \|\P_{N_1}\psi\|_{X^{0,\frac12-\delta}}\|\P_{N_2}\varphi\|_{X^{0,\frac12-\delta}}.
\end{equation}
\textup{(ii)} We have 
\begin{equation}\label{LinfL4bil2}
\big\|\mathbf{1}_{\substack{k\neq 0\\ k\neq -2k_1}}(\ft{\P_{N_1}\psi})_{k_1}(\tau_1) \cj{(\ft{\varphi})_{k_1+k}}(\tau_1+\tau)\big\|_{L^2_{\tau}L^1_{\tau_1}\ell^1_{k,k_1}} \leq CN_1^\varepsilon \|\P_{N_1}\psi\|_{X^{0,\frac12-\delta}}\|\varphi\|_{X^{0,\frac12-\delta}}.
\end{equation}
\textup{(iii)} If $N\ll N_1\sim N_2$ and $I_1,I_2\subset \R$ are some intervals of size $|I_1|+|I_2|\les N$ then
\begin{multline}\label{LinfL4bil3}
\big\|\mathbf{1}_{|k|\sim N}\mathbf{1}_{k\neq 0}(\ft{\P_{I_1}\P_{N_1}\psi})_{k_1}(\tau_1) \cj{(\ft{\P_{I_2}\P_{N_2}\varphi})_{k_1+k}}(\tau_1+\tau)\big\|_{L^2_{\tau}L^1_{\tau_1}\ell^1_{k,k_1}} \\\leq CN^\varepsilon \|\P_{I_1}\P_{N_1}\psi\|_{X^{0,\frac12-\delta}}\|\P_{I_2}\P_{N_2}\varphi\|_{X^{0,\frac12-\delta}}.
\end{multline}
\end{lemma}

\begin{remark}\label{REM:Strichartz}\rm
Using Plancherel theorem and Lemma~\ref{LEM:bilinear}~(ii), one indeed obtains the space-time estimate
\begin{align*}
\|\P_N\psi\|_{L^\infty_xL^4_t}^2&\leq \big\|\P_N\psi\cj{\P_N\psi}\big\|_{L^\infty_x L^2_t}\\
& \lesssim \big\|\mathbf{1}_{k_1^2\neq (k+k_1)^2}(\ft{\P_N\psi})_{k_1}(\tau_1)\cj{(\ft{\P_N\psi})}_{k+k_1}(\tau+\tau_1)\big\|_{L^2_\tau L^1_{\tau_1}\ell^1_{k,k_1}}\\
&\qquad+\big\|(\ft{\P_N\psi})_{k_1}(\tau_1)\cj{(\ft{\P_N\psi})}_{k_1}(\tau+\tau_1)\big\|_{L^2_\tau L^1_{\tau_1}\ell^1_{k_1}}\\
&\qquad+\big\|\mathbf{1}_{k_1\neq 0}(\ft{\P_N\psi})_{k_1}(\tau_1)\cj{(\ft{\P_N\psi})}_{-k_1}(\tau+\tau_1)\big\|_{L^2_\tau L^1_{\tau_1}\ell^1_{k_1}}\\
&\les N^\epsilon\|\P_N\psi\|_{X^{0,\frac12-\delta}}^2,
\end{align*}
where the resonant terms are dealt with Cauchy-Schwarz inequality in $k_1$ and a bilinear refinement of Sobolev inequality in $\tau_1,\tau_2$, as in Case 3 of the proof of Proposition~\ref{PROP:bilinear} below. In particular this implies the $L^\infty_xL^4_T$ Strichartz estimate for $T\in(0,1]$ by taking $\psi(t)=\eta(t)e^{it\dx^2}\psi_0$ for some smooth bump function $\eta$.
\end{remark}

\begin{proof}[Proof of Lemma~\ref{LEM:bilinear}]
We start with the proof of \eqref{LinfL4bil}. Decomposing dyadically in the conormal variable and using Cauchy-Schwarz inequality, we can estimate
\begin{align}\label{bilinear}
&\big\|(\ft{\P_{N_1}\psi})_{k_1}(\tau_1) \cj{(\ft{\P_{N_2}\varphi})_{k_1+k}}(\tau_1+\tau)\big\|_{L^2_{\tau}L^1_{\tau_1}\ell^1_{k,k_1}}\notag\\
&\qquad\leq  \sum_{L_1,L_2}\big\|(\ft{\P_{N_1}\Q_{L_1}\psi})_{k_1}(\tau_1) \cj{(\ft{\P_{N_2}\Q_{L_2}\varphi})_{k_1+k}}(\tau_1+\tau)\big\|_{L^2_{\tau}L^1_{\tau_1}\ell^1_{k,k_1}}\notag\\
&\qquad\leq \sum_{L_1,L_2} \|\P_{N_1}\Q_{L_1}\psi\|_{L^2_{t,x}} \|\P_{N_2}\Q_{L_2}\varphi\|_{L^2_{t,x}}\big\|\mathbf{1}_{A_\tau}(\tau_1,k_1,k_2)\big\|_{L^\infty_\tau L^2_{\tau_1}\ell^2_{k_1,k_2}},
\end{align}
where
$$A_\tau=\big\{(\tau_1,k_1,k_2)\in\R\times\Z^2,~|k_1|\sim N_1,~|k_2|\sim N_2,~
|\tau_1-k_1^2|\sim L_1,~|\tau_1+\tau-k_2^2|\sim L_2\big\}.$$
From the conditions $|\tau_1-k_1^2|\sim L_1$ and $|\tau_1+\tau-k_2^2|\sim L_2$, we have that $\tau_1$ lies in a set of size $L_1\wedge L_2$, and $|\tau+k_1^2-k_2^2|\les L_1\vee L_2$. Namely, 
$$k_1^2-k_2^2\in J_\tau \text{ with } J_\tau=[-\tau-C(L_1\vee L_2),-\tau+C(L_1\vee L_2)]$$ for some constant $C>0$.
Therefore, since $N_1\not\sim N_2$, using Lemma~\ref{LEM:counting}~(i) we can estimate for arbitrary $\nu>0$ and $\theta\in(0,1)$:
\begin{align*}
&\big\|\mathbf{1}_{A_\tau}(\tau_1,k_1,k_2)\big\|_{L^\infty_\tau L^2_{\tau_1}\ell^2_{k_1,k_2}}\\
&\qquad\les (L_1\wedge L_2)^\frac12 \sup_{\tau\in\R}\#\big\{(k_1,k_2)\in\Z^2,~|k_1|\sim N_1,~|k_2|\sim N_2,~|\tau+k_1^2-k_2^2|\les L_1\vee L_2\big\}^\frac12\\
&\qquad\les (L_1\wedge L_2)^\frac12 \min\big((L_1\vee L_2)(N_1\wedge N_2)^\nu,(L_1\vee L_2)^\frac12(N_1\wedge N_2)\big)^\frac12\\
&\qquad\les (L_1L_2)^{\frac12-\frac\theta8}(N_1\wedge N_2)^{\frac\theta2+\nu(1-\theta)}.
\end{align*}
Plugging this bound in \eqref{bilinear} and summing using Cauchy-Schwarz inequality, we get
\begin{align*}
&\big\|(\ft{\P_{N_1}\psi})_{k_1}(\tau_1) \cj{(\ft{\P_{N_2}\varphi})_{k_1+k}}(\tau_1+\tau)\big\|_{L^2_{\tau}L^1_{\tau_1}\ell^1_{k,k_1}}\\
&\qquad\les \sum_{L_1,L_2}(L_1L_2)^{\frac12-\frac\theta8}(N_1\wedge N_2)^{\frac\theta2+\nu(1-\theta)}\|\P_{N_1}\Q_{L_1}\psi\|_{L^2_{t,x}} \|\P_{N_2}\Q_{L_2}\varphi\|_{L^2_{t,x}}\\
&\qquad\les (N_1\wedge N_2)^{\frac\theta2+\nu(1-\theta)}\|\P_{N_1}\psi\|_{X^{0,\frac12-\frac\theta9}} \|\P_{N_2}\varphi\|_{X^{0,\frac12-\frac\theta9}}.
\end{align*}
This is enough for \eqref{LinfL4bil} by taking $0<\theta,\nu\ll \varepsilon$ and then $\delta=\frac{\theta}{9}$.

The proof of \eqref{LinfL4bil2} is similar: now we need to estimate
$$\sup_{\tau\in\R}\#\big\{(k_1,k_2)\in\Z^2,~k_1\neq k_2,~k_1\neq -k_2,~|k_1|\sim N_1,~|\tau+k_1^2-k_2^2|\les L_1\vee L_2\big\}.$$
We can now invoke Lemma~\ref{LEM:counting}~(ii), the condition $m\neq 0$ thereof being ensured by $k_1\neq k_2$ and $k_1\neq -k_2$. This gives the bound
\begin{align*}
&\big\|\mathbf{1}_{A_\tau}(\tau_1,k_1,k_2)\big\|_{L^\infty_\tau L^2_{\tau_1}\ell^2_{k_1,k_2}}\les (L_1\wedge L_2)^\frac12 \min\big((L_1\vee L_2)N^\nu,(L_1\vee L_2)^\frac12N\big)^\frac12
\end{align*}
which is again enough for \eqref{LinfL4bil2}.
 
At last, the proof of \eqref{LinfL4bil3} follows the same lines as before, estimating now
\begin{multline*}
(L_1\wedge L_2)^\frac12\sup_{\tau\in\R}\#\big\{(k_1,k_2)\in\Z^2,~k_1\neq k_2,~|k_1-k_2|\sim N\ll |k_1|\sim|k_2|\sim N_1,~k_1\in I_1,~k_2\in I_2,\\ |\tau+k_1^2-k_2^2|\les L_1\vee L_2\big\}^\frac12\\ \les (L_1\wedge L_2)^\frac12\min\big((L_1\vee L_2)N^\nu,N^2\big)^\frac12\les (L_1L_2)^{\frac12-\frac\theta4}N^{\nu\frac{(1-\theta)}2+\theta}
\end{multline*}
thanks to Lemma~\ref{LEM:counting}~(iii), the condition $m\neq 0$ being ensured by $k_1\neq k_2$ and $|k_1+k_2|\sim N_1\neq 0$. This is again enough for \eqref{LinfL4bil3}.
\end{proof}

\begin{proposition}\label{PROP:bilinear}
For any $s>0$, there exists $\delta>0$ small enough, and $\theta,C>0$, such that for any $u\in L^2(\R;\R^d)$, $\mu\in\F\ell^{s,\infty}(\T;\R^d)$, $T\in(0,1)$, and $\psi\in X^{s,\frac12+\frac\delta2}_T(\T;\C))$, it holds
\begin{equation}\label{keybil}
\Big\|\int_0^te^{-i(t-t')\dx^2}\big\{(u\cdot\mu)\psi\big\}(t')dt'\Big\|_{X^{s,\frac12+\frac\delta2}_T}\leq C T^{\theta}\big(\sum_{j=1}^d\|u_j\|_{L^2}\|\mu_j\|_{\F\ell^{s,\infty}}\big)\|\psi\|_{X^{s,\frac12-\delta}_T}.
\end{equation}
\end{proposition}
\begin{proof}
By Lemma~\ref{LEM:linear}~(iv) and~(ii), extension and duality, we have
\begin{align*}
\Big\|\int_0^te^{-i(t-t')\dx^2}\big\{(u\cdot\mu)\psi\big\}(t')dt'\Big\|_{X^{s,\frac12+\frac\delta2}_T}&\les \big\|\jb{\nabla}^s\big\{(u\cdot\mu)\psi\big\}\big\|_{X^{0,-\frac12+\frac\delta2}_T}\\
&\les T^\frac\delta2 \big\|\jb{\nabla}^s\big\{(u\cdot\mu)\psi\big\}\big\|_{X^{0,-\frac12+\delta}_T}\\
&\les T^\frac\delta2\sup_{\|w\|_{X^{0,\frac12-\delta}}\leq 1}   \int_\R\int_\T\jb{\nabla}^s\big\{(u\cdot\mu)\Psi\big\}\cj wdxdt
\end{align*}
where $\Psi$ is an extension of $\psi$ from $[-T,T]$ to $\R$ satisfying $\|\Psi\|_{X^{s,\frac12-\delta}}\leq 2\|\psi\|_{X^{s,\frac12-\delta}_T}$.

We then decompose dyadically
\begin{align*}
\int_\R\int_\T \jb{\nabla}^s\big\{(u\cdot\mu)\Psi\big\}\cj wdxdt&= \sum_{\substack{N_0,N_1,N_2}}\int_\R\int_\T \jb{\nabla}^s\big\{(u\cdot\P_{N_0}\mu) \P_{N_1}\Psi\big\}\cj{\P_{N_2}w}dxdt.
\end{align*}
Note that using Plancherel identity, the integral on $\T$ can be written as
\begin{align*}
\sum_{\substack{k,k_1,k_2\in\Z\\k+k_1=k_2}} (\ft{\P_{N_0}\mu})_{k} (\ft{\P_{N_1}\Psi})_{k_1}\jb{k_2}^s\cj{(\ft{\P_{N_2}w})_{k_2}},
\end{align*}
and from the localizations of $k_1$, $k_2$ and $k$ together with the convolution constraint $k_2=k+k_1$, we must have $N_0^*\sim N_1^*\gtrsim N_2^*$.

We now perform a case-by-case analysis on the bilinear interaction. 

\medskip
\textbf{Case 1: High $\times$ Low.}\\
In this case, $N_1\gg N_2$, therefore $N_0\sim N_1$. Using Plancherel theorem, Cauchy-Schwarz inequality, and Lemma~\ref{LEM:bilinear}~(i), we estimate
\begin{align*}
&\Big|\int_\R\int_\T \jb{\nabla}^s\big\{(u\cdot\P_{N_0}\mu) \P_{N_1}\Psi\big\}\cj{\P_{N_2}w}dxdt\Big|\\
&\qquad\les \sum_{j=1}^d\|u_j\|_{L^2_t}\|\mu_j\|_{\F\ell^{s,\infty}}N_0^{-s}N_2^s\big\|(\ft{\P_{N_1}\psi})_{k_1}(\tau_1) \cj{(\ft{\P_{N_2}w})_{k_1+k}}(\tau_1+\tau)\big\|_{L^2_{\tau}L^1_{\tau_1}\ell^1_{k,k_1}}\\
&\qquad \les \sum_{j=1}^d\|u_j\|_{L^2_t}\|\mu_j\|_{\F\ell^{s,\infty}}N_0^{-s}N_2^{s+\varepsilon}N_1^{-s} \|\P_{N_1}\psi\|_{X^{s,\frac12-\delta}}\|\P_{N_2}w\|_{X^{0,\frac12-\delta}}.
\end{align*}
This is enough to sum on $N_0\sim N_1\gg N_2$ using Cauchy-Schwarz inequality, provided that $s>\varepsilon$.

\medskip
\textbf{Case 2: Low $\times$ High.}\\
In this case, $N_2\gg N_1$, therefore $N_0\sim N_2$. Proceeding exactly as above would now yield the bound 
$$\sum_{N_0\sim N_2\gg N_1}N_1^{\varepsilon-s} \|\P_{N_1}\psi\|_{X^{s,\frac12-\delta}}\|\P_{N_2}w\|_{X^{0,\frac12-\delta}},$$
which is barely not summable in $N_2$ (after bounding the sum on $N_0$, for fixed $N_2$, by $O(1)$). Therefore we need to modify the previous argument to avoid this logarithmic divergence. We rewrite this contribution as 
\begin{align*}
&\sum_{\substack{N_0,N_1,N_2\\ N_0\sim N_2\gg N_1}}\int_\R\int_\T \jb{\nabla}^s\big\{(u\cdot\P_{N_0}\mu) \P_{N_1}\Psi\big\}\cj{\P_{N_2}w}dxdt\\
&=\sum_{\substack{N_0,N_1,N_2\\ N_0\sim N_2\gg N_1}}\int_\R\int_\T u\cdot\big\{\jb{\nabla}^s\big(\P_{N_0}\mu \P_{N_1}\Psi\big)-\big(\jb{\nabla}^s\P_{N_0}\mu\big)\P_{N_1}\Psi\big\}\cj{\P_{N_2}w}dxdt\\
&\qquad+\sum_{N_1}\int_\R\int_\T \big(u\cdot\jb{\nabla}^s\mu\big)\P_{N_1}\Psi\cj{\P_{\gg N_1}w}dxdt.
\end{align*}
As for the second term, using Cauchy-Schwarz inequality as in Case 1, we find 
\begin{align*}
&\Big|\int_\R\int_\T \big(u\cdot\jb{\nabla}^s\mu\big)\P_{N_1}\Psi\cj{\P_{\gg N_1}w}dxdt\Big|\\
&\qquad\les \sum_{j=1}^d\|u_j\|_{L^2_t}\|\mu_j\|_{\F\ell^{s,\infty}}\big\|\mathbf{1}_{\substack{k\neq 0\\k\neq -2k_1}}(\ft{\P_{N_1}\psi})_{k_1}(\tau_1) \cj{(\ft{\P_{\gg N_1}w})_{k_1+k}}(\tau_1+\tau)\big\|_{L^2_{\tau}L^1_{\tau_1}\ell^1_{k,k_1}}.
\end{align*}
Note that the conditions $k\neq 0$ and $k\neq -2k_1$ follow from the Fourier supports of $\P_{N_1}\Psi$ and $\P_{\gg N_1}w$. We are therefore entitled to use Lemma~\ref{LEM:bilinear}~(ii) to bound this last term by
\begin{align*}
\sum_{j=1}^d\|u_j\|_{L^2_t}\|\mu_j\|_{\F\ell^{s,\infty}}N_1^{\varepsilon-s}\|\P_{N_1}\psi\|_{X^{s,\frac12-\delta}}\|w\|_{X^{0,\frac12-\delta}}.
\end{align*}
This is now enough to sum on $N_1$ using Cauchy-Schwarz inequality, provided that $s>\varepsilon$.

As for the commutator above, we use Plancherel theorem, Cauchy-Schwarz inequality, and the mean value inequality together with the Fourier support of $\P_{N_0}\mu$, $\P_{N_1}\Psi$ and $\P_{N_2}w$ to estimate
\begin{align*}
&\Big|\int_\R\int_\T u\cdot\big\{\jb{\nabla}^s\big(\P_{N_0}\mu \P_{N_1}\Psi\big)-\big(\jb{\nabla}^s\P_{N_0}\mu\big)\P_{N_1}\Psi\big\}\cj{\P_{N_2}w}dxdt\Big|\\
&\qquad=\Big|\int_\R\int_\R \ft{u}(\tau)\cdot\sum_{k,k_1\in\Z}\big(\jb{k+k_1}^s-\jb{k}^s\big)(\ft{\P_{N_0}\mu})_k(\ft{\P_{N_1}\Psi})_{k_1}(\tau_1)\cj{(\ft{\P_{N_2}w})_{k+k_1}}(\tau+\tau_1)d\tau_1d\tau\Big|\\
&\qquad\les N_1N_2^{-1}\sum_{j=1}^d\|u_j\|_{L^2_t}\|\mu_j\|_{\F\ell^{s,\infty}}\big\|(\ft{\P_{N_1}\Psi})_{k_1}\cj{(\ft{\P_{N_2}w})_{k+k_1}}(\tau+\tau_1)\big\|_{L^2_\tau L^2_{\tau_1}\ell^1_{k,k_1}}.
\end{align*}
Using again Lemma~\ref{LEM:bilinear}~(i), we therefore get
\begin{align*}
&\Big|\int_\R\int_\T u\cdot\big\{\jb{\nabla}^s\big(\P_{N_0}\mu \P_{N_1}\Psi\big)-\big(\jb{\nabla}^s\P_{N_0}\mu\big)\P_{N_1}\Psi\big\}\cj{\P_{N_2}w}dxdt\Big|\\
&\qquad\les N_1^{1+\varepsilon-s}N_2^{-1}\sum_{j=1}^d\|u_j\|_{L^2_t}\|\mu_j\|_{\F\ell^{s,\infty}}\|\P_{N_1}\Psi\|_{X^{s,\frac12-\delta}}\|\P_{N_2}w\|_{X^{0,\frac12-\delta}},
\end{align*}
which is also summable in $N_0\sim N_2\gg N_1$ using Cauchy-Schwarz inequality, provided again that $s>\varepsilon$.

\medskip
\textbf{Case 3: High$\times$High $\to$ High.}
In this case we have $N_0\sim N_1\sim N_2$, and we have to deal with the resonances. We therefore decompose
\begin{align*}
&\int_\R\int_\T \jb{\nabla}^s\big\{(u\cdot\P_{N_0}\mu) \P_{N_1}\Psi\big\}\cj{\P_{N_2}w}dxdt\\
&\qquad = \int_\R\int_\R \ft{u}(\tau)\cdot\sum_{\substack{k,k_1\in\Z\\ k\neq 0\\ k\neq -2k_1}} \jb{k+k_1}^s(\ft{\P_{N_0}\mu})_k (\ft{\P_{N_1}\Psi})_{k_1}(\tau_1)\cj{(\ft{\P_{N_2}w})_{k+k_1}}(\tau+\tau_1)d\tau_1d\tau\\
&\qquad +  \int_\R\int_\R \ft{u}(\tau)\cdot\sum_{k_1\in\Z} \jb{k_1}^s(\ft{\P_{N_0}\mu})_0 (\ft{\P_{N_1}\Psi})_{k_1}(\tau_1)\cj{(\ft{\P_{N_2}w})_{k_1}}(\tau+\tau_1)d\tau_1d\tau\\
&\qquad + \int_\R\int_\R \ft{u}(\tau)\cdot\sum_{k_1\in\Z^*} \jb{k_1}^s(\ft{\P_{N_0}\mu})_{-2k_1} (\ft{\P_{N_1}\Psi})_{k_1}(\tau_1)\cj{(\ft{\P_{N_2}w})_{-k_1}}(\tau+\tau_1)d\tau_1d\tau.
\end{align*}
The first term above is estimated as before by using Cauchy-Schwarz inequality and Lemma~\ref{LEM:bilinear}~(ii):
\begin{align*}
&\Big|\int_\R\int_\R \ft{u}(\tau)\cdot\sum_{\substack{k,k_1\in\Z\\ k\neq 0\\ k\neq -2k_1}} \jb{k+k_1}^s(\ft{\P_{N_0}\mu})_k (\ft{\P_{N_1}\Psi})_{k_1}(\tau_1)\cj{(\ft{\P_{N_2}w})_{k+k_1}}(\tau+\tau_1)d\tau_1d\tau\Big|\\
&\qquad\les N_2^{s}N_0^{-s}N_1^{\varepsilon-s}\sum_{j=1}^d\|u_j\|_{L^2_t}\|\mu_j\|_{\F\ell^{s,\infty}}\|\P_{N_1}\Psi\|_{X^{s,\frac12-\delta}}\|\P_{N_2}w\|_{X^{0,\frac12-\delta}}.
\end{align*}
This is enough to sum on $N_0\sim N_1\sim N_2$ using Cauchy-Schwarz inequality.

As for the first resonant term, proceeding similarly as in the proof of Lemma~\ref{LEM:bilinear} (see also Remark~\ref{REM:Strichartz}), we estimate
\begin{align*}
&\Big| \int_\R\int_\R \ft{u}(\tau)\cdot\sum_{k_1\in\Z} \jb{k_1}^s(\ft{\P_{N_0}\mu})_0 (\ft{\P_{N_1}\Psi})_{k_1}(\tau_1)\cj{(\ft{\P_{N_2}w})_{k_1}}(\tau+\tau_1)d\tau_1d\tau\Big|\\
&\qquad\les \sum_{L_1,L_2}(L_1\wedge L_2)^\frac12 N_1^s \sum_{j=1}^d\|u_j\|_{L^2_t}\|\mu_j\|_{\F\ell^{0,\infty}}\|\P_{N_1}\Q_{L_1}\Psi\|_{L^2_{t,x}}\|\P_{N_2}Q_{L_2}w\|_{L^2_{t,x}}\\
&\qquad\les \sum_{j=1}^d\|u_j\|_{L^2_t}\|\mu_j\|_{\F\ell^{0,\infty}}\|\P_{N_1}\Psi\|_{X^{s,\frac14+\varepsilon}}\|\P_{N_2}w\|_{X^{0,\frac14+\varepsilon}}.
\end{align*}
This is also enough to sum on $N_0\sim N_1\sim N_2$. The last term is even easier due to the extra $s$ regularity:
\begin{align*}
&\Big|\int_\R\int_\R \ft{u}(\tau)\cdot\sum_{k_1\in\Z^*} \jb{k_1}^s(\ft{\P_{N_0}\mu})_{-2k_1} (\ft{\P_{N_1}\Psi})_{k_1}(\tau_1)\cj{(\ft{\P_{N_2}w})_{-k_1}}(\tau+\tau_1)d\tau_1d\tau\Big|\\
&\qquad\les N_0^{-s}\sum_{j=1}^d\|u_j\|_{L^2_t}\|\mu_j\|_{\F\ell^{s,\infty}}\|\P_{N_1}\Psi\|_{X^{s,\frac14+\varepsilon}}\|\P_{N_2}w\|_{X^{0,\frac14+\varepsilon}}.
\end{align*}

\medskip
\textbf{Case 4: High$\times$High $\to$ Low.}
In this case, $N_1\sim N_2\gg N_0$.

The case $N_0=0$, i.e. the resonance $k=0$, is estimated as in the previous case, therefore we focus on the case $N_0\geq 1$. In order to still benefit from the $s$ regularity of $\mu$ and therefore to only have a loss on the lowest frequency $N_0$ in the Strichartz type estimate, we need to make a further decomposition of the Fourier support of $\P_{N_1}\Psi$ and $\P_{N_2}w$. For $j\in\Z$ let then $I_j:=[jN_0,(j+1)N_0)$. We then have 
\begin{align*}
&\int_\R\int_\T \jb{\nabla}^s\big\{(u\cdot\P_{N_0}\mu) \P_{N_1}\Psi\big\}\cj{\P_{N_2}w}dxdt\\
&\qquad =\sum_{j_1,j_2\in\Z}\int_\R\int_\R \ft{u}(\tau)\cdot\sum_{k,k_1\in\Z} \jb{k+k_1}^s(\ft{\P_{N_0}\mu})_k (\ft{\P_{I_{j_1}}\P_{N_1}\Psi})_{k_1}(\tau_1)\cj{(\ft{\P_{I_{j_2}}\P_{N_2}w})_{k+k_1}}(\tau+\tau_1)d\tau_1d\tau.
\end{align*}
Moreover, from the convolution constraint $k_2=k_1+k$ with $|k|\sim N_0\ll |k_1|\sim N_1\sim N_2\sim|k_2|$ and $k_\ell\in I_{j_\ell}$, $\ell=1,2$, we find that $N_0\sim |k_2-k_1|\geq (|j_2-j_1|-1)N_0$. Thus $|j_2-j_1|\les 1$ and in particular the decomposition is almost orthogonal. Using then Cauchy-Schwarz inequality and Lemma~\ref{LEM:bilinear}~(iii) (note that we have $k\neq 0$ due to $N_0\geq 1$), we find
\begin{align*}
&\Big|\sum_{\substack{j_1,j_2\in\Z\\|j_2-j_1|\les 1}}\int_\R\int_\R \ft{u}(\tau)\cdot\sum_{k,k_1\in\Z} \jb{k+k_1}^s(\ft{\P_{N_0}\mu})_k (\ft{\P_{I_{j_1}}\P_{N_1}\Psi})_{k_1}(\tau_1)\cj{(\ft{\P_{I_{j_2}}\P_{N_2}w})_{k+k_1}}(\tau+\tau_1)d\tau_1d\tau\Big|\\
&\qquad\les \sum_{j=1}^d\|u_j\|_{L^2_t}\|\mu_j\|_{\F\ell^{s,\infty}}N_2^sN_0^{\epsilon-s}N_1^{-s}\sum_{\substack{j_1,j_2\in\Z\\|j_2-j_1|\les 1}}\|\P_{I_{j_1}}\P_{N_1}\Psi\|_{X^{s,\frac12-\delta}}\|\P_{I_{j_2}}\P_{N_2}w\|_{X^{0,\frac12-\delta}}.
\end{align*}
This is enough to sum on $j_1,j_2$ using Cauchy-Schwarz inequality and the almost orthogonality, and then on $N_1\sim N_2$ using also Cauchy-Schwarz inequality and almost orthogonality, and on $N_0$ provided again that $s>\varepsilon$.
\end{proof}

\subsubsection{Estimates for the nonlinear term}

In this subsection, we recall the multilinear estimates of Bourgain \cite{Bourgain} which are needed to estimate the nonlinearity in \eqref{NLS} when $\kappa\neq 0$. 

\begin{proposition}[\cite{Bourgain}]\label{PROP:Bourgain}
For any $s>s_c=\max(0,1-\frac1m)$, there exists $\delta>0$, such that for any choice of signs $\iota_j\in\{+,-\}$, $j=1,\dots,2m+1$, there exists $\theta,C>0$ such that for any $0<T\leq 1$ and any $\psi_1,\dots,\psi_{2m+1}\in X^{s,\frac12+\delta}_T$, it holds
\begin{align*}
    \Big\|\int_0^te^{i(t-t')\dx^2}\Big(\prod_{j=1}^{2m+1}\psi_j^{\iota_j}\Big)dt'\Big\|_{X^{s,\frac12+\delta}_T}\leq C T^{\theta}\prod_{j=1}^{2m+1}\|\psi_j\|_{X^{s,\frac12+\delta}_T}.
\end{align*}
\end{proposition}

\subsection{Proof of Theorem~\ref{THM:LWP}}
We finally turn to the proof of our well-posedness result.

\begin{proof}[Proof of Theorem~\ref{THM:LWP}]~

We start with the more general case (iv) of Theorem~\ref{THM:LWP}. Let $\psi_0\in H^s(\T;\C)$, $s>s_c(m)$, $u_j\in L^2(\R;\R)$ and $\mu_j\in\F\ell^{s,\infty}(\T;\R)$, $j=1,\dots,d$. Let $C_0$ be the constant in Lemma~\ref{LEM:linear}~(iii), and define
$$R = 2C_0\|\psi_0\|_{H^s}.$$ 
Moreover, let $$T = K^{-1}\big(1+\sum_{j=1}^d\|u_j\|_{L^2_t}\|\mu_j\|_{\F\ell^{s,\infty}}+|\kappa| R^{2m}\big)^{-\frac1\theta},$$
with $\theta=\theta_1\wedge \theta_2>0$, $\theta_1$ being given in Proposition~\ref{PROP:bilinear} and $\theta_2$ in Proposition~\ref{PROP:Bourgain}, and $K$ a constant to be chosen later, independent of $T$, $R$, and $\sum_{j=1}^d\|u_j\|_{L^2_t}\|\mu_j\|_{\F\ell^{s,\infty}}$. Note that $T\leq 1$ with this definition.

For $\delta=\delta_1\wedge\delta_2>0$ with $\delta_1$ as in Proposition~\ref{PROP:bilinear} and $\delta_2$ as in Proposition~\ref{PROP:Bourgain}, let us define the map
\begin{align*}
\Gamma_{T,\psi_0,u,\mu} : \psi\in X^{s,\frac12+\delta}_T \mapsto e^{- it\dx^2}\psi_0-i\int_0^te^{ -i(t-t')\dx^2}\big\{(u(t')\cdot\mu )\psi(t') +\kappa|\psi(t')|^{2m}\psi(t')\big\}dt',
\end{align*}

We will show that $\Gamma$ defines a contraction on the ball
$$B=\big\{\psi\in X^{s,\frac12+\delta}_T,~~\|\psi\|_{X^{s,\frac12+\delta}_T}\leq R\big\}.$$
Indeed, for $\psi\in B$, using Lemma~\ref{LEM:linear}~(iii), Proposition~\ref{PROP:bilinear}, and Proposition~\ref{PROP:Bourgain}, we have
\begin{align*}
\big\|\Gamma_{T,\psi_0,u,\mu} (\psi)\big\|_{X^{s,\frac12+\delta}_T}&\leq C_0\|\psi_0\|_{H^s}+ C_1T^\theta \Big\{\sum_{j=1}^d\|u_j\|_{L^2_t}\|\mu_j\|_{\F\ell^{s,\infty}}\|\psi\|_{X^{s,\frac12+\delta}_T}+|\kappa|\|\psi\|_{X^{s,\frac12+\delta}_T}^{2m+1}\Big\}\\
&\leq \frac12R+C_1T^\theta \Big\{\sum_{j=1}^d\|u_j\|_{L^2_t}\|\mu_j\|_{\F\ell^{s,\infty}}+|\kappa| R^{2m}\Big\}R.
\end{align*}
Therefore, taking $K\geq (2C_1)^{\frac1\theta}$, we have with our choice of $T$ that $\Gamma$ maps $B$ to itself. Similarly, using the bilinearity of the control term and the multilinearity of the nonlinear part, we find for any $\psi,\varphi\in B$ and $u_j,v_j\in L^2(\R;\R)$, $\mu_j,\nu_j\in\F\ell^{s,\infty}(\T;\R)$, $j=1,\dots,d$, that
\begin{align*}
&\big\|\Gamma_{T,\psi_0,u,\mu} (\psi)-\Gamma_{T,\varphi_0,v,\nu}(\varphi)\big\|_{X^{s,\frac12+\delta}_T}\\
&\qquad\leq C_0\|\psi_0-\varphi_0\|_{H^s}+C_2T^\theta \sum_{j=1}^d\Big\{\|u_j\|_{L^2_t}\|\mu_j\|_{\F\ell^{s,\infty}}\|\psi-\varphi\|_{X^{s,\frac12+\delta}_T}\\
&\qquad\qquad+\|u_j-v_j\|_{L^2_t}\|\mu_j\|_{\F\ell^{s,\infty}}\|\varphi\|_{X^{s,\frac12+\delta}_T}+\|v_j\|_{L^2_t}\|\mu_j-\nu_j\|_{\F\ell^{s,\infty}}\|\varphi\|_{X^{s,\frac12+\delta}_T}\Big\}\\
&\qquad\qquad\qquad+C_2T^\theta\big(\|\psi\|_{X^{s,\frac12+\delta}_T}+\|\varphi\|_{X^{s,\frac12+\delta}_T}\big)^{2m}\|\psi-\varphi\|_{X^{s,\frac12+\delta}_T}.
\end{align*}
Up to taking $K$ larger depending on $C_2$, this shows that $\Gamma_{T,\psi_0,u,\mu}$ defines a contraction on $B$, and thus has a unique fixed-point in $B$. The uniqueness in the whole class $X^{s,\frac12+\delta}_T$ follows from repeating the contraction property on a larger ball including any two possible solutions $\psi,\varphi$ for a possibly smaller time $\tau$, and then iterating the estimate to cover the whole time interval $[-T,T]$. The previous estimate also shows the local Lipschitz property of the flow with respect to $(\psi_0,u_j,\mu_j)\in H^s(\T;\C)\times L^2(\R;\R^d)\times \F\ell^{s,\infty}(\T;\R^d)$.

\medskip 
Moreover, since $\Gamma$ is real analytic as a polynomial in $(u,\mu,\psi_0)$, and since the contraction estimate shows that the differential in $\psi$ of $(u,\mu,\psi_0,\psi)\mapsto\psi-\Gamma(\psi)$ is invertible, then the smoothness of $(u,\mu,\psi_0)\mapsto \psi$ follows from the analytic implicit function theorem.

\medskip
We now prove the invariance of the mass. Our multilinear estimates show propagation of regularity, namely, if $\psi_0\in H^\sigma(\T;\C)$ and $\mu\in \F\ell^{\sigma,\infty}(\T;\R^d)$ for some $\sigma\geq s>0$, then the unique local solution $\psi \in X^{s,\frac12+\delta}_T$ actually belongs to $X^{\sigma,\frac12+\delta}_T\subset C([0,T];H^\sigma(\T;\C))$. Let then $\psi_{0,n}=\P_{\leq n}\psi_0\to \psi_0$ in $H^s$ and $\mu_n = \P_{\leq n}\mu \to \mu$ in $\F\ell^{s-\epsilon,\infty}$ for all $0<\epsilon<s$. In particular the continuity of the flow map in $H^{s-\epsilon}$ ensures that the unique solution $\psi_n$ emanating from $(u,\mu_n,\psi_{0,n})$ converges to $\psi$ in $X^{s-\epsilon,\frac12+\delta}_T$, and in particular in $C([0,T];L^2)$. By propagation of regularity, $\psi_n$ is smooth and is therefore a classical solution to \eqref{NLS} and we can compute
\begin{align*}
\frac{d}{dt}\frac12\|\psi_n(t)\|_{L^2}^2 = \Re\langle  \psi,\dt\psi\rangle_{L^2} = \Im\langle \psi, \dx^2\psi  +(u\cdot\mu)\psi +\kappa|\psi|^{2m}\psi\rangle_{L^2}.
\end{align*}
Since $u$ and $\mu$ are real-valued, we deduce that the imaginary part of the inner product above cancels out and the mass is invariant for $\psi_n$, and therefore for $\psi$ by taking $n\to\infty$.

\bigskip
It remains to discuss the globalization of solutions under the smallness assumption on $u_j$ and $\psi_0-\psi_\reff(0)$. Without loss of generality, we take $T\geq 1$. Let $R=\|e_0\|_{H^s}+1$ and $0<\eta\leq 1$ to be determined later such that $\eta\sum_{j=1}^d\|\mu_j\|_{\F\ell^{s,\infty}}\leq 1$. Then by the local theory with $r_1=1$ and $r_2=\|\mu\|_{\F\ell^{s,\infty}}$, there is $\tau=\tau(R)>0$ such that, provided that $\|u_j\|_{L^2_T}\leq \eta\leq 1$, for any data $\psi_0$ of size $\|\psi_0\|_{H^s}\leq R$, the solution exists on $[0,\tau]$. By the Lipschitz regularity of the flow map, there is $C(R)\geq 1$ such that
$$\|\psi-\psi_\reff\|_{L^\infty_\tau H^s}\leq C(R)\big(\|\psi_0-\psi_\reff(0)\|_{H^s}+\|u\|_{L^2_\tau}\big).$$
Therefore, if $\|\psi_0-\psi_\reff(0)\|_{H^s}\leq \eta$ then $\|\psi(\tau)-\psi_\reff(\tau)\|_{H^s}\leq 2C(R)\eta$. So if we take 
$$M=\lceil T/\tau\rceil \text{ and } \eta \sim \frac{1}{(M+1)C(R)^M},$$
then after $J\leq M-1$ iterations we have
$$\|\psi(J\tau)\|_{H^s}\leq \|\psi(J\tau)-\psi_\reff(J\tau)\|_{H^s}+\|e_0\|_{H^s} \leq (J+1)C(R)^J\eta +\|e_0\|_{H^s}\leq R,$$
so we can still apply the local theory on $[J\tau,(J+1)\tau]$. This extends the solution to $[0,T]$.

\medskip
Let us now discuss cases (i) and  (ii): $\kappa=0$. In these cases, since $T$ only depends on invariant quantities, we can iterate the contraction argument on any interval $[kT,(k+1)T]$, $k\in\Z$, hereby extending $\psi$ to a global mild solution in $C(\R;H^s(\T;\C))$, which satisfies the extra property of belonging to $X^{s,\frac12+\delta}_I$ for any time interval $I$ of size $T$ as above.

\medskip
In case (ii) ($\kappa=0$ and $s>\frac13$), let us discuss the unconditional uniqueness. Let $\psi_1,\psi_2\in C([0,T];H^s(\T;\C))$ be mild solutions of \eqref{NLS} with $\kappa=0$, sharing the same initial data, and let $\varphi=\psi_1-\psi_2$, which satisfies
\begin{align*}
    \varphi(t)=i\int_0^te^{-i(t-t')\dx^2}(u\cdot\mu)\varphi dt'.
\end{align*}
Then, 
\begin{align*}
    \|\varphi(T)\|_{L^2}^2&=-2\Im\int_0^T\Big\langle (u\cdot\mu)\varphi,i\int_0^te^{-i(t-t')\dx^2}\big\{(u\cdot \mu)\varphi\big\}dt'\Big\rangle_{L^2_x}dt\\
    &=-2\Im\int_0^T\big\langle (u\cdot\mu)\varphi,\varphi\big\rangle_{L^2_x}dt=0,
\end{align*}
provided that we can make sense of the first line, i.e. the space-time duality pairing. If we only know that the term in the right-hand side of the bracket belongs to $C([0,T];H^s(\T;\C))$, we need to check whether $(u\cdot\mu)\varphi\in L^1([0,T];H^{-s}(\T;\C))$. We have using Cauchy-Schwarz inequality
\begin{align*}
\|\mu_j\varphi\|_{H^{-s}}^2&=\sum_{k\in\Z}\jb{k}^{-2s}\Big|\sum_{k_1\in\Z}(\ft\mu_j)_{k-k_1}\ft\varphi_{k_1}\Big|^2\\
&\les \|\mu_j\|_{\F\ell^{s,\infty}}^2\|\varphi\|_{H^s}^2\sum_{k,k_1\in\Z}\jb{k}^{-2s}\jb{k-k_1}^{-2s}\jb{k_1}^{-2s},
\end{align*}
and the last sum is finite provided that $s>\frac13$: indeed, using symmetry,
\begin{align*}
\sum_{k,k_1\in\Z}\jb{k}^{-2s}\jb{k-k_1}^{-2s}\jb{k_1}^{-2s}\les \sum_{\substack{k,k_1\in\Z\\ |k|\sim |k-k_1|\gg |k_1|}}\jb{k}^{-4s}\jb{k_1}^{-2s}\les \sum_{k_1\in\Z}\jb{k_1}^{1-6s} <\infty,
\end{align*}
where in the second step we used that $s>\frac14$, while the last step holds if $s>\frac13$.

\medskip
It remains to discuss cases (iii) and (iv), namely, the unconditional uniqueness. Proceeding as before, we now have
\begin{align*}
    \|\varphi(T)\|_{L^2}^2&=-2\Im\int_0^T\Big\langle (u\cdot\mu)\varphi+\kappa\big(f(\psi_1)-f(\psi_2)\big),i\int_0^te^{-i(t-t')\dx^2}\big\{(u\cdot\mu)\varphi+\kappa \big(f(\psi_1)-f(\psi_2)\big)\big\}dt'\Big\rangle_{L^2_x}dt\\
    &=-2\Im\int_0^T\big\langle (u\cdot\mu)\varphi+\kappa\big(f(\psi_1)-f(\psi_2)\big),\varphi\big\rangle_{L^2_x}dt=-2\kappa\Im\int_0^T\big\langle f(\psi_1)-f(\psi_2),\varphi\big\rangle_{L^2_x}dt,
\end{align*}
where $f(z)=|z|^{2m}z$, and the cancellation at the last step is justified as above by making sense of the duality pairing provided as above that $s>\frac13$. We now need to estimate the remaining term. Note that in order for the space-time duality pairing ${\displaystyle \int_0^T\big\langle f(\psi_1)-f(\psi_2),\varphi\big\rangle_{L^2_x}dt}$ to make sense, we need $f(\psi_1)-f(\psi_2)\in L^1([0,T];H^{-s}(\T;\C))$. Proceeding as above, provided that $s\geq \frac12- \frac1{2m+2}$, we have by Sobolev inequality
\begin{align*}
\|f(\psi_1)-f(\psi_2)\|_{L^1_TH^{-s}}&\les T \|\varphi(|\psi_1|+|\psi_2|)^{2m}\|_{L^\infty_TL^{\frac{2m+2}{2m+1}}}\\ 
&\les T\|\varphi\|_{L^\infty_TL^{2m+2}}(\|\psi_1\|_{L^\infty_T L^{2m+2}}+\|\psi_2\|_{L^\infty_T L^{2m+2}})^{2m}\\
&\les T\|\varphi\|_{L^\infty_TH^s}(\|\psi_1\|_{L^\infty_T H^s}+\|\psi_2\|_{L^\infty_T H^s})^{2m}.
\end{align*}
Thus the duality bracket makes sense for $s\geq \frac12-\frac1{2m+2}$. However, to close the estimate in order to recover uniqueness, we need the stronger estimate
\begin{align*}
\|\varphi\|_{L^\infty_TL^2}^2&\les  \Big|\int_0^T\big\langle f(\psi_1)-f(\psi_2),\varphi\big\rangle_{L^2_x}dt\Big|\\
&\les T\|\varphi\|_{L^2}^2\big(\|\psi_1\|_{L^\infty_T L^\infty}+\|\psi_2\|_{L^\infty_T L^\infty}\big)^{2m}\\
&\les T\|\varphi\|_{L^\infty_TL^2}^2\big(\|\psi_1\|_{L^\infty_T H^s}+\|\psi_2\|_{L^\infty_T H^s}\big)^{2m},
\end{align*}
but now the last step holds provided that $s>\frac12$. Taking $\tau>0$ small enough depending on $\|\psi_1\|_{L^\infty_T H^s}+\|\psi_2\|_{L^\infty_T H^s}$, we thus find that $\varphi=0$ on $[0,\tau]$, and we can iterate the previous argument on intervals of size $\tau$ to conclude that $\varphi=0$ on $[0,T]$. This shows unconditional uniqueness in case (iii) and concludes the proof of Theorem~\ref{THM:LWP}.
\end{proof}

\section{Small-time global approximate controllability in low regularity spaces}\label{Sec:global-approx}
The aim of this section is to establish short-time global approximate controllability, which, combined with small-time local exact controllability, yields Theorem~\ref{THM:control_global_intro}.
To obtain this short-time global approximate result, we will refine the techniques of \cite{BP04-2} at two levels: first, we treat the equation in presence of a polynomial nonlinearity\footnote{In \cite{BCP}, the authors deal with \eqref{NLS} where the polynomial nonlinearity is replaced by a logarithmic one. As they stick to an $L^2$ framework, their techniques do not transfer to our setting.}; second, our results are proved in an $H^s$ setting instead of the $L^2$ setting considered in the aforementionned work, which is what enables us to deal with the nonlinearity.

\bigskip

 Although Theorem~\ref{THM:control_global_intro} is stated in the one-dimensional setting, global approximate controllability can still be ensured for the following nonlinear Schrödinger equation defined on the $N$-dimensional torus $\T^N$:
\begin{equation}\label{NLS-N}\begin{split}\begin{cases}
i\dt \psi -\Delta  \psi = (u(t)\cdot\mu(x))\psi  +\kappa |\psi|^{2m}\psi,~~&(t,x)\in\R\times\T^N,\\
\psi(t=0,x)=\psi_0(x),~~&x\in\T^N.
\end{cases}\end{split}\end{equation}

Here, $\kappa\in\R$, the operator $-\Delta$ is the self-adjoint Laplacian in $L^2(\T^N;\C)$ with domain $H^2(\T^N;\C)$, and $u(t)\cdot\mu(x)$ is defined by the vector-valued function $\mu\in C^\infty(\T^N;\R^d)$ and $u:[0,T]\rightarrow \R^d$. Let us define, for $j\leq N-1$, the vectors 
$$\mathfrak e_j=(\mathfrak e^1_j,...,\mathfrak e^N_j),\ \ \ \ \ \mathfrak e^k_j=\delta_{j,k};\ \ \ \ \ \  \mathfrak e_N=(1,...,1).$$

The arguments developped in \cite{BP04-2} are indeed partly based on saturation techniques, which were introduced in \cite{Vahagn} for \eqref{NLS-N}. We start by stating a slight adaptation of \cite[Theorem A]{Vahagn}.

\begin{proposition}\label{THM:contr-phase}
Let $s\in \R^*$ with $s>\frac{N}{2}$ and $\theta\in C^\infty(\T^N;\R)$. Then for every initial state $\psi_0 \in H^s(\T^N;\C)$ and $\eta > 0$, there exist $T \in [0,\eta]$ and a piecewise constant control $u$ defined in $(0,T)$ such that the corresponding solution $\psi(t)$ to \eqref{NLS-N} with data $\psi_0$ and control $u$ verifies
$$\|\psi(T) -  e^{i\theta}\psi_0\|_{H^s} < \eta.$$
\end{proposition}
\begin{proof}The proof closely follows that of \cite[Theorem A]{Vahagn}. The only distinction occurs when $s > \frac{N}{2}$ is not an integer and remains below the smallest integer strictly exceeding $\frac{N}{2}$. In this regime, one may use the fact that $H^s$ is an algebra, avoiding the need for integration by parts and interpolation between integer-order Sobolev spaces, as done in the cited prrof.
\end{proof}

A direct consequence of this proposition is the following corollary, which is used for Theorem~\ref{THM:control_global_intro_eigen}.

\begin{corollary}\label{COR:control_modes}
Let $s'\in (0,1/2)$. Then for every $n,m\in\Z^N$ and $\eta > 0$, there exist $T \in [0,\eta]$ and a piecewise constant control $u$ defined in $(0,T)$ such that the corresponding solution $\psi(t)$ to \eqref{NLS-N} with data $e^{i n\cdot x}$ and control $u$ verifies
$$\|\psi(T) -  e^{i m\cdot x}\|_{H^{s'}} < \eta.$$
\end{corollary}
\begin{proof}
The result follows from Proposition~\ref{THM:contr-phase} by considering $\psi_0=e^{in\cdot x}$. Then it is enough to take $\theta$ a smooth approximation of $x\mapsto (n-m)\cdot x$ in $H^{s'}(\T^N,\R)$ chosen such that $\|e^{i(n-m)\cdot x}-e^{i\theta}\|_{H^{s'}}$ is small. Note that this approximation is possible since $s'<1/2$.
\end{proof}

\begin{proposition}\label{THM:control_app}
Let $N\in\N^*$, $\kappa\in\R$. Consider the equation \eqref{NLS-N} in the presence of $d=2N+1$ controls $$\mu = \big( 1, \sin(\mathfrak e_1\cdot x),...,\sin(\mathfrak e_N\cdot x), \cos(\mathfrak e_1\cdot x),...,\cos(\mathfrak e_N\cdot x)\big).$$ 
%The equation \eqref{NLS-N} is small-time globally approximately controllable in $H^s(\T^N;\C)$ for every 
%Let $s>N/2$. % in the following sense.
Let $s>\frac{N}{2}$. For any $\eta>0$ and any
$$\psi_1,\psi_2\in X:=\{\psi= \rho e^{i\theta}\ :\ \rho \in C^\infty(\T^N;\R_+^*),\ \theta\in C^\infty(\T^N;\R),\ \|\rho\|_{L^2}=1\},$$
  there exist $T\in (0,\eta]$ and $u\in L^2((0,T);\R^{2N+1})$ such that there exists $\psi$ a unique mild solution to the equation \eqref{NLS-N}, satisfying
$$\psi(0)=\psi_1,\qquad  \|\psi(T)-\psi_2\|_{H^s}<\eta.$$
\end{proposition}

When $N=1$, we can actually improve the previous result to control a larger set of states.

\begin{proposition}\label{THM:control_app_1}
Let $N=1$ and $\kappa\in\R$. Consider the equation \eqref{NLS-N} in the presence of $d=3$ controls $$\mu = \big( 1, \sin(x),\cos(x)\big).$$ %The equation \eqref{NLS-N} is small-time globally approximately controllable in $H^s(\T^N;\C)$ for every $s>1/2$ in the following sense.
Let $s_1>1/2$ and $1/2<s<\min\{s_1,1\}.$ For any $\eta>0$ and any 
$$\psi_1,\psi_2\in \widetilde X:=\{\psi= \rho e^{i\theta}\ :\ \rho \in H^{s_1}(\T;\R_+^*),\ \theta\in H^{s_1}(\T;\R),\ \|\rho\|_{L^2}=1\},$$
  there exist $T\in (0,\eta]$ and $u\in L^2((0,T);\R^{3})$ such that there exists $\psi$ a unique mild solution to the equation \eqref{NLS-N}, satisfying
$$\psi(0)=\psi_1,\qquad  \|\psi(T)-\psi_2\|_{H^{s}}<\eta.$$
\end{proposition}

\subsection{Some intermediate results}

Proposition~\ref{THM:control_app} and Proposition~\ref{THM:control_app_1} are obtained by extending \cite[Theorem~3]{BP04-2}, which establishes the result in $L^2$, to higher regularity Sobolev spaces. The proof relies on a few technical results, which we provide in what follows.   Let us start with a standard local well-posedness result, which can be found for instance in \cite[Proposition 1.1]{Vahagn}.

\begin{lemma}\label{Lemma:WPHs}
 Consider the \eqref{NLS-N} and $s>\frac{N}{2}$. Then, for any $u \in L^2_{loc}(\R^+;\R^d)$, and any $m\in\N^*$, and any $\kappa\in\R$, then the equation \eqref{NLS-N} is locally well-posed in $H^s(\T^N;\C)$ in the following sense. There exist $T^*\in (0,+\infty]$ and a unique local mild solution $$\psi\in C^0([0,T^*); H^s(\T^N;\C))$$ to \eqref{NLS-N} with initial condition $\psi_0\in H^s(\T^N;\C)$, that is a solution to
the Duhamel formula: 
\begin{align*}
 \psi(t)=e^{- it\Delta }\psi_0-i\int_0^te^{ -i(t-t')\Delta }\big\{(u(t')\cdot\mu(x))\psi(t')+\kappa|\psi(t')|^{2m}\psi(t')\big\}dt',
\end{align*}
with $t\in [0,T^*)$ and satisfying $\|\psi(t)\|_{L^2}=\|\psi_0\|_{L^2}$ for any $t\in [0,T^*)$.
\end{lemma}
We can now introduce the following notation.
\begin{definition}
Let $s>\frac{N}{2}$, $\phi$ be a smooth function on $\T^N$ and $n\in\N^*$. Let $t>0$. We introduce, when it is defined, the  conjugated endpoint map
\[\begin{array}{rcl}
\mathcal R(t,\cdot,\cdot) :  H^s(\T^N)\times L^2_{loc}(\R_+) & \longrightarrow & H^s(\T^N) \\
(\psi_0,u) & \longmapsto & \psi(t)
\end{array}\]
where $\psi$ is solution to \eqref{NLS-N}.    
\end{definition}
For $f\in C^1(\T^N;\C^N)$, we denote by $T_{ f}$ the unitary transport generator 
\begin{equation}\label{DefOfTransports}
\begin{array}{rcl}
T_{ f}: H^1(\T^N;\C)& \longrightarrow & L^2(\T^N;\C) \\ \psi & \longmapsto & 
 \div(f\psi) +f\cdot \nabla \psi.
\end{array}
\end{equation}
We also denote the associated transport group by
\begin{align}\label{DefOfU}
U(t,f):=e^{t T_{f}}.
\end{align}
For $g\in C^2(\T^N;\C)$, we denote by
\begin{equation}\label{DefOfConjugateOperator}
    -\Delta_g =- e^{-ig}\Delta e^{ig} = -\Delta+|\nabla g|^2 -2i \nabla g\cdot \nabla - i\Delta g= -\Delta+|\nabla g|^2 -i T_\nabla g.
\end{equation}

\begin{lemma}\label{Lemma:B_n}
 Let $s>\frac{N}{2}$ and $\phi\in C^{s+3}(\T^N;\R)$. There exists $C>0$ such that for any $R>0$, there exists $\varepsilon_0>0$ such that for any $\psi_0\in\ H^s(\T^N;\C)$ with $\|\psi_0\|_{H^s}\leq R$, any $\varepsilon\in(0,\varepsilon_0)$, and any $n\in \N^*$, the solution $\mathcal R(t,e^{i\sqrt n \phi}\psi_0,0)$ is well-defined on the time interval $[0,\varepsilon n^{-\frac{3}{2}}]$, with the estimate
 \begin{align}\label{bound11}\|  e^{-i \varepsilon\frac{|\nabla\phi|^2}{\sqrt{n}}-i  \phi\sqrt{n}}\mathcal R(\varepsilon n^{-\frac{3}{2}},e^{i\sqrt n \phi}\psi_0,0) \|_{H^s} \leq e^{\varepsilon\frac{C}{n}}  \|\psi_0\|_{H^s}( 1+\varepsilon\frac{C}{{n}^{3/2}}\|\psi_0\|_{H^s}^{2m})  ,\quad \forall n\in\N^*.\end{align}
\end{lemma}
\begin{proof}
Let $T^*_n>0$ be the maximal time of existence of $\psi(t):=\mathcal R(t,e^{i\sqrt n \phi}\psi_0,0)$ the solution of \eqref{NLS-N}. Let $\tau>0$. We define the following time and space depending phase 
\[
H_\tau(t,x)=t\,\tau\, n|\nabla \phi(x)|^2+n^\frac{1}{2}\phi(x). \]
Let
\begin{equation}\label{DefOfZeta}
\zeta=e^{-i H_\tau}\mathcal R(\tau \cdot ,e^{i\sqrt n \phi}\psi_0,0).
\end{equation}
which satisfies %the following locally well-posed equation in $H^{s+2}(\T^N;\C)$
\[
\begin{cases}
i  \partial_t  \zeta
-\tau \Delta_{H_\tau} \zeta -(\partial_t H_\tau )\zeta =\tau \kappa |\zeta|^{2m}\zeta & \text{ in } (0,\tau^{-1}T^*_n)\times \T^N, \\
\zeta_{|_{t=0}} = \psi_0& \text{ in } \T^N.
\end{cases} \]
Expanding the first equation by using \eqref{DefOfTransports} and \eqref{DefOfConjugateOperator}, we obtain
\begin{multline*}
i\partial_t  \zeta-\tau\Delta \zeta + 4t^2n^2\tau^3|\nabla^2\phi \nabla \phi|^2\zeta - itn\tau^2 T_{\nabla |\nabla \phi|^2} \zeta -in^\frac{1}{2}\tau T_{ \nabla \phi}\zeta \\
+2t n^{\frac{3}{2}}\tau^2 (\nabla\phi\cdot\nabla|\nabla\phi|^2)\zeta
= \tau \kappa |\zeta|^{2m}\zeta.
\end{multline*}
%Let $\alpha\in \N^N$ such that $|\alpha|\leq s$. We multiply by $\partial_x^{2\alpha}\overline \zeta$ and integrate by parts to obtain 
We proceed with an estimation of the energy in $[0,2]$ and 
\begin{multline*}
\frac{d}{2dt}\|   \zeta \|^2_{H^s} + 4t^2n^2\tau^3\Im\langle  |\nabla^2\phi\nabla\phi|^2 \zeta, \zeta \rangle_{H^s} -tn\tau^2 \Re \Big\langle  \big[|\Delta|^{\frac{s}{2}}, T_{\nabla |\nabla \phi|^2}\big]\zeta, |\Delta|^{\frac{s}{2}}  \zeta \Big\rangle_{L^2} \\
-n^{\frac{1}{2}}\tau\Re\Big\langle \big[|\Delta|^{\frac{s}{2}}, T_{ \nabla \phi}\big] \zeta, |\Delta|^{\frac{s}{2}} \zeta \Big\rangle_{L^2} + 4 t n^{\frac{3}{2}}\tau^2\Im \langle  \nabla \phi \cdot \nabla^2\phi\nabla\phi  \zeta,  \zeta \rangle_{H^s} =\kappa\tau \Im\langle |\zeta|^{2m}\zeta,  \zeta \rangle_{H^s.}
\end{multline*} 
In the last identity, the commutator terms are controlled by Kato--Ponce type estimates, yielding a bound of the scalar product in $L^2$ in terms of the $C^{s+3}$ norm of $\phi$, and the $H^s$ norm of $\zeta$, with no loss of regularity in $\zeta$.
Hence, for any $t \in \big[0,\min\{2,T^*_n\tau^{-1}\}\big)$,
\begin{align*}
\frac{d}{dt} \| \zeta\|_{H^s}^2 \lesssim   (n^2\tau^3 + n \tau^2 + n^\frac{1}{2}\tau + n^\frac{3}{2}\tau^2) \|\zeta\|_{H^s}^2 +   \tau \|\zeta\|_{H^s}^{2m+2},
\end{align*}
Using the Gr\"onwall inequality, there exists $C>0$ such that, for any $t \in [0,\min\{2,T^*_n\tau^{-1}\})$,
\begin{equation*}
  \| \zeta(t)\|_{H^s}\leq e^{Ct(n^2\tau^3 + n \tau^2 + n^\frac{1}{2}\tau + n^\frac{3}{2}\tau^2)}  \big(\|\psi_0\|_{H^s}+C \tau\int_0^t \|\zeta(r)\|_{H^s}^{2m+1}dr\big) .
\end{equation*}
For $\varepsilon\in(0,1)$, we fix the scaling $\tau=\varepsilon n^{-\frac{3}{2}}$, and %by density a density argument, 
we have the existence of $C>0$ such that%, for every $\psi_0\in H^s(\T;\C)$,
\begin{align}\label{gronw-1}
  \| \zeta(t)\|_{H^s}\leq    e^{t\varepsilon\frac{ C}{n}}( \|\psi_0\|_{H^s}+\varepsilon  \frac{ C}{n^{3/2}}\int_0^t \|\zeta(r)\|_{H^s}^{2m+1}dr) .
\end{align}
Now, it is possible to notice that, there exists $\varepsilon\in(0,1)$ such that, for every $n\in\N^*$, we have 
\begin{align*}
\tau_n=\sup\left\{t_0< \min \{2,\varepsilon^{-1} T^*_n n^{\frac{3}{2}}\},\;\|\zeta( t)\|_{H^s}<2\|\psi_0\|_{H^s}\;,\forall t\in [0,t_0]\right\} >\frac{3}{2}.\end{align*} 
Indeed, assume this is not true. Then, for any $\varepsilon\in(0,1)$ there exists $n_\varepsilon\in\N^*$ such that  
\begin{align}\label{property-1}
\tau_{n_\varepsilon}=\sup\left\{t_0< \min \{2,\varepsilon^{-1} T^*_{n_\varepsilon }n_\varepsilon^{\frac{3}{2}}\},\;\|\zeta( t)\|_{H^s}<2\|\psi_0\|_{H^s}\;,\forall t\in [0,t_0]\right\} \leq\frac{3}{2}.\end{align} 
%an increasing sequence $(n_k)_{k\in\N^*}\subset\N^*$ such that $\lim_{k\to+\infty}n_k=+\infty$ and such that $\tau_{n_k}<1$ for every $k\in\N^*$.
Hence, \eqref{gronw-1} implies  
\begin{align*}
2\|\psi_0\|_{H^s} =  \| \zeta(\tau_{n_\varepsilon})\|_{H^s} \leq &    e^{\tau_{n_\varepsilon}\varepsilon\frac{ 1}{n_\varepsilon}}  ( \|\psi_0\|_{H^s}+\varepsilon\frac{C}{{n_\varepsilon}^{3/2}}\int_0^{\tau_{n_\varepsilon}} \|\zeta(r)\|_{H^s}^{2m+1}dr) \\
\leq & e^{\varepsilon  C}   ( \|\psi_0\|_{H^s}+\varepsilon\frac{C}{{n_\varepsilon}^{3/2}} (2\|\psi_0\|_{H^s})^{2m+1} ) .
\end{align*}
But we obtain a contradiction for $\varepsilon$ sufficiently small only depending on $ \|\psi_0\|_{H^s}$. % < C^{-1}\min \big\{\ln(\tfrac{3}{2}), \tfrac{1}{8}(2\|\psi_0\|_{H^s})^{-2m}\big\} .$$
% \begin{align*}
% 2\|\psi_0\|_{H^s} <2\|\psi_0\|_{H^s}.
% \end{align*}
As a result, the function $\zeta$ is well-defined in $[0,1]\subset[0,3/2]\subset[0,\tau_n]$. Hence, $T^*_n>\varepsilon n^{-\frac{3}{2}}$ and then the first claim is proved. The final bound \eqref{bound11} follows by evaluating \eqref{gronw-1} at $t=1$ along with \eqref{property-1} since%. Indeed, by taking $\varepsilon$ sufficiently small (only depending on $ \|\psi_0\|_{H^s}$), %there exists $C_1>0$ such that
\begin{align*}
  \| \zeta(1)\|_{H^s}&\leq e^{\varepsilon \frac{C}{n}}   ( \|\psi_0\|_{H^s}+\varepsilon\frac{C}{{n}^{3/2}} (2\|\psi_0\|_{H^s})^{2m+1}).%\leq  e^{\varepsilon\frac{ C}{n}}  R( 1+\varepsilon\frac{C}{{n}^{3/2}}R^{2m}).%\leq e^{\varepsilon\frac{2 C_1}{n}} (R+1)^{2m+1}  .
\end{align*}
\end{proof}

 We present now the following property which is ensured by the proof of Lemma \ref{Lemma:B_n}.
\begin{corollary}\label{Coro:B_n}
 Let $s>\frac{N}{2}$ and $\phi\in C^{s+3}(\T^N;\R)$. There exists $C>0$ such that for any $R>0$, there exists $\varepsilon_0>0$ such that for any $\psi_0\in\ H^s(\T^N;\C)$ with $\|\psi_0\|_{H^s}\leq R$, any $\varepsilon\in(0,\varepsilon_0)$, and any $n\in \N^*$, we have
 \begin{align}\label{bound12}\|  e^{-i(t\,n\,|\nabla\phi|^2 + n^{1/2}\phi)}\mathcal R(t,e^{i\sqrt n \phi}\psi_0,0) \|_{H^s} \leq e^{\varepsilon\frac{C}{n}}  R( 1+\varepsilon\frac{C}{{n}^{3/2}}R^{2m})  ,\quad t\in [0,\varepsilon n^{-\frac{3}{2}}].\end{align}
\end{corollary}

From Lemma \ref{Lemma:B_n}, the following map is well defined.

\begin{definition} 
Let $s>\frac{N}{2}$, $\phi\in C^{s+3}(\T^N;\R)$, and $n\in\N^*$. Consider $R>0$ and $\varepsilon_0>0$ verifying Lemma \ref{Lemma:B_n}. We define
\begin{equation}\label{def:Bn} 
\begin{array}{rcl}
B_n : \big\{\psi_0\in H^s,\;\|\psi_0\|_{H^s} \leq R\big\}\times (0,\varepsilon_0] & \longrightarrow &  H^s(\T^N;\C)\\
(\psi_0,\varepsilon) & \longmapsto & e^{-i\varepsilon \frac{|\nabla\phi|^2}{\sqrt{n}}-i  \phi\sqrt{n}}\mathcal R(\varepsilon n^{-\frac{3}{2}},e^{i \sqrt{n}\phi}\psi_0,0).
\end{array}
\end{equation}
%Note that this is well-defined for any $0< \varepsilon\leq\varepsilon_0 $ with $\varepsilon_0$ defined by Lemma \ref{Lemma:B_n}.
\end{definition}

Recall \eqref{DefOfU}. We now present the second technical lemma leading to the proof of Proposition \ref{THM:control_app}.

\begin{lemma}\label{Lemma:B_n_lip}
 Let $s>\frac{N}{2}$ and $\phi\in C^{s+3}(\T^N;\R)$. There exists $C>0$ such that for any $R>0$, there exists $\varepsilon_1>0$ such that for any $\psi_1,\psi_2\in\ H^s(\T^N;\C)$ with $\|\psi_1\|_{H^s},\|\psi_2\|_{H^s}\leq R$, any $0<\varepsilon\leq \varepsilon_1$, and any $n\in \N^*$, we have
 \begin{align}\label{bound13}\|  B_n(\psi_1,\varepsilon) -B_n(\psi_2,\varepsilon) \|_{H^s} \leq e^{\varepsilon\frac{ C}{n}\big(\tfrac{ R^{2m}}{\sqrt{n}}+1\big) }  \|\psi_1-\psi_2\|_{H^s}  ,\quad \forall n\in\N^*.\end{align}
\end{lemma}
\begin{proof}
We proceed similarly to the proof of Lemma \ref{Lemma:B_n}. We define
\begin{equation}\label{DefOfZeta2}
\zeta_1=B_n(\psi_1,\varepsilon),\quad\zeta_2=B_n(\psi_2,\varepsilon),\quad\zeta=\zeta_1-\zeta_2.
\end{equation}
The function $\zeta$ satisfies the following equation which we consider for $t\in [0,1]$
\begin{multline*}
i\partial_t  \zeta-\tau\Delta \zeta + 4t^2n^2\tau^3|\nabla^2\phi \nabla \phi|^2\zeta - itn\tau^2 T_{\nabla |\nabla \phi|^2} \zeta -in^\frac{1}{2}\tau T_{ \nabla \phi}\zeta \\
+2n^{\frac{3}{2}}\tau^2 (\nabla\phi\cdot\nabla|\nabla\phi|^2)\zeta
= \tau \kappa |\zeta_1|^{2m}  \zeta_1-\tau \kappa |\zeta_2|^{2m} \zeta_2.
\end{multline*}
%Let $\alpha\in \N^N$ such that $|\alpha|\leq s$. We multiply by $\partial_x^{2\alpha}\overline \zeta$ and integrate by parts to obtain 
Thanks Corollary \ref{Coro:B_n}, there exists $C>0$ such that, for $\varepsilon_1$ small and any $0<\varepsilon\leq \varepsilon_1$, we have $$\|\zeta_1\|_{H^s}\leq e^{\varepsilon\frac{C}{n}} R( 1+\varepsilon\frac{C}{{n}^{3/2}}R^{2m}),\quad\|\zeta_2\|_{H^s}\leq e^{\varepsilon\frac{C}{n}} R( 1+\varepsilon\frac{C}{{n}^{3/2}}R^{2m}),$$
and then, for $\varepsilon_1$ sufficiently small and $t\in [0,1]$,
\begin{align*}
\| |\zeta_1|^{2m}  \zeta_1- |\zeta_2|^{2m} \zeta_2\|_{H^s}\leq C e^{\varepsilon\frac{C_1}{n}(2m)} R^{2m}( 1+\varepsilon\frac{C}{{n}^{3/2}}R^{2m})^{2m}\|\zeta\|_{H^s}\lesssim R^{2m}   \|\zeta\|_{H^s}.
\end{align*}
We proceed with an estimation of the energy as in the mentioned proof with the scaling $\tau=\varepsilon n^{-\frac{3}{2}}$. For any $t \in \big[0,1]$,
\begin{align*}
\frac{d}{dt} \| \zeta\|_{H^s}^2 &\lesssim \frac{\varepsilon }{n} \|\zeta\|_{H^s}^2 + \frac{\varepsilon }{n^{3/2}}R^{2m}   \|\zeta\|_{H^s}^2\lesssim \frac{ \varepsilon}{n} \Big(1 + \frac{ R^{2m}}{\sqrt{n}}\Big)\|\zeta\|_{H^s}^2.%,
\end{align*}
Using the Gr\"onwall inequality, we conclude the proof. % since for any $t \in [0,1]$, we have
%\begin{align*}
%  \| \zeta(t)\|_{H^s}\leq    e^{t\varepsilon\frac{ C}{n}}( \|\zeta(0)\|_{H^s}+\varepsilon  \frac{ C}{n^{3/2}} R^{2m}\int_0^t \|\zeta(r)\|_{H^s}dr) .
%\end{align*}
%We observe that thanks Corollary \ref{Coro:B_n}, for $\varepsilon_1$ sufficiently small, there exists $C>0$ such that
%\begin{align*}
%  \| \zeta(t)\|_{H^s}\leq    e^{t\varepsilon\frac{ C}{n}}( \|\zeta(0)\|_{H^s}+\varepsilon  \frac{ C}{n^{3/2}}\int_0^t \|\zeta(r)\|_{H^s}^{2m+1}dr+\varepsilon  \frac{ C}{n^{3/2}}  e^{\varepsilon\frac{C}{n}(2m+1)} R^{(2m+1)^2}) .
%\end{align*}
%As in the proof of Lemma \ref{Lemma:B_n}, up to reducing $\epsilon_1,$ there exists $C>0$ such that
%\begin{align*}
%  \| \zeta(t)\|_{H^s}\leq    e^{t\varepsilon\frac{ C}{n}}( \|\zeta(0)\|_{H^s}+\varepsilon\frac{C}{{n}^{3/2}} (2\|\zeta_0\|_{H^s})^{2m+1}+\varepsilon  \frac{ C}{n^{3/2}}  e^{\varepsilon\frac{C}{n}(2m+1)} R^{(2m+1)^2}) .
%\end{align*}
%Now, we consider $\varepsilon_1$ sufficiently small so that, for any $0<\varepsilon\leq \varepsilon_1$, we have
%$$C \varepsilon ((2\|\zeta_0\|_{H^s})^{2m+1}+ e^{\varepsilon\frac{C}{n}(2m+1)} R^{(2m+1)^2})<\|\zeta(0)\|_{H^s}.$$
%Now, up to reducing $\epsilon_1,$ the proof is concluded since
%\begin{align*}
 % \| \zeta(1)\|_{H^s}\leq    e^{\varepsilon\frac{ C}{2n} \big(1 + e^{\varepsilon\frac{C}{n}(2m)} R^{(2m+1)(2m)}\big)}%( (1+\frac{ 1}{n^{3/2}})\|\zeta(0)\|_{H^s}) \leq    e^{\frac{ 2}{n}}
 % \|\zeta(0)\|_{H^s} .
%\end{align*}
\end{proof}

\begin{lemma}\label{Lemma:U_n}
 Let $s>\frac{N}{2}$ and $\phi\in C^{s+4}(\T^N;\R)$. Then, there exists $C>0$ such that for any $R>0$ there exists $\varepsilon_1>0$ such that for any $\psi_0\in H^{s+2}(\T^N;\C)$ with $\|\psi_0\|_{H^{s+2}}\leq R $ and $\|\psi_0\|_{L^2}=1$, any $0<\varepsilon\leq \varepsilon_1$, and any $n\in\N^*$, we have
\begin{align*}
  \big\| B_n(\psi_0,\varepsilon)- U\big(\tfrac{\varepsilon}{n}, \nabla \phi\big)\psi_0 \big\|_{H^s}& \leq \frac{C%\varepsilon
  }{n}.%\max\Big\{   \|\psi_0\|_{H^{s+2}}^{2m+1}.\|\psi_0\|_{H^{s+2}}.%\Big\}.
\end{align*}\end{lemma}
\begin{proof}
%Consider $\psi_0\in H^{s+2}(\T^N;\C)$ such that $\|\psi_0\|_{L^2}=1$. 
Consider $R>0$ and $\varepsilon_0>0$ verifying Lemma \ref{Lemma:B_n}. Let $\zeta$ be defined as \eqref{DefOfZeta} in the proof of Lemma \ref{Lemma:B_n}, with $\tau=\varepsilon n^{-3/2}$, $0<\varepsilon\leq \varepsilon_0$. Note that $\zeta= B_n(\psi_0,\varepsilon)$. From Lemma \ref{Lemma:B_n}, the following functions are defined in $C^0([0,\tfrac{3}{2}];H^s)$:
$$\varphi(t)=\zeta(t)-\xi(t),\quad\xi(t):=U(\varepsilon\tfrac{t}{n}, \nabla \phi)\psi_0=e^{\frac{\varepsilon t}{n}T_{\nabla\phi}}\psi_0.$$ The purpose is to estimate $\|\varphi(t=1)\|_{H^s}$. As $ \zeta$ solves %the following locally well-posed equation in $H^{s+2}(\T^N;\C)$
\[
\begin{cases}
i  \partial_t  \zeta
=\tau \Delta_{H_\tau} \zeta +(\partial_tH_\tau )\zeta +\tau \kappa |\zeta|^{2m}\zeta & \text{ in } (0,\tfrac{3}{2})\times \T^N \\
\zeta_{|_{t=0}} = \psi_0& \text{ in } \T^N, 
\end{cases} \]
we deduce that $\varphi$ satisfies 
\begin{equation}\label{EDP}
   \begin{cases} i\partial_t\varphi=\tau\Delta_{H_\tau}( \varphi+\xi)+(\partial_t H_\tau)(\varphi +\xi)+\kappa\tau |\varphi+\xi|^{2m}(\varphi+\xi) -\varepsilon 
in^{-1}T_{\nabla \phi}\xi & \text{ in } (0,\tfrac{3}{2})\times \T^N \\
   \varphi_{|_{t=0}} = 0 & \text{ in } \T^N.
   \end{cases}
\end{equation}
% \begin{equation}\begin{split}\label{EDP}
% \partial_t \varphi
% =&
%  \frac{1}{n}T_{\nabla \phi}(\varphi+w(t))
% -\frac{i}{n^{3/2}}\Delta(\varphi+w(t))-\frac{it}{n^{3/2}}\nabla\phi\cdot\nabla M\,(\varphi+w(t))
% +\frac{t}{2 n^{2}}T_{\nabla M}(\varphi+w(t))\\
% &+\frac{it^2}{4n^{5/2}}|\nabla M|^2(\varphi+w(t)) - \frac{1}{n}T_{\nabla\phi}w-\frac{i}{n^{3/2}} \kappa |\varphi+w(t)|^{2m} (\varphi+w(t))\\
% &=
%  \frac{1}{n}T_{\nabla \phi} \varphi -\frac{i}{n^{3/2}}\Delta(\varphi+w)-\frac{it}{n^{3/2}}\nabla\phi\cdot\nabla M\,(\varphi+w)+\frac{t}{2 n^{2}}T_{\nabla M}(\varphi+w)\\
%  &
% +\frac{it^2}{4n^{5/2}}|\nabla M|^2(\varphi+w) -\frac{i}{n^{3/2}} \kappa |\varphi+w(t)|^{2m} (\varphi+w(t)).
% \end{split}\end{equation}
As in the proof of Lemma \ref{Lemma:B_n}, we shall perform an energy estimate for $t\in[0,\tfrac{3}{2}]$. We emphasize that for every $r\in [s,s+2]$, there exists $C_r>0$ such that 
\[\|\xi\|_{H^{r}}=\|e^{-\frac{t}{n}T_{\nabla\phi}}\psi_0\|_{H^{r}}\leq C_r\|\psi_0\|_{H^{r}}\] 
since $\phi\in  C^{s+4}(\T^N;\R)$. %Let $\alpha\in (\N^*)^N$ such that $|\alpha|\leq s$. 
%By using the Young inequality, we make an estimation of the energy in $[0,1]$.
We recall that $1=\|\psi_0\|_{L^2}\leq \|\psi_0\|_{H^{s+2}}$ and then
\begin{align*}
\frac{d}{dt} \|\varphi\|_{H^s}^2
=-2 \Im \langle \varphi,i \partial_t\varphi \rangle_{H^s} \lesssim \frac{\varepsilon}{n}\|\varphi\|_{H^s}^2+\frac{\varepsilon}{n^2} \|\psi_0\|_{H^{s+2}}^{4m+2} +\frac{\varepsilon}{ n^{3/2}}\|\varphi\|_{H^s}^{2m+2}.
\end{align*} 
Using the Gr\"onwall's inequality for $t\in [0,\tfrac{3}{2}]$, there exists $C>0$ such that
\begin{align}\|\varphi(t)\|_{H^s}^2 \leq \frac{C\varepsilon}{ n^{3/2}} \Big(\frac{1}{n^{1/2}} \|\psi_0\|_{H^{s+2}}^{4m+2}    + \int_0^t \|\varphi(r)\|_{H^s}^{2m+2}dr\Big)\label{gronw-2}
\end{align} 
Now, we want to prove the existence of $0<\varepsilon_1\leq \varepsilon_0$ such that, for any $n\in\N^*$ and any $0<\varepsilon\leq\varepsilon_1$, we have 
\begin{align}\label{property}\tau_n=\sup\{t_0\in [0,\tfrac{3}{2}],\;\|\varphi( t)\|_{H^s}^2<\tfrac{ \|\psi_0\|_{H^{s+2}}^{2}}{ n},\; t\in [0,t_0]\}> \frac{5}{4}.\end{align} 
% First, notice that for $n$ sufficiently large, $\tau_n>0$ thanks to \eqref{gronw-2}, since for $n$ large,
% \begin{align*}\|\varphi(t=0)\|_{H^s}^2&\leq \frac{C}{ n^{2}}   \|\psi_0\|_{H^{s+2}}^{4m+2} <\frac{ \|\psi_0\|_{H^{s+2}}^2 }{ n} .
%\end{align*}
Assume that \eqref{property} is not verified. Then, for any $\varepsilon_1>0$ there exists $n_{\varepsilon_1}\in\N^*$ and $0<\varepsilon\leq \varepsilon_1$ such that  
\[
\tau_{n_{\varepsilon_1}}=\sup\{t\in [0,\tfrac{3}{2}],\;\|\varphi(\tilde t)\|_{H^s}^2<\tfrac{ \|\psi_0\|_{H^{s+2}}^{2}}{n_{\varepsilon_1}},\;\tilde t\in [0,t]\}\leq \frac{5}{4}.
\] Now, \eqref{gronw-2} implies  
\begin{align*}
\frac{ \|\psi_0\|_{H^{s+2}}^{2} }{ n_{\varepsilon_1}}=  \| \varphi(\tau_{n_{\varepsilon_1}})\|_{H^s}^2  \leq &\frac{C\varepsilon}{ {n_{\varepsilon_1}}^{3/2}} \Big(\frac{1}{{n_{\varepsilon_1}}^\frac{1}{2}} \|\psi_0\|_{H^{s+2}}^{4m+2}    + \int_0^{\tau_{n_{\varepsilon_1}}} \|\varphi(r)\|_{H^s}^{2m+2}dr\Big)\\
\leq &  \frac{C\varepsilon}{ {n_{\varepsilon_1}}^{3/2}} \Big(\frac{1}{{n_{\varepsilon_1}}^{1/2}} \|\psi_0\|_{H^{s+2}}^{4m+2}    + \frac{\|\psi_0\|_{H^{s+2}}^{2m+2}}{n_{\varepsilon_1}^{m+1}}\Big).
\end{align*}
For $\varepsilon_1>0$ sufficiently small such that we have 
$$\frac{C\varepsilon_1}{ {n_{\varepsilon_1}}^{1/2}}   \Bigg(\frac{\|\psi_0\|_{H^{s+2}}^{4m}}{{n_{\varepsilon_1}}^\frac{1}{2}} + \frac{\|\psi_0\|_{H^{s+2}}^{2m}}{n_{\varepsilon_1}^{m+1}}\Bigg)  <1,$$ we obtain a contradiction.
Hence, we conclude the proof by taking $\varepsilon$ sufficiently small since
\begin{align*}\|\varphi(t=1)\|_{H^s}^2&\lesssim  \varepsilon\|\psi_0\|_{H^{s+2}}^{4m+2} \Big(\frac{1}{ n^{2}}   + \frac{1}{n^{m+1}}\frac{1}{ n^{3/2}} \Big)\lesssim \frac{1}{ n^{2}}   .
\end{align*} 
\end{proof}

As a consequence of the proof of Lemma \ref{Lemma:U_n} we have the following corollary. 

\begin{corollary}\label{Coro:U_n}
 Let $s>\frac{N}{2}$ and $\phi\in C^{s+4}(\T^N;\R)$. Then, there exists $C>0$ such that for any $R>0$ there exists $\varepsilon_1>0$ such that for any $\psi_0\in H^{s+2}(\T^N;\C)$ with $\|\psi_0\|_{H^{s+2}}\leq R $ and $\|\psi_0\|_{L^2}=1$, any $0<\varepsilon\leq \varepsilon_1$, and any $n\in\N^*$, we have

\begin{align*}
  \Big\| e^{-i(tn|\nabla\phi|^2 + n^{1/2}\phi)}\mathcal R(t,e^{i  \phi\sqrt{n}}\psi_0,0)- U\big(\tfrac{\varepsilon}{n}, \nabla \phi\big)\psi_0 \Big\|_{H^s}& \leq \frac{C%\varepsilon
  }{n}%\max\Big\{\|\psi_0\|_{H^{s+2}}%^{2m+1}
 ,\quad \forall t\in [0,\varepsilon n^{-\frac{3}{2}}].%,\|\psi_0\|_{H^{s+2}}\Big\}.
\end{align*}\end{corollary}

\begin{lemma}\label{Lemma:U_n_bis}
 Let $s>\frac{N}{2}$ and $\phi\in C^{s+4}(\T^N;\R)$. Then, there exists $C>0$ such that for any $R>0$ there exists $\varepsilon_1>0$ such that for any $\psi_0\in H^{s+2}(\T^N;\C)$ with $\|\psi_0\|_{H^{s+2}}\leq R $ and $\|\psi_0\|_{L^2}=1$, any $0<\varepsilon\leq \varepsilon_1$, and any $n\in\N^*$, we have
\begin{align*}
 \Big\| B_n(\psi_0,\varepsilon)- U\big(\tfrac{\varepsilon}{n}, \nabla \phi\big)\psi_0 \Big\|_{L^2}& \leq  
\frac{C }{n^{3/2}}  .%\|  \psi_0   \|_{ H^{s+2} } .
\end{align*}
\end{lemma}
\begin{proof}
The strategy of the proof is the same as the one of Lemma \ref{Lemma:U_n}. We repeat the first part of the proof until \eqref{EDP} with the only difference that this time we make an estimation of the energy in $L^2(\T^N;\C)$. 
%uniformly in $n$ by Corollary~\ref{Coro:U_n}
By using the Young inequality, we make an estimation of the energy in $[0,1]$ 
\begin{align*}
\frac{d}{dt}  \|\varphi\|_{L^2}^2
&\lesssim \varepsilon \frac{\|  \psi_0\|_{H^{2}}^2+ \|  \psi_0\|_{L^2}^{4m+2}}{n^{3}}+\varepsilon\|\varphi\|_{L^2}^2+\frac{\varepsilon}{ n^{3/2}}\big(\|\zeta\|_{H^s}^{2m}+\|\xi\|_{H^s}^{2m}\big)\|\varphi\|_{L^2}^{2}\\
&\lesssim   \varepsilon\frac{ \|  \psi_0\|_{H^{s+2}}^{(2m+1)(2m+2)}}{n^{3}}+\varepsilon\|\varphi\|_{L^2}^2,
\end{align*} 
where we use $\|\xi\|_{H^s}\leq C\|\psi_0\|_{H^s}$ since $U$ is unitary on $H^s$, and $\|\zeta\|_{H^s}\leq C$ uniformly in $n$ by Corollary~\ref{Coro:U_n}. The proof is then concluded as the proof of Lemma~\ref{Lemma:U_n} with $\varepsilon$ sufficiently small. Indeed, by the Gronwall inequality,
$$\| \varphi(t=1)\|_{L^2}^2\lesssim \frac{\varepsilon}{n^{3}} \|  \psi_0\|_{H^{s+2}}^{(2m+1)(2m+2)} \lesssim \frac{1}{n^{3}} . $$

\end{proof}

Before moving on with the next lemma, we need to introduce the following notation.

\begin{definition} 
%Let $s>\frac{N}{2}$, $\phi\in C^{s+3}(\T^N;\R)$, and $n\in\N^*$. Consider $R>0$ and $\varepsilon_0>0$ verifying Lemma \ref{Lemma:B_n}. We define by iteration, for every $m\geq 2$,
%\begin{equation*}
%\begin{array}{rcl}
%B_n^m : \big\{\psi_0\in H^s,\;\|\psi_0\|_{H^s} \leq R\big\}\times (0,\varepsilon_0] & \longrightarrow &  H^s(\T^N;\C)\\
%(\psi_0,\varepsilon) & \longmapsto & B_n(B_n^{m-1}(\psi_0,\varepsilon),\varepsilon),
%\end{array}
%\end{equation*}
%with the convention $B_n^1(\psi_0,\varepsilon)=B_n(\psi_0,\varepsilon)$.
Consider $R>0$ and $\varepsilon_0>0$ verifying Lemma \ref{Lemma:B_n}. The iterates $B_n^m$ are well-defined on $\big\{\psi_0\in H^s(\T^N;\C),\;\|\psi_0\|_{H^s} \leq e^{2C}R\big\}\times (0,\varepsilon_0]$ for all $m\leq n/\varepsilon$, as shown in the proof of Lemma~\ref{Lemma:B_n-U_n-pre}. We define by iteration, for every $m\geq 2$,
\begin{equation*}
\begin{array}{rcl}
B_n^m : \big\{\psi_0\in H^s,\;\|\psi_0\|_{H^s} \leq R\big\}\times (0,\varepsilon_0] & \longrightarrow & H^s(\T^N;\C)\\
(\psi_0,\varepsilon) & \longmapsto & B_n(B_n^{m-1}(\psi_0,\varepsilon),\varepsilon),
\end{array}
\end{equation*}
with the convention $B_n^1(\psi_0,\varepsilon)=B_n(\psi_0,\varepsilon)$.
\end{definition}

\begin{lemma}\label{Lemma:B_n-U_n-pre}
Let $s>\tfrac{N}{2}$, $s_1\in\big(\tfrac{N}{2},s\big)$, and $\phi\in C^{s+5}(\T^N;\R)$. For any $\eta>0$ and $R>0$, there exists $\varepsilon_1>0$ such that for any $\psi_0\in H^{s+2}(\T^N;\C)$ with $\|\psi_0\|_{H^{s+2}}\leq R $ and $\|\psi_0\|_{L^2}=1 $, any $0<\varepsilon\leq \varepsilon_1$ such that $\varepsilon^{-1}\in\N^*$, there exists $\tilde n$ such that for every $n\geq\tilde n$, we have
% Let $s\in \N^*$ with $s>\frac{N}{2}$ and $\phi\in C^{s+4}(\T^N;\R)$. Consider $\psi\in H^{s}(\T^N;\C)$ such that $\|\psi_0\|_{L^2}=1$. Then, for every $\eta>0$, there exists $\tilde n\in\N^*$, such that, for every $n>\tilde n$, we have 
$$ \| B_n^{\tfrac{n}{\varepsilon}}(\psi_0,\varepsilon)- U(1, \nabla\phi)\psi_0 \|_{H^{s_1}}<\eta.$$
\end{lemma}
\begin{proof}
%Let us consider at first $\psi_0\in H^{s+3}(\T^N;\C)$ such that $\|\psi_0\|_{L^2}=1$. We generalize the result for $\psi_0\in H^{s}(\T^N;\C)$ in the last part of the proof. Thanks to the Lemma \ref{Lemma:B_n}, for $n$ large, we have% and \ref{Lemma:B_n-U_n}, we have
First, consider $C>0$ the constant defined in Lemma \ref{Lemma:B_n}.  We take $\tilde n$ sufficiently large such that, for every $n\geq \tilde n,$ we have
$$\frac{1}{{n}^{1/2}}(e^{2 C} R)^{2m}\leq 1.$$
Now, we use Lemma \ref{Lemma:B_n} with respect to the radius $e^{2 C} R$. There exists $\epsilon_1$ such that, for every $0<\epsilon\leq \epsilon_1$, thanks to the choice of $n$ verifying the previous inequality, we have
\begin{align*}\| B_n(\psi_0,\varepsilon)\|_{s+2}&\leq e^{\varepsilon\frac{C}{n}}  R( 1+\varepsilon\frac{C}{{n}^{3/2}}R^{2m}) \leq e^{\varepsilon\frac{C}{n}}  R( 1+\varepsilon\frac{C}{{n}^{3/2}}(e^{2C} R)^{2m}) \\
&\leq e^{\varepsilon\frac{C}{n}}  R( 1+\varepsilon\frac{C}{n})\leq e^{2\varepsilon\frac{C}{n}} R,   \end{align*}
\begin{align*}\| B_n(B_n(\psi_0,\varepsilon),\varepsilon)\|_{s+2}&\leq e^{\varepsilon\frac{C}{n}}  (e^{2\varepsilon\frac{C}{n}} R)( 1+\varepsilon\frac{C}{{n}^{3/2}}(e^{2\varepsilon\frac{C}{n}} R)^{2m}) \leq e^{\varepsilon\frac{C}{n}}  (e^{2\varepsilon\frac{C}{n}} R)( 1+\varepsilon\frac{C}{{n}^{3/2}}(e^{2C} R)^{2m})\\
&\leq e^{\varepsilon\frac{C}{n}}  (e^{2\varepsilon\frac{C}{n}} R)( 1+\varepsilon\frac{C}{n})\leq e^{4\varepsilon\frac{C}{n}} R,   \end{align*}
\begin{align*}\Longrightarrow\ \| B_n^k(\psi_0,\varepsilon)\|_{s+2}\leq e^{2k\varepsilon\frac{C}{n}} R\leq e^{2C} R\quad \quad \forall k\leq \tfrac{n}{\varepsilon}-1 .\end{align*}
%$$\| B_n^{k}(\psi_0,\varepsilon)\|_{s+2}\leq \big(e^{\varepsilon\frac{C}{n}}R( 1+\varepsilon\frac{C_1}{{n}^{3/2}}R^{2m}) \big)^{k}\leq e^C(R+1)^{(2m+1)k}.$$
Second, we note that $\| B_n^{k}(\psi_0,\varepsilon)\|_{L^2}=1$ for every $k\leq \tfrac{n}{\varepsilon}-1$. We apply Lemma \ref{Lemma:U_n} with radius $e^{2C}R$ and, possibly shrinking $\epsilon_1$, for every $0<\epsilon\leq \epsilon_1$, we have
\begin{equation*}\begin{split} \| B_n^{\tfrac{n}{\varepsilon}}(\psi_0,\varepsilon)- U(1, \nabla\phi)\psi_0 \|_{H^{s}}&\leq \| B_n(B_n^{\tfrac{n}{\varepsilon}-1}(\psi_0,\varepsilon),\varepsilon)- U(\tfrac{\varepsilon}{n}, \nabla\phi)B_n^{\tfrac{n}{\varepsilon}-1}(\psi_0,\varepsilon) \|_{H^{s}}\\
&+\| U(\tfrac{\varepsilon}{n}, \nabla\phi)B_n^{\tfrac{n}{\varepsilon}-1}(\psi_0,\varepsilon)- U(1, \nabla\phi)\psi_0 \|_{H^{s}}\\
&\lesssim \frac{C}{n}+\| U(\tfrac{\varepsilon}{n}, \nabla\phi)B_n^{\tfrac{n}{\varepsilon}-1}\psi_0- U(1, \nabla\phi)\psi_0 \|_{H^{s}}\\
&\leq \frac{C}{n}+e^{C\tfrac{\varepsilon}{n}}\| B_n^{\tfrac{n}{\varepsilon}-1}\psi_0- U(\tfrac{\varepsilon}{n}, \nabla\phi)^{\tfrac{n}{\varepsilon}-1}\psi_0 \|_{H^{s}}.
\end{split}\end{equation*}
By iterating these computations, we obtain the existence of $C>0$ such that
\begin{equation}\begin{split}\label{ineq_1} \| B_n^{\tfrac{n}{\varepsilon}}(\psi_0,\varepsilon)- U(1, \nabla\phi)\psi_0 \|_{H^{s}}&\leq \frac{C}{n}\sum_{l=0}^{\tfrac{n}{\varepsilon}-1}e^{C\tfrac{\varepsilon}{n}l}\lesssim  \frac{1}{\varepsilon}.
\end{split}\end{equation}
We repeat the same approach by considering the $L^2-$norm and by using Lemma \ref{Lemma:U_n_bis},
\begin{equation}\begin{split}\label{ineq_2} \| B_n^{\tfrac{n}{\varepsilon}}(\psi_0,\varepsilon)- U(1, \nabla\phi)\psi_0 \|_{L^2}&\leq \frac{C}{n^{\tfrac{3}{2}}}\sum_{l=0}^{\tfrac{n}{\varepsilon}-1}e^{C\tfrac{\varepsilon}{n}l}\lesssim  \frac{1}{\varepsilon\sqrt{n}}.
\end{split}\end{equation}
Now, by interpolation, there exist $a_1,a_2>0$ such that 
\begin{align*} \| B_n^{\tfrac{n}{\varepsilon}}(\psi_0,\varepsilon)- U(1, \nabla\phi)\psi_0 \|_{H^{s_1}}\leq\| B_n^{\tfrac{n}{\varepsilon}}(\psi_0,\varepsilon)- U(1, \nabla\phi)\psi_0 \|_{H^{s}}^{a_1}\| B_n^{\tfrac{n}{\varepsilon}}(\psi_0,\varepsilon)- U(1, \nabla\phi)\psi_0 \|_{L^2}^{a_2}\end{align*}
and then, thanks to \eqref{ineq_1} and \eqref{ineq_2}, there exists $C>0$ such that, for $n$ large,
\begin{align}\label{limit} \| B_n^{\tfrac{n}{\varepsilon}}(\psi_0,\varepsilon)- U(1, \nabla\phi)\psi_0 \|_{H^{s_1}}\lesssim \frac{1}{\varepsilon^{a_1+a_2}n^\frac{a_2}{2}} \xrightarrow{n\rightarrow \infty} 0.\end{align}
\end{proof}

\begin{lemma}\label{Lemma:B_n-U_n}

Let $s>\tfrac{N}{2}$ and $\phi\in C^{s+5}(\T^N;\R)$. For any $\eta>0$, for any $\psi_0\in H^{s}(\T^N;\C)$ such that $\|\psi_0\|_{L^2}=1$, there exists $\varepsilon_1$ such that, for any $0<\varepsilon\leq \varepsilon_1$ such that $\varepsilon^{-1}\in\N^*$, there exists $\tilde n$ such that for every $n\geq\tilde n$, we have
$$ \| B_n^{\tfrac{n}{\varepsilon}}(\psi_0,\varepsilon)- U(1, \nabla\phi)\psi_0 \|_{H^{s }}<\eta.$$
\end{lemma}
\begin{proof}
There exists $C_1>0$ such that, for every $f,g \in H^s(\T^N;\C)$, we have
$$\| U(1, \nabla\phi) f  - U(1, \nabla\phi)g \|_{H^{s }}\leq e^{C_1}\|   f  -  g \|_{H^{s }}.$$
We take $\eta_1\in (0,1]$ so that $\eta_1e^{C_1}<\tfrac{\eta}{4}$. Consider $\tilde\psi_0\in H^{s+3}(\T^N;\C)$ such that $\|\widetilde\psi_0\|_{L^2}=1 $ and $$\|\psi_0-\tilde\psi_0\|_{H^s}< \eta_1.$$
We use Lemma \ref{Lemma:B_n_lip}: there exists $C_2>0$ such that, for $R:=\|\psi_0\|_{H^{s}}+1$, there exists $\varepsilon_1>0$ sufficiently small such that, for every $0<\varepsilon\leq\varepsilon_1$, there exists $\tilde n\in\N^*$ such that for every $n\geq \tilde n$,
$$ \| B_n^{\tfrac{n}{\varepsilon}}(\psi_0,\varepsilon)- B_n^{\tfrac{n}{\varepsilon}}(\tilde\psi_0,\varepsilon) \|_{H^{s }}\leq  \big(e^{\varepsilon\frac{ C_2}{n} \big(\tfrac{ R^{2m}}{\sqrt{n}}+1\big)}\big)^{\tfrac{n}{\varepsilon}}\eta_1\leq e^{2C_2 }\eta_1%\leq e^{ C }\eta_1<\tfrac{\eta}{4}
.$$
By taking $\tilde\psi_0$ close enough to $\psi_0$ and then shrinking $\eta_1$, we have $e^{2C_2 }\eta_1<\tfrac{\eta}{4}$ and%for every $0<\varepsilon\leq\varepsilon_1$ and %$ \|  U(1, \nabla\phi)\psi_0 - U(1, \nabla\phi)\tilde \psi_0 \|_{H^{s }}\leq e^{C_1}\eta_1<\tfrac{\eta}{4}$. Then
$$ \| B_n^{\tfrac{n}{\varepsilon}}(\psi_0,\varepsilon)- U(1, \nabla\phi)\psi_0 \|_{H^{s }}<\tfrac{\eta}{4}+\| B_n^{\tfrac{n}{\varepsilon}}(\tilde\psi_0,\varepsilon)- U(1, \nabla\phi)\tilde\psi_0 \|_{H^{s }}+ \tfrac{\eta}{4}.$$
We use Lemma \ref{Lemma:B_n-U_n-pre} for the function $\tilde \psi_0\in H^{s+3}(\T^N;\C)$ by considering %$R:=\|\tilde \psi_0\|_{H^{s+3}}$ and 
$s_1\in [s,s+1)$ and, up to decreasing $\varepsilon_1$ enough, for any $0<\varepsilon\leq \varepsilon_1$ such that $\varepsilon^{-1}\in\N^*$, there exists $\tilde n$ such that for every $n\geq \tilde n$, we have
$$ \| B_n^{\tfrac{n}{\varepsilon}}(\psi_0,\varepsilon)- U(1, \nabla\phi)\psi_0 \|_{H^{s }}< \tfrac{\eta}{4}+\| B_n^{\tfrac{n}{\varepsilon}}(\tilde\psi_0,\varepsilon)- U(1, \nabla\phi)\tilde\psi_0 \|_{H^{s_1 }}+\tfrac{\eta}{4}<  \eta .$$
\end{proof}

\subsection{Proof of Proposition \ref{THM:control_app}}

We provide now the following definitions.

\begin{definition}
We denote $\Diff(\T^N)$ denotes the group of diffeomorphisms of $\T^N$ which are $C^\infty$, isotopic to the identity.
\end{definition}

\begin{definition}
For $P \in \Diff(\T^N)$, the unitary operator on $L^2(\T^N;\C)$
associated with $P$ reminder by
$$
\Lc_P \psi = |J_P|^{1/2} (\psi \circ P),
$$
where $|J_P| = \det(J_P)$ is the determinant of the Jacobian matrix of $P$.
Then $\|\Lc_P \psi\|_{L^2} = \|\psi\|_{L^2}$.
\end{definition}

\begin{definition}
 Let $s>\frac{N}{2}$. Consider an unitary operator $U$ on $L^2(\T^N;\C)$ such that $U:\psi \in H^s(\T^N;\C) \longmapsto \Lc_P(\psi e^{i\theta} )\in H^s(\T^N;\C)$ with $P \in \Diff(\T^N)$ and $\theta \in C^\infty(\T^N;\R)$. The operator $U$ is said $H^s$-STAR if, for every initial state $\psi_0 \in H^s(\T^N;\C)$ such that $\|\psi_0\|_{L^2}=1$ and $\eta > 0$, there exist $T \in [0,\eta]$ and a piecewise constant control $u$ defined in $(0,T)$ such that the corresponding solution $\psi(t)$ to \eqref{NLS-N} verifies
$$\|\psi(T) -  U\psi_0\|_{H^s} < \eta.$$
\end{definition}

\begin{definition}
$\vec(\T^N)$ denotes the space of %globally Lipschitz 
smooth vector fields on $\T^N$.
\end{definition} 

\begin{definition}
For $f \in \vec(\T^N)$, $\phi_f^t$ denotes the flow associated with $f$ at time $t$: for every $x_0 \in \T^N$, $x(t) = \phi_f^t(x_0)$ is the solution of
the ODE
$$
\dot{x}(t) = f(x(t))
$$
associated with the initial condition $x(0) = x_0$.
\end{definition} 

\begin{definition}
  Let us define  $\Gc := \big\{ \nabla \varphi \ :\  \varphi \in C^\infty (\T^N;\C)\big\}$.
\end{definition}

\begin{lemma}\label{Lemma:composition}
Let $s>\frac{N}{2}$. The composition and the strong limit of $H^s$-STAR operators are $H^s$-STAR operators.
\end{lemma}
\begin{proof}

\smallskip

\noindent
{\bf Stability by composition.} Let $L_1,L_2$ be $H^s$-STAR operators. Let $\psi_0\in H^s(\T^N;\C)$ such that $\|\psi_0\|_{L^2}=1$ and $\eta_1,\eta_2>0$. %Recall the definition provided in \eqref{def:Rn}.
For $j=1,2$, there exist $T_j>0$ and $u_j$ piecewise constant control in $(0,T_j)$, with $j=1,2$, such that
\begin{equation*}\begin{split}
\|R(T_2,L_1\psi_0,u_2)-L_2L_1\psi_0\|_{H^s}<\eta_2,\ \ \ \ \ \ \
\|R(T_1,\psi_0,u_1)-L_1\psi_0\|_{H^s}<\eta_1.
\end{split}\end{equation*}
Denote now $u$ the control in $(0,T_1+T_2)$ obtained by concatenating $u_1$ and $u_2$. We consider $\eta_2<\frac{\eta}{2}$ and, thanks to the continuity with respect to the initial state for the dynamics of \eqref{NLS-N} (see for instance \cite[Proposition 1.1]{Vahagn}), there exists $C>0$ (depending on $T_2$ and $u_2$) such that 
\begin{align*}
\|R(T_1+T_2,\psi_0,u)-L_2L_1\psi_0\|_{H^s}
&\leq
\|R(T_2,R(T_1,\psi_0,u_1),u_2)
- R(T_2,L_1\psi_0,u_2)\|_{H^s}\\
&+ \|R(T_2,L_1\psi_0,u_2)-L_2L_1\psi_0\|_{H^s} \\
&\leq
C \|R(T_1,\psi_0,u_1)-L_1\psi_0\|_{H^s} + \frac{\eta}{2}\leq
C \eta_1 + \frac{\eta}{2}
\end{align*}
To ensure the property, it is sufficient to program the first control phase so that $\eta_1 <\frac{\eta}{2 C}.$

\smallskip

\noindent
{\bf Stability by strong convergence.} Let $(L_n)_{n\in\N^*}$ be a sequence of $H^s$-STAR operators and let 
%$L$ be an unitary operator on $L^2(\T^N;\C)$ such that $L:H^s(\T^N;\C)\rightarrow H^s(\T^N;\C)$ and $L:X\rightarrow X$ which is the strong limit in $H^s(\T^N;\C)$ of $(L_n)_{n\in\N^*}$. 
$L$ be an unitary operator on $L^2(\T^N;\C)$ such that $L:\psi \in H^s(\T^N;\C) \longmapsto \Lc_P(\psi e^{i\theta} )\in H^s(\T^N;\C)$ with $P \in \Diff(\T^N)$ and $\theta \in C^\infty(\T^N;\R)$. Assume that $L$ is the strong limit in $H^s(\T^N;\C)$ of $(L_n)_{n\in\N^*}$, {\it i.e}
$$
\|(L_n-L)\psi\|_{H^s}\xrightarrow[]{n\rightarrow\infty} 0,\ \ \ \ \ \ \forall  \psi\in H^s(\T^N;\C).$$
Let $\psi_0\in H^s(\T^N;\C)$. By definition of the strong convergence, there exists $n\in\N^*$ such that
$$\|(L_n-L)\psi_0\|_{H^s}<\eta/2.$$ 
Now, thanks to the definition of $H^s-$STAR operator, there exist $T\in[0,\eta]$ and a piecewise constant control $u$ in $(0,T)$ such that
$$\|R(T,\psi_0,u)-L_n\psi_0\|_{H^s}<\eta/2.$$
The property is finally proved by using together the two previous inequalities so that 
$$\|R(T,\psi_0,u)-L\psi_0\|_{H^s}
\leq
\|R(T,\psi_0,u)-L_n\psi_0\|_{H^s}
+\|(L_n-L)\psi_0\|_{H^s}
<\eta.$$
\end{proof}

\begin{lemma}\label{Prop:contr-flow}
 Let $s>\frac{N}{2}$. Consider $f=\nabla\varphi\in\Gc$ and $P=\phi^1_f$. Then $\Lc_P$ is $H^s-$STAR. 
\end{lemma}
\begin{proof}
Note at first that the exponentials $e^{i \varepsilon\frac{|\nabla\phi|^2}{\sqrt{n}}}$, $e^{-i  \phi\sqrt{n}}$ and $e^{i  \phi\sqrt{n}}$ which can be found in the definition of $B_n$ provided in \eqref{def:Bn} are $H^s-$STAR thanks to Proposition \ref{THM:contr-phase} for every $\phi$ sufficiently smooth and $\varepsilon>0$. For any $\varepsilon^{-1}\in\N^*$ and $n\in\N^*$, the map $B_n^{n/\varepsilon}$ is reachable in time $\frac{n}{\varepsilon}\cdot \varepsilon n^{-3/2} = n^{-1/2}\rightarrow 0$ as $n\rightarrow\infty$. By Lemma \ref{Lemma:B_n-U_n}, $B_n^{n/\varepsilon}\rightarrow U(1,\nabla\phi)$ in $H^s$ as $n\to\infty$. Since $\Lc_P=U(1, f/2)$, the operator $\Lc_P$ is $H^s-$STAR by the stability under strong convergence of Lemma \ref{Lemma:composition}.
\end{proof}

\begin{proposition}\label{Prop:algebra}
Let $s>\frac{N}{2}$. Then $\Ac^s$ is a Lie subalgebra of $ \vec (\T^N)$ for
$$\Ac^s:= \{ f \in  \vec(\T^N) \ :\ \forall t \in \R,\ \mathcal{L}_{\phi_f^t} \text{ is } H^s\text{-STAR} \}.$$
\end{proposition}
\begin{proof}
The result is a consequence of the stability by summation proved in Lemma \ref{Lemma:sum-1} and the stability by commutation ensured by Lemma \ref{Lemma:commut-1}
\end{proof}

\begin{corollary}\label{Coro:algebra}
 Let $s>\frac{N}{2}$. Then, $\Ac^s$ contains $Lie(\Gc).$
\end{corollary}
\begin{proof}
Thanks to Lemma \ref{Prop:contr-flow}, $\Ac^s$ contains $ \Gc.$ The result follows from Proposition \ref{Prop:algebra}.

\end{proof}

\begin{lemma}\label{Prop:flow-density}
 Let $s>\frac{N}{2}$. For every $f\in \vec(\T^N)$, there exists $(f_k)_{k\in\N^*}\subset Lie(\Gc)$ such that $(\Lc_{\phi^1_{f_k}})_{k\in\N^*}$ converge strongly in the $H^s-$topology to $\Lc_{\phi^1_{f}}$.
\end{lemma}
\begin{proof}
The proof is essentially based on arguments developed to prove \cite[Theorem 18]{BP04-2}, which, in the case of the torus, is ensured by \cite[Proposition 27]{BP04-2}. Indeed, thanks to the first five steps of the proof of \cite[Proposition 27]{BP04-2} and to the Fejér's theorem, %applied to the function $f$ to 
we obtain the existence of a sequence of trigonometric polynomial vector fields such that, since $f \in C^\infty(\mathbb{T}^N)$,
$$\|f_k -f\|_{W^{s+1,\infty}}\xrightarrow[]{k\rightarrow \infty} 0.$$
%Now, $\|T_{f_k}-T_{f}\|_{L(H^2,H^1)}\xrightarrow[]{k\rightarrow \infty}$ since, for every $\psi_0\in H^2$ such that $\|\psi_0\|_{H^2}=1$, we have
%\begin{align*}\big\|\nabla \big(T_{f_k}-T_{f}\big)\psi_0\big\|_{L^2}&=\Big\|2\nabla \big( (f_k-f)\cdot\nabla\psi_0\big) +\nabla \big(div(f_k-f)\psi_0)\Big\|_{L^2} \\
%&\leq 2\|f_k -f\|_{W^{1,\infty}} \|\psi_0\|_{H^1}+2\|f_k -f\|_{L^{\infty}} \|\psi_0\|_{H^2}\\
%&+\|f_k -f\|_{W^{2,\infty}}\|\psi_0\|_{L^2}+\|f_k -f\|_{W^{1,\infty}}\|\psi_0\|_{H^1}. \end{align*}
%Also, there exists $C>0$ such that $\|T_{f_k}\|_{L(H^2,H^1)}<C$. 
%Now, there exists $C>0$ such that $\|e^{t T_{f}}\psi_0\|_{H^2}\leq C$ when $t\in [0,2]$. 
Assume at first $\psi_0\in H^{s+1}(\T^N;\C)$. Consider $\psi(t)=(e^{t T_{f_k}}-e^{t T_f})\psi_0$ solving in $H^s(\T^N;\C)$ the following equation
$$\partial_t \psi= T_{f_k}e^{t T_{f_k}}\psi_0 -T_{f}e^{t T_{f}}\psi_0= T_{f_k} \psi -(T_{f}-T_{f_k})e^{t T_{f}}\psi_0.$$
We estimate its $H^s-$norm when $t\in [0,2]$. Since $f\cdot\nabla$ is skew-adjoint up to lower-order terms, thanks to the Young inequality, there exist $C_1,C_2,C_3>0$ such that
\begin{align*}\partial_t \|\psi\|_{H^s}^2
&\leq C_1\Big( \|f_k\|_{W^{s+1,\infty}}\|\psi\|_{H^s} + \|f-f_k\|_{W^{s+1,\infty}}\|\psi_0\|_{H^{s+1}}\Big)\| \psi\|_{H^s}\\
&\leq  C_2\Big( \|\psi\|_{H^s}^2 + \|f-f_k\|_{W^{s+1,\infty}} \| \psi\|_{H^s}\Big)\leq  C_3\Big( \|\psi\|_{H^s}^2 + \|f-f_k\|_{W^{s+1,\infty}}^2\|\psi_0\|_{H^{s+1}}^2 \Big).\end{align*}
We apply the Gronwall inequality and we see that, in $t\in[0,2]$, there exists $C>0$ such that
$$\|\psi\|_{H^s}^2 \leq \|f-f_k\|_{W^{s+1,\infty}}^2\|\psi_0\|_{H^{s+1}}^2 e^C$$
$$\Longrightarrow \ \ \ \|\psi(t=1)\|_{H^s} \leq C\|f-f_k\|_{W^{s+1,\infty}}\|\psi_0\|_{H^{s+1}}\xrightarrow[]{k\rightarrow\infty}0.$$
Since $e^{tT_{f_k}}$ is unitary on $H^s$, the proof follows by density of $H^{s+1}$ in $H^s$.
\end{proof}

\begin{lemma}\label{Prop:control-flow-general}
 Let $s>\frac{N}{2}$. For every $f\in \vec(\T^N)$, $\Lc_{\phi^1_{f}}$ is $H^s-$STAR.
\end{lemma}
\begin{proof}
 The proof is a consequence of Corollary \ref{Coro:algebra} and Lemma \ref{Prop:flow-density}.
\end{proof}

\begin{lemma}\label{Prop:control-diffeo}
 Let $s>\frac{N}{2}$. For every $P\in \Diff(\T^N)$, the operator $\Lc_P$ is $H^s-$STAR.
\end{lemma}
\begin{proof}
As showed in the proof of \cite[Theorem 13]{BP04-2}, there exists $(f_k)_{k\leq K}\subset \vec(\T^N)$ such that $\Lc_P=\Lc_{\phi^1_{f_1}}\circ...\circ \Lc_{\phi^1_{f_K}}.$
The result is then ensured thanks to Lemma \ref{Lemma:composition} and Lemma \ref{Prop:control-flow-general}.
\end{proof}

\begin{proof}[Proof of Proposition \ref{THM:control_app}]
 Assume $\psi_1=\rho_1 e^{i\theta_1}, \psi_2=\rho_2 e^{i\theta_2}\in  X$. Thanks to \cite[Proposition~32]{BP04-2}, there exists a diffeomorphism $\Dc$ of the torus, of class $C^\infty$ and isotopic to the identity verifying the following property such that $\Lc_\Dc \rho_1 = \rho_2$ and then $  e^{i\theta_2}\Lc_\Dc e^{-i\theta_1}\psi^1=\psi_2$. The statement follows since $e^{i\theta_2}$, $\Lc_\Dc$ and $e^{-i\theta_1}$ are $H^s-$STAR thanks to Proposition \ref{THM:contr-phase} and Lemma \ref{Prop:control-diffeo}. The composition of the three is still $H^s-$STAR thanks to Lemma \ref{Lemma:composition}. 
\end{proof}

We now focus on the proof of Proposition \ref{THM:control_app_1} bet before we need the following proposition. 

\begin{proposition}\label{Prop:control-diffeo-2}
 Let $\psi,\varphi \in H^{s_1}(\T;\R_+^*)$ with $s_1>1/2$ such that $\|\psi\|_{L^2}=\|\varphi\|_{L^2}=1$ and $\frac{1}{2}< s<\min\{s_1,1\}$.  For every $\eta > 0$, there exists $  \D\in \Diff(\T)$ such that 
$$\|\Lc_{ \D}\psi-\varphi\|_{H^{s}} < \eta.$$
\end{proposition}
The proof of Proposition \ref{Prop:control-diffeo-2} is postponed to the Appendix \ref{proof_diffeo}.

\begin{proof}[Proof of Proposition \ref{THM:control_app_1}]
 Assume $\psi_1=\rho_1 e^{i\theta_1}, \psi_2=\rho_2 e^{i\theta_2}\in  \widetilde X$. By density and thanks to Proposition \ref{Prop:control-diffeo-2}, for $\eta_1,\eta_2,\eta_3>0$, there exist $\widehat \theta_1,\widehat \theta_2\in C^\infty(\T;\R)$ and $\Dc\in \Diff(\T)$ such that  
 $$\|  e^{i\theta_2} -  e^{i\widehat \theta_2}\|_{H^{s}}\leq \eta_1,\ \ \ \ \|\rho_2   -\Lc_\Dc\rho_1\|_{H^{s}}\leq \eta_2, \ \ \ \ \|e^{-i\theta_1}-e^{-i\widehat\theta_1}\|_{H^s}\leq \eta_3,$$
 \begin{align*}\|  e^{i\theta_2}\rho_2 -  e^{i\widehat \theta_2}\Lc_\Dc\rho_1 \|_{H^{s}}&\leq  \|  e^{i\theta_2}\rho_2 -  e^{i\widehat \theta_2}\rho_2 \|_{H^{s}} +\|  e^{i\widehat \theta_2}\rho_2 -  e^{i\widehat \theta_2}\Lc_\Dc\rho_1 \|_{H^{s}}\\
& \leq  \eta_1 \|   \rho_2 \|_{H^{s}} +\|  e^{i\widehat \theta_2} \|_{H^{s}}\eta_2.\end{align*}
As $\rho_2$ is fixed, we first set $\widehat \theta_2$ so that $\eta_1 \|   \rho_2 \|_{H^{s}}< \frac{\eta}{4}$ and after $\Dc$ so that $\|  e^{i\widehat \theta_2} \|_{H^{s}}\eta_2<\frac{\eta}{4}$. Thus, 
 \begin{align}\label{inne}\|  \psi_2 -  e^{i\widehat \theta_2}\Lc_\Dc\rho_1 \|_{H^{s}}=\|  e^{i\theta_2}\rho_2 -  e^{i\widehat \theta_2}\Lc_\Dc\rho_1 \|_{H^{s}}&<\frac{\eta}{2}.\end{align}
Note that $  e^{i\widehat \theta_2} \Lc_\Dc (e^{-i\widehat\theta_1} \psi_1 ) =e^{i\widehat \theta_2}\Lc_\Dc\rho_1 $. Since $e^{i\widehat\theta_2}$, $\Lc_\Dc$ and $e^{-i\widehat\theta_1}$ are $H^{s}-$STAR thanks to Proposition \ref{THM:contr-phase} and Lemma \ref{Prop:control-diffeo}, the composition of the three is still $H^{s}-$STAR thanks to Lemma \ref{Lemma:composition}. This property implies the existence of $u\in L^2((0,T);\R^{3})$ such that the corresponding mild solution $\psi$ to the equation \eqref{NLS-N} satisfies
$$\psi(0)=\psi_1,\qquad  \|\psi(T)-e^{i\widehat \theta_2}\Lc_\Dc\rho_1\|_{H^{s}}<\frac{\eta}{2}.$$
Finally, the last inequality, together with \eqref{inne}, infers
$$ \|\psi(T)-\psi_2\|_{H^{s}}\leq \|\psi(T)-e^{i\widehat \theta_2}\Lc_\Dc\rho_1\|_{H^{s}}+ \|  e^{i\widehat \theta_2}\Lc_\Dc\rho_1 -e^{i\theta_2}\rho_2  \|_{H^{s}}< \eta .$$
\end{proof}

\section{A preliminar local exact controllability result}\label{Sec:control_ground_proof}

The aim of this section is to ensure the small-time local exact controllability around % ensured by Theorem \ref{THM:control}. In this section, we consider different reference trajectories or hypotheses on the problem. 
%
%
%
%The first result on the subject is the following which is stated around 
the constant reference trajectory $\psi_\reff$ defined in \eqref{psireff}.

\begin{proposition}\label{THM:contr_ground}
Let $m\in\N^*$. Assume \begin{equation*}\kappa\not\in\Big\{\frac{2^{m-1}\pi^m}{m}(j^2+k^2):\ j,k\in\N^*,\ j\neq k\Big\}.\end{equation*} Let one of the following sets of hypotheses hold.

\begin{enumerate}[label=\textup{(\roman*)}, leftmargin=2em, rightmargin=2em]
\item Assume $\sqrt{2\widetilde \kappa  } \notin \N^*$. Consider \eqref{NLS} in the presence of $d=2$ controls $\mu = (\mu^{\mathsf D},\mu^{\mathsf N})$. Assume that the conditions \eqref{condmu}-\eqref{condmu2} are satisfied with $s>0$.

\item Assume $\sqrt{2\widetilde \kappa } \in \N^*$. Consider \eqref{NLS} in the presence of $d=3$ controls $\mu = (\mu^{\mathsf D},\mu^{\mathsf N},1)$. Assume that the conditions \eqref{condmu}-\eqref{condmu2} are satisfied with $s>0$.
\end{enumerate}

\noindent
%Then, \eqref{NLS} is small-time locally controllable around the reference trajectory $\psi_\reff$ in $H^s(\T;\C)$ in the following sense.
For any $T>0$, there exists $\eta>0$ such that, for any $\psi_1,\psi_2\in H^s(\T;\C)$ so that 
 $$\|\psi_1\|_{L^2}=\|\psi_2\|_{L^2}=1,\ \ \ \ \ \ \ \|\psi_1-\psi_\reff(0)\|_{H^s}\leq \eta,\ \ \text{ and }\ \ \|\psi_2-\psi_\reff(0)\|_{H^s}\leq \eta,$$ 
 there exists $u\in L^2((0,2T);\R^d)$ with $\|u\|_{L^2_{2T}}\les_T \eta$ such that there exists a unique mild solution $\psi$ to \eqref{NLS} (in the sense of Theorem \ref{THM:LWP_intro}), satisfying 
 $$\psi(0)=\psi_1,\qquad\psi(2T)=\psi_2.$$
 %There exists $\eta>0$ such that for any $T>0$ and $\psi_1\in H^s(\T;\C)$ so that 
 %$$\|\psi_1\|_{L^2}=1\ \ \text{ and }\ \ \|\psi_1-\psi_\reff(T)\|_{H^s}\leq \eta,$$ 
 %there exists $u \in L^2((0,T);\R^d)$ such that there exists a unique mild solution $\psi$ to \eqref{NLS} (in the sense of Theorem \ref{THM:LWP_intro}), satisfying 
% $$\psi(0)=\psi_1,\qquad\psi(T)=\psi_1.$$
\end{proposition}

\subsection{Linearized system and proof of Proposition \ref{THM:contr_ground}}

Let $\psi \in  L^2(\T;\R) \cap \Sc$, where $\Sc$ is the unit sphere in $L^2(\T;\R)$. Notice that
$$\Sc= \ker(g), \ \ \ \ \ \text{where}\ \ \ \ \ g(\varphi) = \|\varphi\|_{L^2}^2-1, \ \ g: L^2(\T;\C)\rightarrow\R.$$
The tangent space to $\Sc$ at $\psi$ is $T_{\psi}\Sc=\ker(d_\varphi g(\psi))$ where $d_\varphi g(\psi)\cdot \xi = 2\Re\jb{\psi,\xi}_{L^2}$ and then 
$$\Tc_{\psi}\Sc = \{\xi\in L^2(\T;\C),\ \Re\jb{\xi,\psi}_{L^2}=0\}.$$

We introduce and characterize now the linearization of \eqref{NLS}. To this end, we introduce the end-point map of the solution to~\eqref{NLS} and prove its $C^1$-regularity. The linearization we consider here is given by the equation solved by the differential of the end-point map. % whose proof follows the approach leading to \cite[Proposition 3]{BL}.

\begin{proposition} \label{smooth_p}
Let $u\in L^2((0,T);\R^d) $ with $d\in\N^*$ and $\psi$ the corresponding solution to~\eqref{NLS} in $(0,T)$ with any initial state $\psi_0$. Under the assumptions of Proposition~\ref{THM:LWP_intro}, the mapping
$$F:u\in L^2((0,T);\R^d)
\mapsto \psi(T)\in \Sc\cap  H^s(\T;\C)$$ is $C^1$, and for any $w\in L^2((0,T);\R^d)$, there holds $$\partial_{u} F( u): L^2((0,T);\R^d)\rightarrow  H^s(\T;\C),$$
$$(\partial_{u} F( u)\cdot w\big)(t)=\xi(t)\in \{\xi\in H^s(\T;\C)\ :\ \Re\jb{\xi,\psi}_{L^2}=0\} ,$$ where $\xi$ is the
linearized solution
\begin{equation*}\begin{split}
i\dt \xi(t) -\dx^2 \xi(t) =& (u(t)\cdot\mu(x))\xi(t) +(w(t)\cdot\mu(x))\psi(t)\\
&+m\kappa \psi^2(t)|\psi(t)|^{2(m-1)}\overline\xi+(m+1)\kappa |\psi(t)|^{2m}\xi,~~(t,x)\in\R\times\T.
\end{split}\end{equation*}
\end{proposition}
\begin{proof}
The result is rather standard and follows by the techniques adopted in the proof of \cite[Proposition 3]{BL} and by the well-posedness in $H^s(\T;\C)$ ensured by Proposition~\ref{THM:LWP_intro}. %Notice that the similar results are also ensured in ------- via the same techniques and in --- via similar arguments.
\end{proof}

We recall that $\psi_\reff(t)= \frac1{\sqrt{2\pi}} e^{-i \frac{\kappa}{(2\pi)^m} t}$. We consider the tangent spaces 
\begin{align*}
 \Tc_{\psi_\reff(t)}\Sc=&\left\{\xi\in L^2(\T;\C): \Re(\jb{\xi,\psi_\reff(t)}_{L^2})=0\right\}.
 \end{align*}
Note that $\psi_\reff(t)$ is solution of \eqref{NLS} with control $u=0$ and initial state $\psi_\reff(0)= \frac1{\sqrt{2\pi}}$. The corresponding linearized solution $\xi(t)$ belongs to $\Tc_{\psi_\reff(t)}\Sc$ for every $t\in (0,T)$ and solves
\begin{equation}\begin{split}\begin{cases}\label{NLS-lin-1}
i\dt \xi(t) -\dx^2 \xi(t) = (w(t)\cdot\mu(x)) \frac{e^{-i \frac{\kappa}{(2\pi)^m} t}}{\sqrt{2\pi}} +\frac{m\kappa}{(2\pi)^m} e^{-2 i \frac{\kappa}{(2\pi)^m} t}\overline\xi(t)+\frac{(m+1)\kappa}{(2\pi)^m} \xi(t),\quad \quad &(t,x)\in\R\times\T,\\
\xi(t=0)=h\in \Tc_{\frac1{\sqrt{2\pi}}}\Sc.\\
\end{cases}\end{split}\end{equation}
Let $v=\frac{w}{\sqrt{2\pi}}\in L^2((0,T);\R^d)$ and $\widetilde \kappa = \frac{m\kappa}{(2\pi)^m}$. We act a change of phase $\varphi(t)=e^{ i \frac{\kappa}{(2\pi)^m} t}\xi$ and the function $\varphi$ solves the linearized equation
\begin{equation}\label{NLS-lin}\begin{split}\begin{cases}
i\dt \varphi(t) -\dx^2 \varphi(t)-2 \widetilde \kappa \Re(\varphi(t)) = (v(t)\cdot\mu(x)),\quad \quad &(t,x)\in\R\times\T,\\
\varphi(t=0)=h\in \Tc_{\frac1{\sqrt{2\pi}}}\Sc.\\
\end{cases}\end{split}\end{equation}
Note that the solution $\varphi(t)$ belongs to the following tangent space for every $t\in(0,T)$:
\begin{align*}
 \Tc_{\frac1{\sqrt{2\pi}}}\Sc=&\left\{\varphi\in L^2(\T;\C): \frac1{\sqrt{2\pi}}\Re\left(\int_{\T} \varphi(x)dx\right)=0\right\}.
 \end{align*}

We aim to prove the local exact controllability of \eqref{NLS} around the spatially constant state $\phi(t)$. To this end, we first establish the global exact controllability of \eqref{NLS-lin} in $\Tc_{\frac1{\sqrt{2\pi}}}\Sc \cap H^s(\T;\C)$. The local controllability for the initial bilinear problem \eqref{NLS} then follows from the application of the Generalized Inverse Function Theorem.

\begin{proposition}\label{THM:control-periodic}
Let $m\in\N^*$. Assume \begin{equation}\label{resonances1}\kappa\not\in\Big\{\frac{2^{m-1}\pi^m}{m}(j^2+k^2):\ j,k\in\N^*,\ j\neq k\Big\}.\end{equation} Let one of the following sets of hypotheses hold.

\begin{enumerate}[label=\textup{(\roman*)}, leftmargin=2em, rightmargin=2em]
\item Assume $\sqrt{2\widetilde \kappa } \notin \N^*$. Consider \eqref{NLS} in the presence of $d=2$ controls $\mu = (\mu^{\mathsf D},\mu^{\mathsf N})$. Assume that the conditions \eqref{condmu}-\eqref{condmu2} are satisfied with $s>0$.

\item Assume $\sqrt{2\widetilde \kappa } \in \N^*$. Consider \eqref{NLS} in the presence of $d=3$ controls $\mu = (\mu^{\mathsf D},\mu^{\mathsf N},\widetilde \mu^{\mathsf N})$. Assume that the conditions \eqref{condmu}-\eqref{condmu2} are satisfied with $s>0$ and
\begin{equation*}
\ft{\widetilde{\mu}}^{\mathsf N}_k = 0, \quad \forall k \in \Z^* ; \qquad
\ft{\widetilde{\mu}}^{\mathsf N}_{0} \neq 0.
\end{equation*}
\end{enumerate}
Then, for any $T>0$ and $\varphi_1 \in\Tc_{\frac1{\sqrt{2\pi}}}\Sc \cap H^s(\T;\C)$, there exists a control $u\in L^2((0,T);\R^d)$ such that the corresponding solution $\varphi$ to the problem~\eqref{NLS-lin} with $\varphi(0)=0$ satisfies  $$\varphi(T)= \varphi_1.$$

\end{proposition}
\begin{proof} \smallskip

\noindent
{\bf Preliminaries.} Let us assume at first that we only have two controls and 
$$\mu=(\mu^{\mathsf N},\mu^{\mathsf D}),\ \ \ \ \ \ \text{and}\ \ \ \ \  \ v = (v^{\mathsf N},v^{\mathsf D}) .$$
We begin by %understanding how to express the free dynamics $e^{-i t A}\varphi_0$ in a convenient form, where $A$ is the linear operator $\Delta + 2   \widetilde \kappa \Re$ with domain $H^2(\T;\C)$. To this end, we denote
denoting
$$f=\Re(\varphi), \quad\quad g=i\Im(\varphi).$$ We observe that \eqref{NLS-lin} correspond to the following coupled system
\begin{equation}\label{NLS-lin-sist}\begin{split}\begin{cases}
i\dt f(t) -\dx^2 g(t) = 0,\quad \quad &(t,x)\in\R\times\T,\\
i\dt g(t) - \dx^2 f(t)-2  \widetilde \kappa f(t) = (v(t)\cdot\mu(x)),\quad \quad &(t,x)\in\R\times\T,\\
f(t=0)=g(t=0)=0.\\
\end{cases}\end{split}\end{equation}
We now rewrite the dynamics in terms of Fourier modes, using the decompositions 
$f(t)=\sum_{k\in\Z} f_k(t)\, e_k$ and $g(t)=\sum_{k\in\Z} g_k(t)\, e_k$, 
where $f_k(t) = \jb{f(t),e_k}$ and $g_k(t) = \jb{g(t),e_k}$. 
With this notation, the problem \eqref{NLS-lin-sist} takes the form
\begin{equation}\label{NLS-lin-sist-four}\begin{split}\begin{cases}
i\dt f_k(t) +k^2 g_k(t) = 0,\quad \quad & t \in\R ,\quad\forall k\in \Z,\\
i\dt g_k(t) +(k^2-2  \widetilde \kappa) f_k(t) = (v(t)\cdot\mu_k),\quad \quad & t \in\R ,\quad\forall k\in \Z,\\
f_k(t=0)=g_k(t=0)=0.\\
\end{cases}\end{split}\end{equation}
The last equation can be written in a matrix form as follows
\begin{equation}\label{NLS-lin-sist-matrix}\begin{split}
i\dt \begin{pmatrix}
f_k   \\
g_k  
\end{pmatrix} =
M_k\cdot
\begin{pmatrix}
f_k   \\
g_k  
\end{pmatrix}+\begin{pmatrix}
0   \\
(v(t)\cdot\mu_k)  
\end{pmatrix},\quad \quad & t \in\R ,\quad\forall k\in \Z,
\end{split}\end{equation}
where, for every $k\in \Z$, we have 
$$M_k=\begin{pmatrix}
0 & -k^2 \\
-k^2+ 2  \widetilde \kappa & 0
\end{pmatrix}.$$
We now notice that the spectrum of each matrix depends on the sign of $\widetilde \kappa$ and then of $k^2-2 \widetilde \kappa$. %To this purpose, we distinguish the following two cases.
%
%
%
%
%\noindent
%{\bf Real and complex spectrum.} %Assume that $\widetilde \kappa>0$. 
The characterization of the spectrum of each matrix $M_k$ allows us to study the dynamics of the time depending parameters $f_k(t)$ and $g_k(t)$. We have three distinct situations as characterized below.

\smallskip

\noindent
{\bf Nonzero real spectrum.}   When $k\in\Z^*$ is such that $k^2(k^2-2 \widetilde \kappa)>0$, the spectrum of each $M_K$ is composed by two distinct real eigenvalues
$$\sigma(M_k)=\{\sqrt{k^2(k^2-2 \widetilde \kappa)},-\sqrt{k^2(k^2-2 \widetilde \kappa)}\}$$
and the corresponding eigenvectors are of the form
$$ \begin{pmatrix}
-k^2   \\
\sqrt{k^2(k^2-2 \widetilde \kappa)}  
\end{pmatrix},\ \ \ \ \ \text{and}\ \ \ \  \ \begin{pmatrix}
k^2   \\
\sqrt{k^2(k^2-2 \widetilde \kappa)}  
\end{pmatrix}.$$
We denote $I_1\subset Z$ the set of indices verifying $k^2(k^2-2 \widetilde \kappa)>0$. When $k\in I_1$, 
$$
e^{-i t M_k}\begin{pmatrix}
v_1  \\
v_2  
\end{pmatrix}= 
\begin{pmatrix}
k^2\Big(b e^{-i t \sqrt{k^2(k^2-2 \widetilde \kappa)}} - a e^{i t \sqrt{k^2(k^2-2\widetilde \kappa)}}\Big)   \\
\sqrt{k^2(k^2-2 \widetilde \kappa)} \Big(a e^{i t \sqrt{k^2(k^2-2 \widetilde \kappa)}} + b e^{-i t \sqrt{k^2(k^2-2 \widetilde \kappa)}}\Big) 
\end{pmatrix}$$
where $a,b\in \C$ are defined by initial condition and
$$a=\frac{1}{2\sqrt{k^2(k^2-2 \widetilde \kappa)}}v_2-\frac{1}{2 k^2}v_1,\ \ \ \ \ b=\frac{1}{2 k^2}v_1+\frac{1}{2\sqrt{ k^2(k^2-2 \widetilde \kappa)}}v_2.$$
Finally, we can develop and factor out
$$
e^{-i t M_k}\begin{pmatrix}
v_1  \\
v_2  
\end{pmatrix}= \frac{1}{2}
\begin{pmatrix}
 ( v_1+\frac{k^2}{\sqrt{k^2(k^2-2 \widetilde \kappa)}}v_2) e^{-i t \sqrt{k^2(k^2-2 \widetilde \kappa)}} - (\frac{k^2}{\sqrt{k^2(k^2-2 \widetilde \kappa)}}v_2- v_1) e^{i t \sqrt{k^2(k^2-2 \widetilde \kappa)}}    \\
 (v_2-\frac{\sqrt{k^2(k^2-2 \widetilde \kappa)}}{k^2}v_1) e^{i t \sqrt{k^2(k^2-2 \widetilde \kappa)}} + (\frac{\sqrt{k^2(k^2-2 \widetilde \kappa)}}{k^2}v_1+v_2) e^{-i t \sqrt{k^2(k^2-2 \widetilde \kappa)}}
\end{pmatrix}.$$
We write solutions to \eqref{NLS-lin-sist} by keeping in mind that the initial data is $0$ and, for $t\in [0,T]$,
\begin{equation}\label{sol-real}
\begin{pmatrix}
f_k(t)  \\
g_k(t)  
\end{pmatrix}= \frac{1}{2}
\begin{pmatrix}
 -i\frac{k^2}{\sqrt{k^2(k^2-2  \widetilde \kappa)}}\int_0^t(v(\tau)\cdot\mu_k)  \Big(  e^{i (t-\tau) \sqrt{k^2(k^2-2 \widetilde \kappa)}} - e^{-i (t-\tau) \sqrt{k^2(k^2-2 \widetilde \kappa)}}  \Big) d\tau \\
-i  \int_0^t(v(\tau)\cdot\mu_k)\Big(  e^{-i (t-\tau) \sqrt{k^2(k^2-2 \widetilde \kappa)}} +  e^{i (t-\tau) \sqrt{k^2(k^2-  2\widetilde \kappa)}}\Big)d\tau
\end{pmatrix}.\end{equation}

\smallskip

\noindent
{\bf Purely complex spectrum.}   When $k^2(k^2-2 \widetilde \kappa)<0$, we identify two distinct complex eigenvalues 
$$\sigma(M_k)=\{i\sqrt{k^2(-k^2+2 \widetilde \kappa)}, -i\sqrt{k^2(-k^2+2 \widetilde \kappa)}\}.$$
and the corresponding eigenvectors are of the form
$$ \begin{pmatrix}
-k^2   \\
i\sqrt{k^2(-k^2+2 \widetilde \kappa)}  
\end{pmatrix},\ \ \ \ \ \text{and}\ \ \ \  \ \begin{pmatrix}
k^2   \\
i\sqrt{k^2(-k^2+2\widetilde \kappa)}  
\end{pmatrix}.$$
Notice that this case is verified for at most a finite number of $k\in\Z$ and we call $I_2\subset Z$ the corresponding indices. As above, when $k\in I_2$, 
$$
e^{-i t M_k}\begin{pmatrix}
v_1  \\
v_2  
\end{pmatrix}= \frac{1}{2}
\begin{pmatrix}
 ( v_1+\frac{k^2}{i\sqrt{k^2(-k^2+2 \widetilde \kappa)}}v_2) e^{- t \sqrt{k^2(-k^2+2 \widetilde \kappa)}} - (\frac{k^2}{i\sqrt{k^2(-k^2+2 \widetilde \kappa)}}v_2- v_1) e^{ t \sqrt{k^2(-k^2+2 \widetilde \kappa)}}    \\
 (v_2-\frac{i\sqrt{k^2(-k^2+2 \widetilde \kappa)}}{k^2}v_1) e^{ t \sqrt{k^2(-k^2+2 \widetilde \kappa)}} + (\frac{i\sqrt{k^2(-k^2+2 \widetilde \kappa)}}{k^2}v_1+v_2) e^{- t \sqrt{k^2(-k^2+2 \widetilde \kappa)}}
\end{pmatrix}.$$
and the solutions to \eqref{NLS-lin-sist} for $t\in [0,T]$ are
\begin{equation}\label{sol-complex}
\begin{pmatrix}
f_k(t)  \\
g_k(t)  
\end{pmatrix}= \frac{1}{2}
\begin{pmatrix}
 -\frac{k^2}{\sqrt{k^2(-k^2+2 \widetilde \kappa)}}\int_0^t(v(\tau)\cdot\mu_k)  \Big(  e^{- (t-\tau) \sqrt{k^2(-k^2+2 \widetilde \kappa)}} - e^{ (t-\tau) \sqrt{k^2(-k^2+2 \widetilde \kappa)}}  \Big)d\tau  \\
  -i\int_0^t(v(\tau)\cdot\mu_k)\Big(  e^{ (t-\tau) \sqrt{k^2(-k^2+2 \widetilde \kappa)}} +  e^{- (t-\tau) \sqrt{k^2(-k^2+2 \widetilde \kappa)}}\Big)d\tau
\end{pmatrix}.\end{equation}

\smallskip

\noindent
{\bf The zero eigenvalue.} When $k^2(k^2-2 \widetilde \kappa)=0$, we have one eigenvalue with algebraic multiplicity $2$ and geometric multiplicity $1$
$$\sigma(M_k)=\{0\}$$
and for identify the eigenvectors we need to distinguish between different cases. When $\sqrt{2   \widetilde \kappa}\in \N^*$, we have the $0$ eigenvalue for $k=\pm \sqrt{2  \widetilde \kappa}$ and the corresponding eigenvector is of the form
$$ \begin{pmatrix}
 c  \\
0  
\end{pmatrix},\ \ \ \ \ c\in \C.$$
When $k=0$, we find again the $0$ eigenvalue and the corresponding eigenvector is of the form
$$\begin{pmatrix}
0   \\
c  
\end{pmatrix},\ \ \ \ \ c\in \C.$$
When $\sqrt{2   \widetilde \kappa}\in \N^*$ and $k=\pm \sqrt{2   \widetilde \kappa}$, the free dynamics generate by $M_k$ is given by
$$
e^{-i t M_k}\begin{pmatrix}
v_1  \\
v_2  
\end{pmatrix}= 
\begin{pmatrix}
v_1+ i k^2 v_2 t   \\
v_2
\end{pmatrix}.$$
Now, the solutions to \eqref{NLS-lin-sist} for $t\in [0,T]$ are
\begin{equation}\label{sol-zero-nonzero}
\begin{pmatrix}
f_k(t)  \\
g_k(t)  
\end{pmatrix}= 
\begin{pmatrix}
   k^2\int_0^t(v(\tau)\cdot\mu_k)   (t-\tau ) d\tau  \\
  -i\int_0^t(v(\tau)\cdot\mu_k)d\tau
\end{pmatrix}.\end{equation}
Finally, when $k=0$, we have
$$
e^{-i t M_0}\begin{pmatrix}
v_1  \\
v_2  
\end{pmatrix}= 
\begin{pmatrix}
v_1  \\
v_2- 2i  \widetilde \kappa  v_1 t 
\end{pmatrix}.$$
The solutions to \eqref{NLS-lin-sist} for $t\in [0,T]$ in this case are
\begin{equation}\label{sol-zero-zero}\begin{pmatrix}
f_0(t)  \\
g_0(t)  
\end{pmatrix}= 
\begin{pmatrix}
 0  \\
  -i\int_0^t(v(\tau)\cdot\mu_0)    d\tau  
\end{pmatrix}.\end{equation}
In what follows, we denote the (finite) indices involved in this case $I_3= \Z\setminus (I_1\cup I_2)$.

\smallskip

\noindent
{\bf Some intermediate computations.} We now consider the different formulations of the solutions \eqref{sol-real},\eqref{sol-complex},\eqref{sol-zero-nonzero} and \eqref{sol-zero-zero} in $t=T$ by keeping in mind that if $k\in I_1$ (or $k\in I_2$ respectively), then $-k\in I_1$ (or $-k\in I_2$ respectively).
Let us introduce the sequence of numbers $\{\sigma_k\}_{k\in\Z}$ such that
$$\sigma_{k}=\begin{cases}\sqrt{k^2(k^2-2 \widetilde \kappa)}\ \ \ &\text{if}\ k\in I_1,\\
i\sqrt{k^2(-k^2+2 \widetilde \kappa)}\ \ \ &\text{if}\ k\in I_2,\\
0,\ \ \ &\text{if}\ k\in I_3.
\end{cases}$$
Now, we consider for every $k\in\N$ the following numbers
\begin{equation*}\begin{split}f_k^{\mathsf N}:=f_k(T) + f_{-k}(T),\ \ \ \ \ \ g_k^{\mathsf N}:=g_k(T) + g_{-k}(T),\\
f_k^{\mathsf D}:=f_k(T) - f_{-k}(T),\ \ \ \ \ \ g_k^{\mathsf D}:=g_k(T) - g_{-k}(T),\end{split}\end{equation*}
and, we recall $v(\tau)\cdot(\mu_k\pm\mu_{-k})=v^{\mathsf D}(\tau) (\mu_k^{\mathsf D}\pm\mu_{-k}^{\mathsf D}) + v^{\mathsf N}(\tau) (\mu_k^{\mathsf N}\pm\mu_{-k}^{\mathsf N}) $, which implies:
\begin{equation}\label{decoupling}v(\tau)\cdot(\mu_k+\mu_{-k})=2 v^{\mathsf N}(\tau) \mu_k^{\mathsf N} ,\ \ \ \ \ v(\tau)\cdot(\mu_k-\mu_{-k})=2v^{\mathsf D}(\tau)  \mu_k^{\mathsf D}  . \end{equation}
Now, we know that for every $k\in \N\setminus I_3$, 
\begin{equation*}\label{sol-real-1}\begin{split}
\begin{pmatrix}
f_k^{\mathsf N}  \\
g^{\mathsf N}_k  
\end{pmatrix}&=  
\begin{pmatrix}
 -i\frac{k^2}{\sigma_k}\int_0^T\Big(v(\tau)\cdot(\mu_k+\mu_{-k})\Big)  \Big(  e^{i (T-\tau) \sigma_k} - e^{-i (T-\tau) \sigma_k}  \Big) d\tau \\
-i  \int_0^T(v(\tau)\cdot(\mu_k+\mu_{-k})\Big)\Big(  e^{-i (T-\tau) \sigma_k} +  e^{i (T-\tau) \sigma_k}\Big)d\tau
\end{pmatrix},\\
\begin{pmatrix}
f_k^{\mathsf D}  \\
g^{\mathsf D}_k  
\end{pmatrix}&=  
\begin{pmatrix}
 -i\frac{k^2}{\sigma_k}\int_0^T\Big(v(\tau)\cdot(\mu_k-\mu_{-k})\Big)  \Big(  e^{i (T-\tau) \sigma_k} - e^{-i (T-\tau) \sigma_k}  \Big) d\tau \\
-i  \int_0^T(v(\tau)\cdot(\mu_k-\mu_{-k})\Big)\Big(  e^{-i (T-\tau) \sigma_k} +  e^{i (T-\tau) \sigma_k}\Big)d\tau
\end{pmatrix}.\end{split}\end{equation*}
The identity \eqref{decoupling} allow us to notice that the last two systems can be rewritten as two systems respectively involving the controls $v^{\mathsf N}$ and $v^{\mathsf D}$ so that, for $k\in \N\setminus I_3$:
\begin{equation*}\label{sol-real-2}
\begin{pmatrix}
f_k^j  \\
g^j_k
\end{pmatrix}=  
\begin{pmatrix}
 -2i\frac{k^2}{\sigma_k}\mu_k^j \int_0^T v^j(\tau)   \Big(  e^{i (T-\tau) \sigma_k} - e^{-i (T-\tau) \sigma_k}  \Big) d\tau \\
-2i \mu_k^j \int_0^T v^j(\tau)  \Big(  e^{-i (T-\tau) \sigma_k} +  e^{i (T-\tau) \sigma_k}\Big)d\tau
\end{pmatrix},\ \ \ \ \ \ \ j=N,D.\end{equation*}
We now consider $k\in \N\cap I_3$. When $\sqrt{2   \widetilde \kappa}\in \N^*$ and $k= \sqrt{2  \widetilde \kappa}$,  we have
\begin{equation*}\label{sol-zero-nonzero-1}
\begin{pmatrix}
f_k^j \\
g_k^j 
\end{pmatrix}= 
\begin{pmatrix}
  2  k^2 \mu_k^j \int_0^T v^j(\tau)   (T-\tau ) d\tau  \\
  -2 i \mu_k^j \int_0^T v^j(\tau)d\tau
\end{pmatrix},\ \ \ \ \ \ \ j=N,D.\end{equation*}
Finally, when $k=0$,  we have
\begin{equation*}\label{sol-zero-zero-1}\begin{pmatrix}
f_0^j \\
g_0^j  
\end{pmatrix}= 
\begin{pmatrix}
 0  \\
  -i2\mu_0^j\int_0^Tv^j(\tau)    d\tau  
\end{pmatrix},\ \ \ \ \ \ \ j=N,D.\end{equation*}
Let us introduced the numbers $h_k^{\mathsf N} $ and $h_k^{\mathsf D} $  for $k\in \Z*$ as follows. When $k\in \Z\setminus I_3$, the first are defined by
\begin{equation}\label{mp-geneal-nonzero}\begin{split}
h_k^j&= \frac{\sigma_k}{ 4 k^2} f_k^j + \frac{1}{ 4 } g_k^j= -i \mu_k^j \int_0^T v^j(\tau)    e^{i (T-\tau) \sigma_k}d\tau,\ \ \ \ \ \ \ j=N,D, \ \ \ \ \text{if}\ \ k>0,\\
h_k^j&= \frac{1}{ 4 } g_{-k}^j-\frac{\sigma_k}{ 4 k^2} f_{-k}^j= -i \mu_{-k}^j \int_0^T v^j(\tau)    e^{-i (T-\tau) \sigma_k}d\tau,\ \ \ \ \ \ \ j=N,D \ \ \ \ \text{if}\ \ k<0.
\end{split}\end{equation}
We now consider $k\in \N\cap I_3$. When $\sqrt{2 \widetilde \kappa}\in \N^*$, we denote 
\begin{equation}\label{mp-geneal-zero-nonzero}\begin{split}
h_k^j&= \frac{1}{ 2 k^2 } f_k^j-  \frac{i T}{ 2 } g_{k}^j = -\mu_k^j \int_0^T  \tau\, v^j(\tau)     d\tau,\ \ \ \ \ \ \ j=N,D, \ \ \ \ \text{if}\ \ k=\sqrt{2  \widetilde \kappa},\\
h_k^j&= \frac{1}{ 2 } g_{-k}^j= -i\mu_{-k}^j\int_0^Tv^j(\tau)    d\tau,\ \ \ \ \ \ \ j=N,D \ \ \ \ \text{if}\ \ k=-\sqrt{2  \widetilde \kappa}.
\end{split}\end{equation}
Finally, we introduce the numbers $\widehat h^{\mathsf N}$ and $\widetilde h^{\mathsf N}$ 
\begin{equation}\label{mp-geneal-zero-zero}\begin{split}
\widehat h^{\mathsf N}_0&=  f_0^{\mathsf N} = 0,\\
\widetilde h^{\mathsf N}_0&= \frac{1}{ 2 } g_0^{\mathsf N}= -i\mu_0^{\mathsf N}\int_0^Tv^{\mathsf N}(\tau)    d\tau,\ \ \ \ \ \ \ j=N,D,
\end{split}\end{equation}
and the sequence of vectors $\{h_k\}_{k\in\Z}$ such that $h_k=(h^{\mathsf N}_k,h_k^{\mathsf D})$ for every $k\in \Z^*$, and $h_0=\begin{pmatrix}
  \widehat h^{\mathsf N}_0 \\
\widetilde h^{\mathsf N}_0
\end{pmatrix}$.

\smallskip

\noindent
{\bf Some properties of the controlled space.} Our purpose is to prove the controllability in $\Tc_{\frac{1}{\sqrt{2\pi}}}\Sc\cap  H^s(\T;\C)$ of the lihearized system \eqref{NLS-lin} by studying the same property for the equivalent system \eqref{NLS-lin-sist}. We want to prove that, for every $T>0$ and $\varphi_1$ in such space, there exists a control $v\in L^2((0,T),\R^q)$ such that the corresponding solutions $f$ and $g$ to the problem~\eqref{NLS-lin-sist} with $f(0)=g(0)=0$ satisfy  
$$f(T)= \Re(\varphi_1),\ \ \  \ \ \ g(T)= i\Im(\varphi_1).$$
We write $\Re(\varphi_1)=\sum_{k\in\Z} \Re(\varphi_1)_k e_k$ and $\Im(\varphi_1)=\sum_{k\in\Z} \Im(\varphi_1)_k e_k$ and we observe that, for every $k\in\Z^*$,  
$$\Re(\varphi_1)_k = \overline{\Re(\varphi_1)_{-k}},\ \ \ \ \ \Im(\varphi_1)_k=\overline{\Im(\varphi_1)_{-k}}.$$
The controllability follows the surjectivity, for every $T>0$, of the mapping 
$$v\in L^2((0,T),\R^q)\mapsto   \left\{\begin{pmatrix}
  f_k(T)  \\
  g_k(T)
\end{pmatrix}\right\}_{k\in \Z^*}\in L,\ \ \ \text{with}$$
$$L=\left\{\left\{\begin{pmatrix}
   x_k  \\
  y_k
\end{pmatrix}\right\}_{k\in \Z}:\{ x_k\}_{k\in \Z},\{ y_k\}_{k\in \Z}\in h^s(\Z;\C),\  x_{-k}=\overline{ x_{k}},\  y_{-k}=-\overline{ y_{k}},\  x_0=0,\  y_0\in i\R\right\},$$
Now, we notice that the mapping $$A:\left\{\begin{pmatrix}
  f_k(T)  \\
  g_k(T)
\end{pmatrix}\right\}_{k\in \Z}\mapsto \left\{h_k\right\}_{k\in \Z},$$
considered as an application from $\ell^2(\Z;\C^2)$ to $\ell^2(\Z;\C^2)$ is invertible since $f_0^{\mathsf D}=g_0^{\mathsf D}=0$ by definition
and $\mu_0^{\mathsf D}=0$ by hypothesis. Now, $A$ maps the space $L$ in 
\begin{equation*} \begin{split} 
A(L)=\Bigg\{&\left\{\begin{pmatrix}
  x_k  \\
  y_k
\end{pmatrix}\right\}_{k\in \Z^*}:\{x_k\}_{k\in \Z},\{y_k\}_{k\in \Z}\in h^s(\Z;\C),\  x_{-k}=-\overline{x_{k}},\ \ y_{-k}=\overline{y_{k}},\ \ k\in I_1;\\
& x_{k}\in i\R,\ \ y_{k}\in \R,\ \ k\in I_2;\ \   x_{\sqrt{2  \widetilde \kappa} }\in \R,\ \ y_{\sqrt{2  \widetilde \kappa}}\in i \R,\\
& x_{-\sqrt{2 \widetilde \kappa} }\in i \R,\ \ y_{-\sqrt{2 \widetilde \kappa}}\in \R,\ \ \text{if}\ \ \sqrt{2 \widetilde \kappa}\in \N;\ \ x_0=0,\ y_0\in i\R\Bigg\}.
\end{split}\end{equation*}
Indeed, consider any $\left\{\begin{pmatrix}
  x_k  \\
  y_k
\end{pmatrix}\right\}_{k\in \Z^*}\in L$. We can observe that, for every $k\in I_1\cap \N$, we have
\begin{equation*} \begin{split} 
\tilde x_{-k}&:=\frac{1}{4}(y_{k}+y_{-k})-\frac{\sigma_{k}}{4k^2}(x_{k}+x_{-k})=-\frac{1}{4}(\overline y_{-k}+\overline y_{k})-\frac{\sigma_k}{4k^2}(\overline x_{-k}+\overline x_{k})=-\overline{\tilde x}_{k}, \\
\tilde y_{-k}&:=\frac{1}{4}(y_{k}-y_{-k})-\frac{\sigma_{k}}{4k^2}(x_{k}-x_{-k})=-\frac{1}{4}(\overline y_{-k}-\overline y_{k})-\frac{\sigma_k}{4k^2}(\overline x_{-k}-\overline x_{k})=\overline{\tilde y}_{k}.\\
\end{split}\end{equation*}
When $k\in I_2$, since $\sigma_k\in i\R$, we have
\begin{equation*} \begin{split} 
\tilde x_{k}&:=\frac{\sigma_k}{4k^2}(x_{k}+x_{-k}) + \frac{1}{4}(y_{k}+y_{-k})   =\overline{-\frac{\sigma_k}{4k^2}(x_{-k}+x_k) - \frac{1}{4}(y_{-k}+y_{k})}=-\overline{\tilde x}_{k}, \\
\tilde y_{k}&:=\frac{\sigma_k}{4k^2}(x_{k}-x_{-k}) + \frac{1}{4}(y_{k}-y_{-k})   =\overline{-\frac{\sigma_k}{4k^2}(x_{-k}-x_k) - \frac{1}{4}(y_{-k}-y_{k})}=\overline{\tilde y}_{k}, \\
\tilde x_{-k}&:=\frac{1}{4}(y_{k}+y_{-k}) - \frac{\sigma_{k}}{4k^2}(x_{k}+x_{-k})   =\overline{-\frac{1}{4}(y_{-k}+y_{k}) + \frac{\sigma_{k}}{4k^2}(x_{-k}+x_{k})}=-\overline{\tilde x}_{-k}, \\
\tilde y_{-k}&:=\frac{1}{4}(y_{k}-y_{-k}) - \frac{\sigma_{k}}{4k^2}(x_{k}-x_{-k})   =\overline{-\frac{1}{4}(y_{-k}-y_{k}) + \frac{\sigma_{k}}{4k^2}(x_{-k}-x_{k})}=\overline{\tilde y}_{-k}. \\\end{split}\end{equation*}
Finally, when $\sqrt{2 \widetilde \kappa}\in \N$, we have for $k=\sqrt{2 \widetilde \kappa}\in \N$:
\begin{equation*} \begin{split} 
\tilde x_{k}&:=\frac{1}{2k^2}(x_{k}+x_{-k}) - \frac{iT}{2}(y_{k}+y_{-k})  =\overline{\frac{1}{2k^2}(x_{-k}+x_k)- \frac{i T}{2}(y_{-k}+y_{k}) } =\overline{\tilde x}_{k}, \\
\tilde y_{k}&:=\frac{1}{2k^2}(x_{k}-x_{-k})- \frac{iT}{2}(y_{k}-y_{-k})  =\overline{\frac{1}{2k^2}(x_{-k}-x_k)- \frac{i T}{2}(y_{-k}-y_{k})}=-\overline{\tilde y}_{k}, \\
\tilde x_{-k}&:=\frac{1}{2}(y_{k}+y_{-k})    =\overline{-\frac{1}{2}(y_{-k}+y_{k}) }=-\overline{\tilde x}_{-k}, \\
\tilde y_{-k}&:=\frac{1}{2}(y_{k}-y_{-k})    =\overline{-\frac{1}{2}(y_{-k}-y_{k}) }=\overline{\tilde y}_{-k}. \\
\end{split}\end{equation*}
Now, we consider the mapping 
$$ A_1:\left\{\begin{pmatrix}
  f_k(T)  \\
  g_k(T)
\end{pmatrix}\right\}_{k\in \Z}\mapsto \left\{\ ...\ \begin{pmatrix}
i\frac{h_{-k}^{\mathsf N}}{\mu_{-k}^{\mathsf N}}e^{i T \sigma_k}\\
i\frac{h_{-k}^{\mathsf D}}{\mu_{-k}^{\mathsf D}}e^{i T \sigma_k}
\end{pmatrix},\ ...\ ,
\begin{pmatrix}
i\frac{\widehat h_0}{\mu_0^{\mathsf N}}\\
i\frac{\widetilde h_0}{\mu_0^{\mathsf N}}
\end{pmatrix},\ ...\ \begin{pmatrix}
i\frac{h_{k}^{\mathsf N}}{\mu_k^{\mathsf N}}e^{-i T \sigma_k}\\
i\frac{h_{k}^{\mathsf D}}{\mu_k^{\mathsf D}}e^{-i T \sigma_k}
\end{pmatrix}\ ...\ 
\right\}_{k\in \Z},$$
is also invertible from $h^s(\Z;\C^2)$ to $\ell^2(\Z;\C^2)$ thanks to the hypotheses on $\mu^{\mathsf D}$ and $\mu^{\mathsf N}$. In addition, $$\mu_k^{\mathsf N}= \mu_{-k}^{\mathsf N}\in\R\ \ \ \ \text{ and }\ \ \ \mu_k^{\mathsf D}=-\mu_{-k}^{\mathsf D}\in i\R,$$ which implies that $A_1$ maps the space $L$ in 
\begin{equation*} \begin{split} 
A_1(L)=\Bigg\{\left\{\begin{pmatrix}
  x_k  \\
  y_k
\end{pmatrix}\right\}_{k\in \Z^*}:\ &\{x_k\}_{k\in \Z},\{y_k\}_{k\in \Z}\in \ell^2(\Z;\C),\  x_{-k}=\overline{x_{k}},\ \ y_{-k}=\overline{y_{k}},\ \ k\in I_1;\\
& x_{k},y_{k}\in \R,\ \ k\in I_2;\ \   x_{\sqrt{2 \widetilde \kappa} },y_{\sqrt{2 \widetilde \kappa}}\in i\R,\\
& x_{-\sqrt{2 \widetilde \kappa} }, y_{-\sqrt{2 \widetilde \kappa}}\in \R\ \ \text{if}\ \ \sqrt{2  \widetilde \kappa}\in \N;\ \ x_0=0,\ y_0\in \R\Bigg\}.
\end{split}\end{equation*}

\smallskip

\noindent
{\bf The two moment problems.} Thanks to the observations above, we can study controllability by ensuring the surjectivity of the mapping $A_1\circ F :L^2((0,T),\R^d) \rightarrow \ell^2(\Z;\C^2)$, which corresponds to the solvability of the following moment problems in $A_1(L)$: one involving the control $v^{\mathsf N}$ and the other $v^{\mathsf D}$. When $k\in \Z\setminus I_3$, thanks to \eqref{mp-geneal-nonzero}, we have
\begin{equation}\label{mp-nonzero}\begin{split}
x_k&=  \int_0^T v^{\mathsf N}(\tau)    e^{-i \tau \sigma_k}d\tau, \ \ \ \ \ \ y_k=  \int_0^T v^{\mathsf D}(\tau)    e^{-i \tau \sigma_k}d\tau,\ \  \ \ \ \ \text{if}\ \ k>0,\\
x_k&=   \int_0^T v^{\mathsf N}(\tau)    e^{i \tau \sigma_k}d\tau, \ \ \ \ \ \ y_k=   \int_0^T v^{\mathsf D}(\tau)    e^{i \tau \sigma_k}d\tau,\ \ \ \ \ \ \text{if}\ \ k<0.
\end{split}\end{equation}
When $\sqrt{2 \widetilde \kappa}\in \N^*$, from \eqref{mp-geneal-zero-nonzero},
\begin{equation}\label{mp-zero-nonzero}\begin{split}
x_k&= -i\int_0^T  \tau v^{\mathsf N}(\tau)  d\tau,\ \ \ \ \ \ \ \ \ y_k = i\int_0^T\tau\, v^{\mathsf D}(\tau) d\tau,\ \ \ \ \ \ \ \ \ \text{if}\ \ k=\sqrt{2 \widetilde \kappa},\\
x_k &=    \int_0^Tv^{\mathsf N}(\tau)    d\tau,\ \ \ \ \ \ \ \ \  y_k =    \int_0^Tv^{\mathsf D}(\tau)    d\tau,\ \ \ \ \ \ \ \ \  \text{if}\ \ k=-\sqrt{2 \widetilde \kappa},
\end{split}\end{equation}
and, since $x_0=0$, we have only one last equation which is 
\begin{equation}\label{mp-zero-zero}\begin{split}
y_0&=    \int_0^Tv^{\mathsf N}(\tau)    d\tau.
\end{split}\end{equation}
In other words, we want to prove that, for every $T>0$ and for every $\left\{\begin{pmatrix}
  x_k  \\
  y_k
\end{pmatrix}\right\}_{k\in \Z^*}\in A_1(L)$, there exists $v\in L^2((0,T);\R^q)$ such that all the identities \eqref{mp-nonzero}, \eqref{mp-zero-nonzero} and \eqref{mp-zero-zero} are verified. 
Here, we notice that when $\sqrt{2 \widetilde \kappa}\in \N^*$, we can not expect to ensure the above solvability and then the controllability by using only two controls since it is not true in general that 
$$y_0=   -i  \int_0^Tv^{\mathsf N}(\tau)    d\tau =x_{-\sqrt{2 \widetilde \kappa}}.$$
Nevertheless, when $\sqrt{2 \widetilde \kappa} \not\in \N^*$, we can infer the solvability of the moment problems by using Corollary \ref{PROP:solvability-gen}, which is ensured thanks to the validity of the non-resonant hypothesis \eqref{resonances1} and concludes the proof of the first point of the theorem. The second point, instead, is ensured by repeating the previous theory in the case where we have three controls. In this case, everything is almost the same, since when $k \neq 0$, we have $\mu_k^{\mathsf N} \neq 0$ but $\tilde{\mu}^N_k=0$ with $\tilde{\mu}^N=1$ is our third control. The controllability is then ensured by proving the solvability of the moment problems \eqref{mp-nonzero} and \eqref{mp-zero-nonzero} with $k \in \Z^*$ for the controls $v^{\mathsf N}$ and $v^{\mathsf D}$ again by using Corollary \ref{PROP:solvability-gen}. We define $\widetilde v^{\mathsf N}$ instead by solving the equation  
$$
{\Im(\varphi_1)_0} = -i \mu_0^{\mathsf N} \int_0^T v^{\mathsf N}(\tau) \, d\tau - i  \int_0^T \widetilde v^{\mathsf N}(\tau) \, d\tau
$$  
with $\varphi_1$ the target of the dynamics and the control $v^{\mathsf N}$ previously found for the moment problem with frequencies $k\in\Z^*$.
\end{proof}

\begin{proof}[Proof of Proposition \ref{THM:contr_ground}]
The result follows as a standard application of the Inverse Mapping Theorem to the map $u \mapsto F(u)$ defined in Proposition~\ref{smooth_p}. This strategy has been widely used in the literature, starting from its exploitation in \cite[Proof of Theorem~1]{BL}. We briefly outline how this approach applies in our setting. First, $\psi_{ref}(t)$ denotes the solution of \eqref{NLS} corresponding to the control $u=0$ and the initial condition $\psi_{ref}(0)= \frac1{\sqrt{2\pi}}$. Second, Proposition~\ref{THM:LWP_intro} and~\ref{smooth_p} ensure that the map $F : u \mapsto \psi(T)$ is well defined and of class $C^1$ in $H^s$. Third, the derivative $\partial_u\Psi_T(u=0)\cdot v=\xi(T)$ is given by the endpoint of the mild solution to the linearized equation~\eqref{NLS-lin-1}, which can be rewritten, after a change of variables, in the form \eqref{NLS-lin}. 
Now, the exact controllability of \eqref{NLS-lin}, and hence of the linearized equation~\eqref{NLS-lin-1}, is ensured by Proposition \ref{THM:control-periodic}. The Inverse Mapping Theorem ensures that every $\psi_1$ and $\psi_2$ close to $\psi_\reff(0)$ can both be steered exactly to $\psi_\reff(T)$ with two controlled dynamics in time $T$. The time-reversibility of the Schrödinger equation then allows one to steer $\psi_1$ to $\psi_\reff(T)$ and then, running the second dynamics backward, to $\psi_2$ in total time $2T$, concluding the proof.

\end{proof}

\section{Proofs of Theorem \ref{THM:control_global_intro} and Theorem \ref{THM:control_global_intro_eigen}}\label{Sec:global_proof}

The aim of this section is to establish small-time global exact controllability for the system \eqref{NLS} and then prove Theorem \ref{THM:control_global_intro}. %The main result is formulated in the following theorem, from which Theorem \ref{THM:control_global_intro} is subsequently deduced. The key difference here is that we now make explicit the assumptions on the controls $\mu^{\mathsf D}$ and $\mu^{\mathsf N}$ that guarantee the controllability result.
%\begin{theorem}\label{THM:control_global}
%Let $\kappa\in\R$. Consider the \eqref{NLS} in the presence of $d=5$ controls $$\mu = (\mu^{\mathsf D},\mu^{\mathsf N}, 1, \sin(x), \cos(x)).$$ Assume the conditions \eqref{condmu} and \eqref{condmu2} be verified with $s>0$. Then, the \eqref{NLS} is small-time globally exactly controllable in $H^s(\T;\C)$ for every $s>0$ in the following sense.
% For any $T>0$ and any $\psi_1,\psi_2\in H^s(\T;\C)$ with $\|\psi_1\|_{L^2}=\|\psi_2\|_{L^2}$, there exists $u\in L^2((0,T);\R^5)$ such that there exists $\psi$ a unique mild solution to \eqref{NLS} (in the sense of Theorem \ref{THM:LWP_intro}), satisfying
%$$\psi(0)=\psi_1,\qquad  \psi(T)=\psi_2.$$
%\end{theorem}
%Theorem \ref{THM:control_global_intro}
The result is obtained by combining a suitable small-time local exact controllability result in low-regularity spaces with a small-time global approximate controllability result in the same functional setting (Proposition \ref{THM:control_app_1}). %We now focus on the former, which, together with Proposition \ref{THM:control_app_1}, yields Theorem \ref{THM:control_global_intro}. 
% is provided by Theorem \ref{THM:control} and is proved below. The second is established in the subsequent part of this section.
%The two results and the time-reversibility of the Schrödinger equation ensure Theorem \ref{THM:control_global}.

\begin{proof}[Proof of Theorem~\ref{THM:control_global_intro}] 

 \smallskip

\noindent
We only prove the statement in the case $r=1$, since the general case $r\neq 0$ follows from rescaling. Let $s_1>1/2$, and $\psi_1 = \rho_1 e^{i\theta_1}$ and $\psi_2 = \rho_2 e^{i\theta_2}$, with $\rho_1,\rho_2 \in H^{s_1}(\T;\R_+^*)$ and $\theta_1,\theta_2 \in H^{s_1}(\T;\R)$, satisfying $\|\psi_1\|_{L^2}=\|\psi_2\|_{L^2}=1$. Moreover, let $T>0$, $s\in(0,\min(s_1,1))$. Let $\eta > 0$ be sufficiently small so that Proposition \ref{THM:contr_ground} allows to steer any two states close to $e_0$ in $H^s(\T;\C)$ one into the other in time $\widetilde\tau = T/2$. This is done with 2 or 3 controls, depending on the value of $\sqrt{2\widetilde{\kappa}}$; the third potential needed if $\sqrt{2\widetilde{\kappa}}\in\N^*$ being the constant function.

By Proposition \ref{THM:control_app_1} with $s_1$ and $s'\in(\max(s,1/2),\min(s_1,1))$, there exist $T_1, T_2 \in (0, \min\{\eta/2,\widetilde\tau\}]$ and controls $u_1\in L^2((0,T_1);\R^5)$, $u_2\in L^2((0,T_2);\R^5)$ such that the corresponding mild solutions $\Psi_1$ and $\Psi_2$ of \eqref{NLS} satisfy
$$
\Psi_1(0) = \psi_1,  \qquad  \|\Psi_1(T_1) - e_0\|_{H^{s}}\leq \|\Psi_1(T_1) - e_0\|_{H^{s'}} < \min\{\eta/2,\widetilde\tau\},
$$
$$
\Psi_2(0) = \psi_2, \qquad \|\Psi_2(T_2) - e_0\|_{H^{s}}\leq \|\Psi_2(T_2) - e_0\|_{H^{s'}} < \min\{\eta/2,\widetilde\tau\}.
$$

Applying now Proposition \ref{THM:contr_ground} in $H^s(\T;\C)$, there exist controls $\tilde u_1, \tilde u_2 \in L^2((0,\widetilde\tau);\R^5)$ (with 2 or 3 components equal to 0, depending on the cases in Proposition \ref{THM:contr_ground}) such that the corresponding mild solutions $\tilde\Psi_1$ and $\tilde\Psi_2$ to \eqref{NLS} satisfy
$$
\tilde\Psi_1(0) = \Psi_1(T_1), \qquad \tilde\Psi_1(\widetilde\tau) = e_0,
$$
$$
\tilde\Psi_2(0) = \Psi_2(T_2), \qquad \tilde\Psi_2(\widetilde\tau) = e_0.
$$
We then concatenate the trajectories $\Psi_1$ and $\tilde \Psi_1$ to obtain an evolution steering $\psi_1$ to $e_0$ in time $T_1 + \widetilde\tau$. Since $e_0$ is a stationary solution associated with a suitable constant control $u=(0,0,-\widetilde{\kappa},0,0)$.

We can thus extend the trajectory so as to connect $\psi_1$ to $e_0$ in exactly in a time $T/2$. 
The same construction applies to $\psi_2$ which can be moved to $e_0$ in a time $T/2$. Finally, using the time reversibility of the flow, we construct a trajectory steering $e_0$ to $\psi_2$ in time $T/2$ and the proof is concluded by concatenating these evolutions.

\end{proof}

\begin{proof}[Proof of Theorem \ref{THM:control_global_intro_eigen}] 
The proof is a straightforward adaptation to the one of Theorem \ref{THM:control_global_intro}. We only need to replace the use of Proposition~\ref{THM:control_app_1} to send $\psi_1$ and $\psi_2$ near $e_0$, by Corollary \ref{COR:control_modes} (with $s'=s$) to send $e_n$ and $e_m$ near $e_0$.

\end{proof}

 \begin{remark}
  \rm
  It is likely that our techniques for the low regularity local controllability result around the constant trajectory extend to fractional NLS with dispersion $(-\dx^2)^\alpha$, $\alpha>1$, provided that we work at sufficient regularity $s(\alpha)$. Indeed, the moment method for the controllability of the linearised equation around $\psi_\reff$ applies verbatim in this regime, and the well-posedness should persist if $s$ is large enough.
 \end{remark}

%\begin{proof}[Proof of Theorem \ref{THM:contr_ground-more}]
%The result is very similar to the one of Theorem \ref{THM:contr_ground-more} with the only difference that one has to use the %additional controls $\tilde \mu^N,\mu$
%\end{proof}

\appendix

\section{Auxiliary algebraic properties}

The aim of this appendix is to ensure some technical lemmas which are used in the Section \ref{Sec:global-approx}.
\begin{lemma}
Let $f,g\in \vec(\T^N)$.  
Denote by $[f,g]$ their Lie bracket, which in the flat setting of $\mathbb T^N$
is defined by the formula $[f,g]=f\cdot\nabla g-g\cdot\nabla f$ and componentwise by
$$[f,g]^i =\sum_{j=1}^N f^j\partial_j g^i - g^j\partial_j f^i,\ \ \ \ \ i\leq N.$$
Then
$\ div([f,g])
= f\cdot\nabla(div(g))
- g\cdot\nabla(div(f))$.
\end{lemma}
\begin{proof}
By definition, $div([f,g])= \sum_{i=1}^N\sum_{j=1}^N\partial_i\big(f^j\partial_j g^i - g^j\partial_j f^i\big)$ and
$$
div([f,g])
=\sum_{i=1}^N \sum_{j=1}^N(\partial_i f^j)(\partial_j g^i)
+ f^j\partial_i\partial_j g^i
- (\partial_i g^j)(\partial_j f^i)
- g^j\partial_i\partial_j f^i.
$$
We commute the partial derivatives and we observe that 
$$\sum_{i=1}^N\sum_{j=1}^Nf^j\partial_i\partial_j g^i
= f^j\partial_j(\partial_i g^i)
= f\cdot\nabla(div(g))$$
and $\sum_{i=1}^N\sum_{j=1}^N g^j\partial_i\partial_j f^i
= g\cdot\nabla(div(f))$.
Moreover, the mixed terms cancel, since
$$\sum_{i=1}^N\sum_{j=1}^N(\partial_i f^j)(\partial_j g^i)
=\sum_{i=1}^N\sum_{j=1}^N (\partial_j f^i)(\partial_i g^j).$$
Collecting the terms, we conclude the proof.
\end{proof}

\begin{lemma}
Let $f,g\in \vec(\T^N)$. Then $[f\cdot\nabla,\; g\cdot\nabla] = [f,g]\cdot\nabla.$
\end{lemma}

\begin{lemma}\label{Lemma:comm-vec}
Let $f,g\in \vec(\T^N)$.
Then
$[T_f,T_g] = T_{2[f,g]}= 2T_{[f,g]}$.
\end{lemma}

\begin{proof}
Let us recall that $T_f=2f\cdot\nabla+div(f)$. Notice that we have the identities $[f\cdot\nabla,g\cdot\nabla]=[f,g]\cdot\nabla$
and
$ div([f,g])
= f\cdot\nabla(div(g))-g\cdot\nabla(div(f))$ thanks to the previous lemma.
We compute
$$[T_f,T_g]
=4[f\cdot\nabla,g\cdot\nabla]
+2\big(f\cdot\nabla(div(g))-g\cdot\nabla(div(f))\big)
=4[f,g]\cdot\nabla+2div([f,g]).$$
Hence
$[T_f,T_g]
=2\big(2[f,g]\big)\cdot\nabla+div(2[f,g])\bigr)
= T_{\,2[f,g]}$.
\end{proof}

We now need the following lemma, which follows from Lemma \ref{Lemma:comm-vec} and can be interpreted as a Baker-Campbell-Hausdorff formula (see, for instance, \cite[Theorem 7.4]{Serre} and its proof).

\begin{lemma}\label{Lemma:baker}
Let $f,g\in \vec(\T^N)$. There exist two operators $O_1$ and $O_2$ in $L(H^{s+1}; H^s)$ for every $s\in \R$ with $s>\frac{N}{2}$ such that
%$$e^{tT_f}e^{tT_g}=e^{t T_{f+g} + t^2 O_1},$$
$$e^{tT_f}e^{tT_g}=e^{t T_{f+g} +  t^2T_{[f,g]} + t^3 O_1},$$
$$e^{-tT_f}e^{-tT_g}e^{tT_f}e^{tT_g}=e^{ 2 t^2 T_{ [f,g]} + t^3 O_2}.$$
In addition, for every $s\in \R$ with $s>\frac{N}{2}$, there exists $C>0$, only depending in $f$ and $g$, such that, for every $\psi\in H^{s+1}(\T^N;\C)$, there holds
\begin{align}\label{baker-form}\Re\Big(\langle\psi,O_1\psi\rangle_{H^s}\Big),\Re\Big(\langle\psi,O_2\psi\rangle_{H^s}\Big)\leq C\|\psi\|_{H^s}^2.\end{align}
\end{lemma}
\begin{proof}
Consider $\psi_0\in H^1(\T^N;\C)$ and $\varphi(t)=e^{tT_f}e^{tT_g}\psi_0$. The function $\varphi$ verifies the equation
$$\partial_t\varphi= T_f\varphi + e^{tT_f} T_g e^{-tT_f}\varphi.$$
Now, we denote $ ad_{T_f}^kT_g$ and $ ad_f^kg$ with $k\in\N$ by iteration as follows $$ ad_{T_f}^{0}T_g= T_g ,\ \ \ \ \  ad_f^0g=g,\ \ \ \ \ ad_{T_f}^{k+1}T_g=[T_f,ad_{T_f}^{k}T_g],\ \ \ \ \  ad_f^kg=[f,ad_{f}^{k}g].$$
Thanks to Lemma \ref{Lemma:comm-vec} and since $T_f+T_g=T_{f+g}$, we have
\begin{align*}\partial_t\varphi&=  T_f \varphi + \sum_{k=0}^\infty \frac{1}{k!} ad^k_{t T_f}T_g \varphi\\
&= (T_f+T_g)\varphi + t[T_f,T_g] + t^2 \Big[T_f,\Big[T_f,\sum_{k=0}^\infty \frac{1}{(k+2)!} ad^k_{t T_f}T_g\Big]\Big] \varphi\\
&= (T_f+T_g)\varphi + 2tT_{[f,g]} + t^2 T_{[f,[f,\sum_{k=0}^\infty \frac{2^{k+2}}{(k+2)!} ad^k_{t f}g]]} \varphi.
\end{align*}
Then, the first identity is proved by estimation of the energy and by considering $O_1=3 T_{\big[f,\big[f,\sum_{k=0}^\infty \frac{2^{k+2}}{(k+2)!} ad^k_{t f}g\big]\big]}$.
The second one follows from the same approach by changing the defintion of $\varphi(t)$ in $e^{-tT_f}e^{-tT_g}e^{tT_f}e^{tT_g}\psi_0$.
\end{proof}

\begin{lemma}\label{Lemma:sum}
 Let $s\in \R$ with $s>\frac{N}{2}$ and $f,g\in \vec(\T^N)$. There exists $C>0$ such that 
$$ \Big\| U\big(\frac{1}{n},f\big)U\big(\frac{1}{n},g\big)-U\big(\frac{1}{n},f+g\big) \Big\|_{L(H^{s+1},H^{s})}<\frac{C}{n}.$$
\end{lemma}
\begin{proof}
    We start by applying Lemma \ref{Lemma:baker} and there exists $O\in L(H^{s+1}(\T^N;\C); H^{ s}(\T^N;\C))$ for every $s\geq 0$ such that $e^{\frac{1}{n}T_f}e^{\frac{1}{n}T_g}=e^{\frac{1}{n}(T_f+T_g) + \frac{1}{n^2} O}.$  Consider $\psi_0\in H^{s+1}(\T^N;\C)$ such that $\|\psi_0\|_{H^{s}}\leq \|\psi_0\|_{H^{s+1}}=1$ and denote
    $$\psi(t)= e^{\frac{t}{n}(T_f+T_g) + \frac{t}{n^2} O}\psi_0,$$
$$\varphi(t)=\psi(t)-w(t),\ \ \ \ \ \ \ w(t):=U(n^{-1}t,f+g)\psi_0=e^{\frac{t}{n}T_{f+g}}\psi_0.$$
The function $\varphi$ solves the following well-posed in $H^{s+1}(\T^N;\C)$ partial differential equation
\begin{equation}\begin{split}\label{EDP_11}
\partial_t \varphi
=& \frac{1}{n}(T_f+T_g) \varphi(t) + \frac{1}{n^2} O (\varphi(t)+w(t)).
\end{split}\end{equation}
We proceed by doing an estimation of the energy in $t\in[0,2]$. We recall that $\|\psi(t)\|_{L^2}\leq 1$ and the existence of $C>0$ such that $\|w(t)\|_{H^s}\leq C$ for $t\in [0,2]$. We use \eqref{baker-form} and, thanks to the Young inequality, there exist $C_1,C_2>0$ such that
%\begin{align*}
%\partial_t\| \varphi\|_{H^s}^2&\leq C_1\Big(\frac{1}{n}\| \varphi\|_{H^s}+\frac{1}{n}\|\varphi\|_{L^2} +\frac{1}{n^2}\|\varphi+w\|_{H^s} \Big)\|\varphi\|_{H^s}
%&\leq C_2\Big(\frac{1}{n}\|\varphi\|_{H^s}+\frac{1}{n }\ \|\varphi\|_{H^s}^2\Big)\\
%\leq C_2\Big(\frac{1}{ n^2}+\Big(1+\frac{1}{n }\Big)\|\varphi\|_{H^s}^2\Big).
%\end{align*}
\begin{align*}
\partial_t\| \varphi\|_{H^s}^2&\leq C_1\Big(\frac{1}{n}\| \varphi\|_{H^s}+\frac{1}{n}\|\varphi\|_{L^2} +\frac{1}{n^2}\|\varphi\|_{H^s}+\frac{1}{n^2}\|w\|_{H^{s+1}} \Big)\|\varphi\|_{H^s}
\leq C_2\Big(\frac{1}{ n^2}+\Big(1+\frac{1}{n }\Big)\|\varphi\|_{H^s}^2\Big).
\end{align*}
Thus, thanks to the Gronwall inequality, there exists $C>0$ such that
$$\big\| U\big(\frac{1}{n},f\big)U\big(\frac{1}{n},g\big)\psi_0-U\big(\frac{1}{n},f+g\big)\psi_0 \Big\|_{H^s}=\|\varphi(t=1)\|_{H^s}\leq \frac{C}{n} .$$
\end{proof}

\begin{lemma}\label{Lemma:sum-L2}
Let $f,g\in \vec(\T^N)$. There exists $C>0$ such that 
$$ \Big\| U\big(\frac{1}{n},f\big)U\big(\frac{1}{n},g\big)-U\big(\frac{1}{n},f+g\big) \Big\|_{L(H^1,L^2)}<\frac{C}{n^2}.$$
\end{lemma}
\begin{proof}
The proof is the same of the one of Lemma \ref{Lemma:sum}, with the only difference that one has to work in $L^2$ with $\psi_0\in H^1(\T^N;\C)$ such that $\|\psi_0\|_{H^1}=1$. Also in this case, we define $\varphi$ solving the equation \eqref{EDP_11}. Again, we proceed to make an estimation of the energy in $t\in[0,2]$ by observing that $T_f$ and $T_g$ are skew-adjoint. For $n\in\N^*$ sufficiently large, there exist $C_1,C_2>0$ such that
\begin{align*}
\partial_t\| \varphi\|_{L^2}^2&\leq C_1\Big(\frac{1}{n^2}\|\varphi\|_{L^2}+\frac{1}{n^2}\|w\|_{H^1} \Big)\|\varphi\|_{L^2}\leq C_2\Big(\frac{1}{n^4}+\Big(1+\frac{1}{n^2 }\Big)\|\varphi\|_{L^2}^2\Big).
\end{align*}
The Gronwall inequality concludes the proof again since there exists $C>0$ such that
$$\Big\| U\big(\frac{1}{n},f\big)U\big(\frac{1}{n},g\big)\psi_0-U\big(\frac{1}{n},f+g\big)\psi_0 \Big\|_{L^2}=\|\varphi(t=1)\|_{L^2}\leq \frac{1}{n^2}e^C.$$
\end{proof}

\begin{lemma}\label{Lemma:sum-1}
 Let $s\in \R$ with $s>\frac{N}{2}$ and $f,g\in \vec(\T^N)$. Consider $\psi_0\in H^{s}(\T^N;\C)$. For every $\eta>0$, there exists $\tilde n\in\N^*$, such that, for every $n>\tilde n$, we have 
$$ \Big\| \Big(U\big(\frac{1}{n},f\big)U\big(\frac{1}{n},g\big)\Big)^n\psi_0- U(1,f+g)\psi_0 \Big\|_{H^{s}}<\eta.$$
\end{lemma}
\begin{proof}
Assume $\psi_0\in H^{s+2}(\T^N;\C)$. %The generalization to $\psi_0\in H^{s}(\T^N;\C)$ holds by density arguments as in the proof of Lemma \ref{Lemma:B_n-U_n}. 
The proof follows the same strategy leading to Lemma \ref{Lemma:B_n-U_n}. Indeed, we notice at first
\begin{equation*}\begin{split}
  &\Big\| \big(U(\frac{1}{n},f)U(\frac{1}{n},g)\big)^n\psi_0- \big(U(\frac{1}{n},f+g)\big)^n\psi_0 \Big\|_{H^{s+1}}\\
  &\leq \sum_{k=0}^{n-1}\Big\| \big(U(\frac{1}{n},f)U(\frac{1}{n},g)\big)^k\Big(  U(\frac{1}{n},f)U(\frac{1}{n},g) -  U(\frac{1}{n},f+g)\Big)   \big(U(\frac{1}{n},f+g)\big)^{n-1-k}\psi_0 \Big\|_{H^{s+1}}.
\end{split}\end{equation*}
We recall the existence of $C>0$ such that $\|U(\frac{1}{n}, f)\|_{L(H^{\tilde s})},\|U(\frac{1}{n}, g)\|_{L(H^{\tilde s})}\leq e^{\frac{C}{n}}$ for every ${\tilde s}\in\N$. Now, thanks to Lemma \ref{Lemma:sum}, there exist $C_1,C_2>0$ such that
\begin{equation*}\begin{split}
  &\Big\| \big(U(\frac{1}{n},f)U(\frac{1}{n},g)\big)^n\psi_0- \big(U(\frac{1}{n},f+g)\big)^n\psi_0 \Big\|_{H^{s+1}}\\
%  &\leq \sum_{k=0}^{n-1}e^{\frac{Ck}{n}}\Big\|\Big(  U(\frac{1}{n},f)U(\frac{1}{n},g) -  U(\frac{1}{n},f+g)\Big)   \big(U(\frac{1}{n},f+g)\big)^{n-1-k}\psi_0 \Big\|_{H^{s}}\\
  &\leq \sum_{k=0}^{n-1}e^{\frac{Ck}{n}}\frac{C_1}{n}\big\|\big(U(\frac{1}{n},f+g)\big)^{n-1-k}\psi_0 \big\|_{H^{s+2}}\leq \sum_{k=0}^{n-1}e^{\frac{Ck}{n}}\frac{C_1}{n} e^{\frac{C(n-1-k)}{n}} \| \psi_0 \|_{H^{s+2}} \leq C_2\|\psi_0\|_{H^{s+2}}.
    \end{split}\end{equation*}
We repeat the approach with the $L^2$-norm and thanks to Lemma \ref{Lemma:sum-L2}, there exists $C>0$ so that
\begin{equation*}\begin{split}
  &\Big\| \big(U(\frac{1}{n},f)U(\frac{1}{n},g)\big)^n\psi_0- \big(U(\frac{1}{n},f+g)\big)^n\psi_0 \Big\|_{L^2}  \leq \frac{C}{n}\|\psi_0\|_{H^1}.
    \end{split}\end{equation*}
An interpolation inequality and a density argument, as in the proof of Lemma \ref{Lemma:B_n-U_n}, prove the result in $H^s(\T^N;\C)$.
\end{proof}

\begin{lemma}\label{Lemma:commut}
 Let $s\in \R$ with $s>\frac{N}{2}$ and $f,g\in \vec(\T^N)$. There exists $C>0$ such that 
$$ \Big\| U\big(-\frac{1}{\sqrt{2n}},f\big)U\big(-\frac{1}{\sqrt{2n}},g\big)U\big(\frac{1}{\sqrt{2n}},f\big)U\big(\frac{1}{\sqrt{2n}},g\big)-U\big(\frac{1}{n},[f,g]\big) \Big\|_{L(H^{s+1},H^s)}<\frac{C}{n}.$$
\end{lemma}
\begin{proof}
   Thanks to Lemma \ref{Lemma:baker}, there exists $O\in L(H^{s+1}; H^s)$ for every $s\geq 0$ such that $$e^{-\frac{1}{ \sqrt{2n}}T_f}e^{-\frac{1}{ \sqrt{2n}}T_g}e^{\frac{1}{ \sqrt{2n}}T_f}e^{\frac{1}{ \sqrt{2n}}T_g}=e^{2\big(\frac{1}{ \sqrt{2n}}\big)^2 T_{[f,g]} + t^3 O}=e^{\frac{1}{ n} T_{[f,g]} + \frac{1}{(2n)^{3/2}} O}.$$
Let $\psi_0\in H^{s+1}(\T^N;\C)$ such that $\|\psi_0\|_{H^{s+1}}=1$ and 
    $$\psi(t)= e^{\frac{t}{ n} T_{[f,g]} + \frac{t}{(2n)^{3/2}} O}\psi_0,$$
$$\varphi(t)=\psi(t)-w(t),\ \ \ \ \ \ \ w(t):=U(n^{-1}t ,[f,g])\psi_0=e^{\frac{t}{n}T_{[f,g]}}\psi_0.$$
The function $\varphi$ solves the following well-posed in $H^s(\T^N;\C)$ partial differential equation
\begin{equation}\begin{split}\label{EDP_1}
\partial_t \varphi
=& \frac{1}{n}T_{[f,g]} \varphi(t) + \frac{1}{(2n)^{3/2}} O (\varphi(t)+w(t)).
\end{split}\end{equation}
We proceed by doing an estimation of the energy in $t\in[0,2]$. Thanks to \eqref{baker-form} and the Young inequality, there exist $C_1,C_2>0$ such that
\begin{align*}
\partial_t\| \varphi\|_{H^s}^2&\leq C_1\Big(\frac{1}{n} \|  \varphi\|_{H^s}+\frac{1}{n}\|\varphi\|_{L^2} +\frac{1}{  n ^{3/2}} \|\varphi\|_{H^s}+\frac{1}{  n ^{3/2}}\|w\|_{H^{s+1}} \Big)\| \varphi\|_{H^s}\leq C_2\Big(\frac{1}{n^2}+\Big(1+\frac{1}{  n  }\Big)\| \varphi\|_{H^s}^2\Big).
\end{align*}
The Gronwall inequality finally ensures the existence of $C>0$ such that
\begin{align*} &\Big\| U\big(-\frac{1}{ \sqrt{2n}},f\big)U\big(-\frac{1}{ \sqrt{2n}},g\big)U\big(\frac{1}{ \sqrt{2n}},f\big)U\big(\frac{1}{ \sqrt{2n}},g\big)\psi_0-U\big(\frac{1}{n},[f,g]\big) \psi_0 \Big\|_{H^s}\\
&=\|\varphi(t=1)\|_{H^s}\leq \frac{1}{n}e^C.\end{align*}
\end{proof}

\begin{lemma}\label{Lemma:commut-L2}
Let $f,g\in \vec(\T^N)$. There exists $C>0$ such that 
$$ \Big\| U\big(-\frac{1}{ \sqrt{2n}},f\big)U\big(-\frac{1}{ \sqrt{2 n}},g\big)U\big(\frac{1}{ \sqrt{2 n}},f\big)U\big(\frac{1}{ \sqrt{2 n}},g\big)-U\big(\frac{1}{n},[f,g]\big) \Big\|_{L(H^1,L^2)}<\frac{C}{n^{3/2}}.$$
\end{lemma}
\begin{proof}
The proof is is consequence of the computations leading to Lemma \ref{Lemma:commut}. We consider $\psi_0\in H^1(\T^N;\C)$ such that $\|\psi_0\|_{H^1}=1$ and define $\varphi$ solving the equation \eqref{EDP_1}. We make an estimation of the energy in $t\in[0,2]$ noticing that $T_{[f,g]}$ is skew-adjoint. The Gronwall inequality concludes the proof since there exist $C_1,C_2>0$ such that
\begin{align*}
\partial_t\| \varphi\|_{L^2}^2&\leq C_1 \frac{1}{n^{3/2}}( \|\varphi\|_{L^2}+\|w\|_{H^1})  \|\varphi\|_{L^2}\leq C_2\Big(\frac{1}{ n^3}+\Big(1+\frac{1}{n^{3/2} }\Big)\|\varphi\|_{L^2}^2\Big).
\end{align*}
\end{proof}

\begin{lemma}\label{Lemma:commut-1}
 Let $s\in \R$ with $s>\frac{N}{2}$ and $f,g\in \vec(\T^N)$. For every $\eta>0$ and $\psi_0\in H^{s}(\T^N;\C)$, there exists $\tilde n\in\N^*$, such that, for every $n>\tilde n$, we have 
$$ \Big\| \Big(U\big(-\frac{1}{ \sqrt{2n}},f\big)U\big(-\frac{1}{ \sqrt{2n}},g\big)U\big(\frac{1}{ \sqrt{2n}},f\big)U\big(\frac{1}{ \sqrt{2n}},g\big)\Big)^n\psi_0-U\big(1,[f,g]\big)\psi_0 \Big\|_{H^s }<\eta.$$
\end{lemma}
\begin{proof}
The statement follows as Lemma \ref{Lemma:sum-1} by considering the lemmas \ref{Lemma:commut} and \ref{Lemma:commut-L2}. We start by taking $\psi_0\in H^{s+2}(\T^N;\C)$ and, thanks to the Lemma \ref{Lemma:commut}, we have
\begin{equation*}\begin{split}
\Big\| \Big(U\big(-\frac{1}{ \sqrt{2 n}},f\big)U\big(-\frac{1}{ \sqrt{2n}},g\big)U\big(\frac{1}{ \sqrt{2n}},f\big)U\big(\frac{1}{ \sqrt{2n}},g\big)\Big)^n\psi_0-U\big(1,[f,g]\big)\psi_0 \Big\|_{H^{s+1} }\leq C\|\psi_0\|_{H^{s+2}}.
    \end{split}\end{equation*}
The same approach with respect to the $L^2$-norm infers thanks to Lemma \ref{Lemma:commut-L2}
\begin{equation*}\begin{split}
\Big\| \Big(U\big(-\frac{1}{ \sqrt{2 n}},f\big)U\big(-\frac{1}{ \sqrt{2n}},g\big)U\big(\frac{1}{ \sqrt{2n}},f\big)U\big(\frac{1}{ \sqrt{2n}},g\big)\Big)^n\psi_0-U\big(1,[f,g]\big)\psi_0 \Big\|_{L^2 }\leq \frac{C}{ n^\frac{1}{2} } \|\psi_0\|_{H^1} .
    \end{split}\end{equation*}
The proof is concluded via an interpolation inequality and a density argument.
\end{proof}

\section{Proof of Proposition \ref{Prop:control-diffeo-2}}\label{proof_diffeo}

Let $f \in H^{s_1}(\T;\R_+^*)$ with $s_1>1/2$ such that $\|f\|_{L^2}=1$. We introduce the diffeomorphism of the torus $\T$
$$Q_f(x)=2\pi\int_{-\pi}^x f^2(t)dt-\pi ,\ \ \ \ (\text{mod}\ 2\pi),$$
and the diffeomorphism given by $P_f=(Q_f)^{-1}$.

%\begin{lemma}
%Let $f \in H^{\frac{s}{2}}(\T;\R_+^*)$ with $\frac{1}{2} \le s < 1$ such that $\|f\|_{L^2}=1$. Then $P_f$ is $H^{1+s}$ and it is isotopic to the identity.
%\end{lemma}
%\begin{proof}
%The regularity of $P_f$ is established by observing that $P_f'(x)=\frac{1}{f^2(x)}$. Now, $f \in C^0(\T;\R_+^*)$ and we can denote $$0<m=\min_{x\in\T} |f(x)|\leq \max_{x\in\T} |f(x)|=M.$$
%The mapping $t\mapsto 1/t$ is Lipschitz in $[m,M]$ with Lipschitz constant equal to $m^{-2}.$ By studying the Gagliardo seminorm, we can then notice that    
%   since, for $P_f(x,s)=\int_0^x (sf^2(t) +(1-s) dt$, we have $P_f(x,1)=P_f(x)$ and $P_f(x,0)=x$. 
%\end{proof}

\begin{lemma}\label{label1}
 Let $f \in H^{s}(\T;\R_+^*)$ and $1/2 < s < 1$ such that $\|f\|_{L^2}=1$. Then $\Lc_{P_f}:H^{s}(\T;\R_+^*)\rightarrow H^{s}(\T;\R_+^*)$.   
\end{lemma}
\begin{proof}
Note that, for every $\psi\in H^{s}(\T;\R_+^*)$, we have $\phi=\psi\circ P_f\in H^{s}(\T;\R_+^*)$. Indeed, if we study its Gagliardo seminorm, then we observe
\begin{align*}|\phi|_{H^{s}}^2=\int_{\T^2}\frac{|\phi(x)-\phi(y)|^2}{|x-y|^{1+2{s}}}dxdy=\int_{\T^2}\frac{|\psi(x)-\psi(y)|^2}{|Q_f(x)-Q_f(y)|^{1+2{s}}}Q_f'(x)Q_f'(y)dxdy.
\end{align*}
Now, $f \in C^0(\T;\R_+^*)$ and for $0<m=\min_{x\in\T} |f(x)|\leq \max_{x\in\T} |f(x)|=M,$ we have 
$$2\pi m^2|x-y|\leq |Q_f(x)-Q_f(y)|=2\pi\Big|\int_y^xf^2(t)dt\Big|\leq 2\pi M^2|x-y|,$$
and then $Q_f$ is bi-Lipschitz. We use this property in the Gagliardo seminorm and we have
\begin{align}\label{inn}|\phi|_{H^{s}}^2\leq \frac{(2\pi)^{1-2{s}}M^4}{m^{2+4{s}}}\int_{\T^2}\frac{|\psi(x)-\psi(y)|^2}{|x-y|^{1+2{s}}}dxdy=\frac{(2\pi)^{1-2{s}}M^4}{m^{2+4{s}}}|\psi|_{H^{s}}^2,
\end{align}
which implies that $\phi\in H^{s}(\T;\R_+^*).$ Now, the mapping $t \mapsto \frac{1}{t}$ is Lipschitz on $[m, M]$ with Lipschitz constant $m^{-2}$, and we define $g := |J_{P_f}|^{1/2} = \frac{1}{\sqrt{2\pi}}\frac{1}{f\circ P_f}$. We see that $g\in H^{s}(\T;\R_+^*)$ since
\begin{align*}|g|_{H^{s}}^2&=\int_{\T^2}\frac{|g(x)-g(y)|^2}{|x-y|^{1+2{s}}}dxdy\leq \frac{1}{(2\pi)m^{4}}\int_{\T^2}\frac{|f\circ P_f(x)-f\circ P_f(y)|^2}{|x-y|^{1+2{s}}}dxdy\\
&= \frac{1}{(2\pi)m^{4}}\int_{\T^2}\frac{|f(x)-f(y)|^2}{|Q_f(x)-Q_f(y)|^{1+2{s}}}Q_f'(x)Q_f'(y)dxdy\leq \frac{ M^4}{ m^{6+4{s}}  (2\pi)^{2{s}}} |f|_{H^{s}}.\end{align*}
The proof is concluded since $g\cdot \phi \in H^{s}(\T;\R)$ as $H^{s}(\T;\R)$ is an algebra.
\end{proof}

\begin{lemma}\label{Prop:control-diffeo-1}
 Let $f \in H^{s_1}(\T;\R_+^*)$ with $s_1>1/2$ such that $\|f\|_{L^2}=1$ and $\frac{1}{2}< s<\min\{s_1,1\} $.  For every $\eta > 0$ and $\psi\in H^{s_1}(\T;\R_+^*)$, there exists $ \D\in \Diff(\T)$ such that 
$$\|\Lc_{ \D}\psi-\Lc_{P_f}\psi\|_{H^{s}} < \eta.$$
\end{lemma}
\begin{proof}
Let $(\rho_n)$ be a sequence of positive mollifiers on $\T$ such that $\tilde f_n = f * \rho_n$ verifies $ \tilde f_n\in C^\infty(\T;\R_+^*)$ and $\tilde f_n\xrightarrow[H^{s_1}]{n\rightarrow\infty}f.$
We consider $f_n=\tilde f_n/\|\tilde f_n\|_{L^2}$ and we have 
$$ f_n\in C^\infty(\T;\R_+^*),\ \ \ \ \ f_n\xrightarrow[H^{s_1}]{n\rightarrow\infty}f,\ \ \ \ \ \|f_n\|_{L^2}=1.$$
Let $n$ large such that $0<m :=\min_{x\in \T}f_n(x) \leq \max_{x\in \T}f_n(x)=:M$. Denote $s^*=\tfrac{s+\min\{s_1,1\}}{2}$ and
\begin{equation}\begin{split}\label{eqq}\|(\Lc_{P_{f_n}}-\Lc_{P_f})\psi\|_{H^{s^*}}&\leq\Big\|\big(|J_{P_{f_n}}|^{1/2}-|J_{P_f}|^{1/2}\big) \psi\circ P_{f_n}\Big\|_{H^{s^*}}+\Big\||J_{P_{f}}|^{1/2} \big(\psi\circ P_{f_n}-\psi\circ P_f\big)\Big\|_{H^{s^*}}\\
&\leq\Big\|\frac{1}{(f_n\circ P_{f_n})(f\circ P_f)}\Big\|_{H^{s^*}} \|f\circ P_f-f_n\circ P_{f_n}\|_{H^{s^*}} \big\|\psi\circ P_{f_n}\big\|_{H^{s^*}}\\
&+\Big\|\frac{1}{f\circ P_f}\Big\|_{H^{s^*}} \big\|\psi\circ P_{f_n}-\psi\circ P_f\big\|_{H^{s^*}}.\end{split}\end{equation}
Note that $\|\psi\circ P_{f_n}\|_{H^{s^*}}$ is uniformly bounded in $n$ from above and below. Indeed, similarly to the proof of Lemma \ref{label1} (equation \ref{inn}), since $Q_{f_n}$ is bi-Lipschitz, we have that 
$$
|\psi\circ P_{f_n}|_{H^{s^*}}^2= \int_{\T^2}\frac{|\psi(x)-\psi (y)|^2}{|Q_{f_n}(x)-Q_{f_n}(y)|^{1+2{s^*}}}Q_{f_n}'(x)Q_{f_n}'(y)dxdy
\begin{cases*}
\leq \frac{(2\pi)^{1-2s}M^4}{m^{2+4{s^*}}}|\psi|_{H^{s^*}}^2,\\
\geq \frac{(2\pi)^{1-2s}m^4}{M^{2+4{s^*}}}|\psi|_{H^{s^*}}^2 .
\end{cases*}
$$
These computations can be also applied to prove that $\|f_n\circ P_{f_n}\|_{H^{s^*}}$ is uniformly bounded in $n$ from above and below since this property verified for $\|f_n\|_{H^{s^*}}$. Hence, from \eqref{eqq}, there exists $C>0$, not depending on $n$, such that
%\begin{align*}\|\psi\circ P_{f_n}-\psi\circ P_f\|_{H^{s}}%^2&=\int_{\T^2}\frac{|\psi\circ P_{f_n}(x)-\psi\circ P_f(x)-\psi\circ P_{f_n}(y)+\psi\circ P_f(y)|^2}{|x-y|^{1+2s}}dxdy\\
%&\leq 2\int_{\T^2}\frac{|\psi\circ P_{f_n}(x)-\psi\circ P_{f_n}(y)|^2+|\psi\circ P_{f}(x)-\psi\circ P_f(y)|^2}{|x-y|^{1+2s}}dxdy\\
%&\leq C\int_{\T^2}\frac{|P_{f_n}(x)-\psi\circ P_{f_n}(y)|^2+|\psi\circ P_{f}(x)-\psi\circ P_f(y)|^2}{|x-y|^{1+2s}}dxdy
%&\leq \|\psi\circ P_{f_n}\|_{H^{s}}+\|\psi\circ P_{f}\|_{H^{s}}
%\leq \big(\|\psi\circ P_{f_n}\|_{H^{s}}+\|\psi\circ P_{f}\|_{H^{s}}\big)
%<C.
%\end{align*}
\begin{equation*}\begin{split}\|(\Lc_{P_{f_n}}-\Lc_{P_f})\psi\|_{H^{s^*}}&\leq C.\end{split}\end{equation*}
We show $\|\psi\circ P_{f_n}-\psi\circ P_f\|_{L^2}\xrightarrow{n\rightarrow\infty}0$ since $\psi\in C^{0,s-1/2}$. Indeed, there exists $C>0$ such that
\begin{align*}\|\psi\circ P_{f_n}-\psi\circ P_f\|_{L^2}^2&%=\int_\T|\psi\circ P_{f_n}(x)-\psi\circ P_f(x)|^2dx
\leq C \int_\T\Big| P_{f_n}(x)-P_{f}(x)\Big|^{2s-1}dx= C \int_\T\Big| P_{f_n}\circ Q_f(x)-x\Big|^{2s-1}Q_f'(x)dx\\
&\leq (2\pi) M^2 C \int_\T\Big| P_{f_n}\circ Q_f(x)-P_{f_n}\circ Q_{f_n}(x)\Big|^{2s-1}dx.\end{align*}
Since $P_{f_n}'= \frac{1}{2\pi(f_n\circ P_{f_n})^2}$, the Lipschitz constant of $P_{f_n}$ is $\frac{1}{2\pi m^2}$ and then 
\begin{align*}\|\psi\circ P_{f_n}-\psi\circ P_f\|_{L^2}^2&\leq \frac{(2\pi)^{2-2s} M^2}{m^{4s-2}} C \int_\T\Big|  \int_{-\pi}^x (f_n^2-f^2)(t)dt\Big|^{2s-1}dx \xrightarrow{n\rightarrow\infty}0.\end{align*}
Similarly, $\|f_n\circ P_{f_n}-f\circ P_f\|_{L^2}\leq \|f_n\circ P_{f_n}-f_n\circ P_f\|_{L^2}+\|f_n\circ P_{f}-f\circ P_f\|_{L^2} \xrightarrow{n\rightarrow\infty}0$ and
\begin{equation*}\begin{split}\|(\Lc_{P_{f_n}}-\Lc_{P_f})\psi\|_{L^2}&\leq\Big\|\frac{1}{f\circ P_f}\Big\|_{L^\infty} \big\|\psi\circ P_{f_n}-\psi\circ P_f\big\|_{L^2}\\
&+\Big\|\frac{1}{(f_n\circ P_{f_n})(f\circ P_f)}\Big\|_{L^\infty} \|f\circ P_f-f_n\circ P_{f_n}\|_{L^2} \big\|\psi\circ P_{f_n}\big\|_{L^\infty}\xrightarrow{n\rightarrow\infty} 0.\end{split}\end{equation*}
By interpolation $\|(\Lc_{P_{f_n}}-\Lc_{P_f})\psi\|_{H^{s}}\xrightarrow{n\rightarrow\infty}0$ for every $\psi\in H^{s_1}(T;\R_+^*)$. %and the result is then extended by density to $ H^{s}(T;\R_+^*)$. 
Finally, $ P_{f_n} $ is isotopic to the identity since we can introduce the continuous family of mappings with $r \in [0,1]$
$$P_{f_n}^r(x):=rP_{f_n}(x)+(1-r)x\ \ \ \ (mod\ 2\pi), $$
verifying $P_{f_n}^1=P_{f_n}$ and $P_{f_n}^0=Id$. Now, $P_{f_n}^r$ defines a family of diffeomorphisms since 
$$(P_{f_n}^r)'(x) = r (P'_{f_n}(x)-1) +1  = r\Big(\frac{1}{(Q_{f_n}'\circ P_{f_n})(x)}-1\Big) +1  
\begin{cases*} 
\geq  r\big(\frac{1}{2\pi M^2}-1\big) +1\geq   \min\big\{1, \frac{1}{2\pi M^2}\big\},\\
\leq r\big(\frac{1}{2\pi m^2}-1\big) +1\leq \max\big\{1, \frac{1}{2\pi m^2}\big\}.
\end{cases*}
$$
The uniform bounds from above and below by positive constants ensure strict monotonicity and hence imply that it is a global diffeomorphism of $\T$.
\end{proof}

\begin{lemma}\label{Prop:control-diffeo-1-bis}
 Let $f \in H^{s_1}(\T;\R_+^*)$ with $s_1>1/2$ such that $\|f\|_{L^2}=1$ and $\frac{1}{2}< s<\min\{s_1,1\}$.  For every $\eta > 0$ and $\psi\in H^{s_1}(\T;\R_+^*)$, there exists $ \D\in \Diff(\T)$ such that 
$$\|\Lc_{ \D}\psi-\Lc_{Q_f}\psi\|_{H^{s}} < \eta.$$
\end{lemma}
\begin{proof}
We proceed as in the proof of Lemma \ref{Prop:control-diffeo-1} and, for $\psi\in H^{s_1}(\T;\R_+^*)$, we have
\begin{equation*}\begin{split}\|(\Lc_{Q_{f_n}}-\Lc_{Q_f})\psi\|_{H^{s^*}}&\leq \| f  \|_{H^{s^*}} \big\|\psi\circ Q_{f_n}-\psi\circ Q_f\big\|_{H^{s^*}}+  \|f -f_n \|_{H^{s^*}} \big\|\psi\circ Q_{f_n}\big\|_{H^{s^*}}.\end{split}\end{equation*}
Since $P_{f_n}'= \frac{1}{2\pi(f_n\circ P_{f_n})^2},$ we have that $P_{f_n}$ is bi-Lipschitz and
$$
|\psi\circ Q_{f_n}|_{H^{s^*}}^2= \int_{\T^2}\frac{|\psi(x)-\psi (y)|^2}{|P_{f_n}(x)-P_{f_n}(y)|^{1+2{s^*}}}P_{f_n}'(x)P_{f_n}'(y)dxdy
\begin{cases*}
\leq \frac{ M^{2+4{s^*}}}{(2\pi)^{1-2{s^*}}m^{4}}|\psi|_{H^{s^*}}^2\\
\geq \frac{ m^{2+4{s^*}}}{(2\pi)^{1-2{s^*}}M^{4}}|\psi|_{H^{s^*}}^2 ,
\end{cases*}
$$
implying that $\|\psi \circ Q_{f_n}\|_{H^{s^*}}$ is uniformly bounded in $n$. 
Consequently, the same bound holds for 
$$
\|(\mathcal{L}_{Q_{f_n}} - \mathcal{L}_{Q_f})\psi\|_{H^{s^*}},
$$
which is also uniformly controlled in $n$. Now, %We show $\|\psi\circ P_{f_n}-\psi\circ P_f\|_{L^2}\xrightarrow{n\rightarrow\infty}0$ since 
$\psi\in C^{0,s-1/2}$ and there exists $C>0$ such that
\begin{align*}\|\psi\circ Q_{f_n}-\psi\circ Q_f\|_{L^2}^2&\leq C \int_\T\Big|\int_{-\pi}^x (f_n^2-f^2)(t)dt\Big|^{2s-1}dx\xrightarrow{n\rightarrow\infty}0.\end{align*}
Again, $\|f_n\circ Q_{f_n}-f\circ Q_f\|_{L^2}\leq \|f_n\circ Q_{f_n}-f_n\circ Q_f\|_{L^2}+\|f_n\circ Q_{f}-f\circ Q_f\|_{L^2} \xrightarrow{n\rightarrow\infty}0$ and
\begin{equation*}\begin{split}\|(\Lc_{Q_{f_n}}-\Lc_{Q_f})\psi\|_{L^2}&\leq \| f  \|_{L^\infty} \big\|\psi\circ Q_{f_n}-\psi\circ Q_f\big\|_{L^2}+  \|f -f_n \|_{L^2} \big\|\psi\circ Q_{f_n}\big\|_{L^\infty}\xrightarrow{n\rightarrow\infty}0.\end{split}\end{equation*}
By interpolation $\|(\Lc_{Q_{f_n}}-\Lc_{Q_f})\psi\|_{H^{s}}\xrightarrow{n\rightarrow\infty}0$ for every $\psi\in H^{s_1}(T;\R_+^*)$. 
Now, $ Q_{f_n} $ is isotopic to the identity since the continuous family of mappings with $r \in [0,1]$
$$Q_f^r(x):=rQ_f(x)+(1-r)x\ \ \ \ (mod\ 2\pi),$$
verifies $Q_{f_n}^1=Q_{f_n}$ and $Q_{f_n}^0=Id$. Finally, $Q_{f_n}^r$ is a family of diffeomorphisms since 
$$(Q_{f_n}^r)'(x) = r (Q'_{f_n}(x)-1) +1  
\begin{cases*} 
\geq  r ( 2\pi m^2 -1) +1\geq   \min \{1,  2\pi m^2 \},\\
\leq r ( 2\pi M^2 -1) +1\leq \max \{1, 2\pi M^2  \}.
\end{cases*}
$$
The uniform bounds from above and below by positive constants ensure strict monotonicity and hence imply that it is a global diffeomorphism of $\T$.
\end{proof}

We are finally, ready to provide the proof of Proposition \ref{Prop:control-diffeo-2}.
\begin{proof}[Proof of Proposition \ref{Prop:control-diffeo-2}]
%We start by introducing the %function $f:=\frac{\psi}{\sqrt{2\pi}}$ and we notice that $f \in H^{s}(\T;\R_+^*)$. %Indeed, as we showed in the proof of Lemma \ref{label1}, the function $1/\varphi\in H^{s}(\T;\R_+^*)$ when $\varphi\in H^{s}(\T;\R_+^*)$ with $1/2<s<1$. Now, $f=\frac{\psi}{\sqrt{2\pi}\varphi} \in H^{s}(\T;\R_+^*)$ since $H^{s}(\T;\R)$ is an algebra. 
We start by considering the diffeomorphism $P_{\frac{\psi}{\sqrt{2\pi}}}$ and we notice that
\begin{align*}\Lc_{P_{\frac{\psi}{\sqrt{2\pi}}}} \psi &=  |J_{P_{\frac{\psi}{\sqrt{2\pi}}}}|^{1/2}\big( \psi\circ P_{\frac{\psi}{\sqrt{2\pi}}}\big) =   \frac{1}{ (\psi\circ P_{\frac{\psi}{\sqrt{2\pi}}})}\big(\psi\circ P_{\frac{\psi}{\sqrt{2\pi}}}\big)=1.
\end{align*}Now, we observe that $\Lc_{Q_{\frac{\varphi}{\sqrt{2\pi}}}\circ P_{\frac{\psi}{\sqrt{2\pi}}}} 
\psi=\varphi$ since $\Lc_{P_{\frac{\varphi}{\sqrt{2\pi}}}}^{-1}\Lc_{P_{\frac{\psi}{\sqrt{2\pi}}}}
\psi=\varphi$ and
$$\Lc_{P_{\frac{\varphi}{\sqrt{2\pi}}}}^{-1}\Lc_{P_{\frac{\psi}{\sqrt{2\pi}}}}=\Lc_{P_{\frac{\varphi}{\sqrt{2\pi}}}^{-1}}\Lc_{P_{\frac{\psi}{\sqrt{2\pi}}}}
=\Lc_{Q_{\frac{\varphi}{\sqrt{2\pi}}}}\Lc_{P_{\frac{\psi}{\sqrt{2\pi}}}}=\Lc_{Q_{\frac{\varphi}{\sqrt{2\pi}}}\circ P_{\frac{\psi}{\sqrt{2\pi}}}} 
 .$$ From the Lemmas \ref{Prop:control-diffeo-1} and \ref{Prop:control-diffeo-1-bis}, for every $\eta_1,\eta_2>0$, there exist $ \D_1,\D_2\in \Diff(\T)$ such that
\begin{align*}\|\Lc_{ \D_1\circ \D_2}\psi-\Lc_{Q_{\frac{\varphi}{\sqrt{2\pi}}}\circ P_{\frac{\psi}{\sqrt{2\pi}}}}\psi\|_{H^{s}}&\leq \|\Lc_{ \D_1\circ \D_2}\psi-\Lc_{ \D_1\circ P_{\frac{\psi}{\sqrt{2\pi}}}}\psi\|_{H^{s}}+\|\Lc_{ \D_1\circ P_{\frac{\psi}{\sqrt{2\pi}}}}\psi-\Lc_{Q_{\frac{\varphi}{\sqrt{2\pi}}}\circ P_{\frac{\psi}{\sqrt{2\pi}}}}\psi\|_{H^{s}}\\
&\leq \|\Lc_{ \D_1}\|_{L(H^{s})} \|\Lc_{ \D_2}\psi-\Lc_{  P_{\frac{\psi}{\sqrt{2\pi}}}}\psi\|_{H^{s}}+\big\|\big(\Lc_{ \D_1}-\Lc_{Q_{\frac{\varphi}{\sqrt{2\pi}}}}\big)\Lc_{P_{\frac{\psi}{\sqrt{2\pi}}}}\psi\big\|_{H^{s}}\\
&\leq \|\Lc_{ \D_1}\|_{L(H^{s})} \eta_2+\eta_1.\end{align*}
The proof ends considering $\D_1$ so that $\eta_1<\frac{\eta}{2}$ at first and after $\D_2$ so that $\eta_2<\frac{\eta}{2\|\Lc_{ \D_1}\|_{L(H^{s})}}.$

\end{proof}

\section{Some results on the solvability of the moments problems}

In this appendix, we present some results on the solvability of moment problems, which are used in various parts of this work. The techniques and the results presented here are not new, and very similar statements can be found in \cite[Appendix B.2]{BL} and \cite[Appendix A]{Duc20b}, for instance. However, since a direct application of such results to our cases is not possible due to the specificity of our problems, we provide their proofs for the sake of completeness.

\begin{proposition}\label{PROP:solvability-mp-lin-0}
Let $I_1,I_2$ be two sets of distinct indices in $\Z^*$ such that  
$$|I_1|=\infty,\ \ \ |I_2|<\infty,\ \ \ \text{and if} \ \ k\in I_1,\ \ \ \text{then}\ \ \ -k\in I_1.$$
Let $(\lambda_k)_{k\in I_1\cup I_2}\subset\R^* $ be a sequence of distinct numbers such that
$$\lambda_k=-\lambda_{-k},\ \ \ \ \text{if}\ \ k\in I_1.$$
Assume $\inf_{\substack{j,k\in I_1\cup I_2\\j\neq k}}|\lambda_{j}-\lambda_{k}|>0$ and
$$\gamma:=\sup_{\substack{K\subset I_1\\K \text{ finite }}}\inf_{\substack{j,k\in I_1\setminus K\\j\neq k}}|\lambda_{j}-\lambda_{k}|.$$ 
%where $K$ runs over the finite subsets of $\Z$.
Let $T>2\pi/\gamma$. We have the following result. For every $x_0,\widehat x_0\in \R,$ and
$$(x_k)_{k\in I_1\cup I_2}\in \left\{(x_k)_{k\in I_1\cup I_2} \in \ell^2(\C)\ :\ x_{-k}=\overline{x_{k}},\ \ k\in I_1;,\ \ x_{k}\in \R,\ \ k\in I_2\right\},$$
there exists $u\in L^2((0,T),\R)$ such that
\begin{equation*}\begin{split}\label{mp-appen} 
&x_{k}=\int_{0}^Tu(t)e^{-i\lambda_k t}dt,\qquad  \forall
k\in I_1,\\
&x_{k}=\int_{0}^Tu(t)e^{\lambda_k t}dt,\qquad  \forall
k\in I_2,\\
x_0&= \int_0^T  t\, u(t)  dt,\ \ \ \ \ 
\widehat x_0=\int_0^T   u(t)  dt.\\
\end{split} \end{equation*}
\end{proposition}
\begin{proof}
Let $T>2\pi/\gamma$ and $\lambda=0$. In the first part of the proof, we show that the functions $\{t,\ e^{i\lambda  t},\ e^{i\lambda_j t},\ e^{\lambda_k t}:j\in I_1,\ k\in I_2\}$ form a Riesz basis of the space
$$X:=\overline{span\{t,\ e^{i\lambda  t},\ e^{i\lambda_j t},\ e^{\lambda_k t}:j\in I_1,\ k\in I_2\}}^{ L^2}\subseteq L^2((0,T);\C).$$

\begin{itemize}[leftmargin=2em, rightmargin=2em]

\item First, consider the sequence $\{\mu_k\}_{k\in I_1\cup I_2}$ obtained by substituting, in the sequence $\{\lambda_k\}_{k\in I_1\cup I_2}$, the frequencies $\{\lambda_k\}_{k\in I_2}$ with some other real numbers $\{\mu_k\}_{k\in I_2}\subset\R^*$ such that
$$\inf_{\substack{j,k\in I_1\cap I_2\\j\neq k}}|\mu_{j}-\mu_{k}|>0.$$
Now, the functions $\{ e^{i\lambda  t},e^{i\mu_k t}\}_{k\in I_1\cup I_2}$ are a Riesz basis in 
$$\overline{span\{e^{i\lambda  t},e^{i\mu_k t}:k\in I_1\cup I_2\}}^{ L^2}\subseteq L^2((0,T),\C),$$
thanks to \cite[Theorem 4.6]{komornik}, and they are therefore minimal. From \cite[Lemma II.4.11]{avdonin}, the functions $\{ e^{i\lambda  t}, \ e^{i\lambda_j t},\ e^{\lambda_k t}:j\in I_1,\ k\in I_2\}$ are also minimal. 

\item Second, we want to show that the family $\{t,\ e^{i\lambda  t}\ e^{i\lambda_j t},\ e^{\lambda_k t}:j\in I_1,\ k\in I_2\}$ is minimal as well, and we follow the techniques developed in the proof of \cite[Corollary 2]{BL}. In detail, we proceed by contradiction and observe that, if 
$$t\in \widetilde X:=\overline{span\{e^{i\lambda  t},\ e^{i\lambda_j t},\ e^{\lambda_k t}:j\in I_1,\ k\in I_2\}}^{ L^2},$$ 
then by integration, for every $j\leq 2$,
$$t^j \in \overline{span\{e^{i\lambda  t},\ e^{i\lambda_j t},\ e^{\lambda_k t}:j\in I_1,\ k\in I_2\}}^{ C^0}\subseteq \widetilde X.$$
Now, $ span\{e^{i\lambda  t}=1 ,t^j:j\in\N^*,\ j\leq 2 \}$ is dense in $C^0([0,T];\C)$ thanks to the Müntz--Szász theorem, and hence also dense in $L^2((0,T);\C)$.  
This property leads to a contradiction, since it would imply 
$$\widetilde X\supseteq \overline{span\{e^{i\lambda  t}=1\ ,t^j:j\in\N^*,\ j\leq 2 \}}^{ L^2}=L^2((0,T);\C),$$
which is impossible since $$e^{i\theta t}\not\in \widetilde X,\ \ \ \ \ \ \ \forall \theta\in\R^*\setminus\{\lambda_k\}_{k\in I_1\cup I_2},$$ thanks to \cite[Theorem 4.6]{komornik}. Finally, the family $\{t,\ e^{i\lambda  t},\ e^{i\lambda_j t},\ e^{\lambda_k t}:j\in I_1,\ k\in I_2\}$ is minimal.

\item Third we observe that \cite[Theorem 4.5]{komornik} ensures that the functions $\{t,\  e^{i\lambda  t}, e^{i\mu_k t}\}_{k\in I_1\cup I_2}$ form a Riesz basis of the space
$$\overline{span\{t,\ e^{i\lambda  t},\ e^{i\mu_k t}:k\in I_1\cup I_2\}}^{ L^2}\subseteq L^2((0,T);\C).$$
By using the Bari's  theorem (see for instance \cite[Theorem 15]{young}), we also have that $\{t,\ e^{i\lambda  t},\ e^{i\lambda_j t},\ e^{\lambda_k t}:j\in I_1,\ k\in I_2\}$ is a Riesz basis in $X$.

\end{itemize}

We are finally ready to use the property that the sequence $\{f_k\}_{k\in\Z}$, obtained by reordering the function $\{t,\ e^{i\lambda  t},\ e^{i\lambda_j t},\ e^{\lambda_k t}:j\in I_1,\ k\in I_2\}$, is a Riesz Basis in $X$. Indeed, this property infers that $$M:u\in X\longmapsto \big\{\jb{
u,f_k}_{L^2}\}_{k\in\Z}\in \ell^2$$ is invertible.
Now, for every $x_0,\widehat x_0\in \R,$ and
$$(x_k)_{k\in I_1\cup I_2}\in \left\{(x_k)_{k\in I_1\cup I_2} \in \ell^2(\C)\ :\ x_{-k}=\overline{x_{k}},\ \ k\in I_1;,\ \ x_{k}\in \R,\ \ k\in I_2\right\},$$ %We define $(y_k)_{k\in\Z}\in\ell^2$ as follows.
%\begin{itemize}[leftmargin=2em, rightmargin=2em]
% \item For every $k\in\Z$ such that $\mu_k=\lambda_j$ with $j\in\Z^*$, we denote $y_k=x_j$.
% \item For every $k\in\Z$ such that $\mu_k=-\lambda_j$ with $j\in\Z^*$, we denote $y_k=\overline{x_j}$.
% \item For $k\in\Z$ such that $\mu_k=\lambda_0$, we denote $y_k=x_0$.
%\end{itemize}
the invertibility of the map $M$ yields the existence of $u\in
X\subseteq L^2((0,T);\C)$
such that
\begin{equation}\label{mp-id-1}\begin{split}
&x_{k}=\int_{0}^Tu(t)e^{-i\lambda_k t}dt,\qquad  \forall
k\in I_1,\\
&x_{k}=\int_{0}^Tu(t)e^{\lambda_k t}dt,\qquad  \forall
k\in I_2,\\
x_0&= \int_0^T  t\, u(t)  dt,\ \ \ \ \ 
\widehat x_0=\int_0^T   u(t)  dt.
\end{split}\end{equation}
To conclude the proof, we just need to show that the control is real by using the properties of the numbers $\{x_k\}_{k\in I_1\cup I_2}$, $x_0$, $\widehat x_0$, and $\{\lambda_k\}_{k\in I_1\cup I_2}$. To this purpose, from \eqref{mp-id-1}, we deduce
\begin{equation}\label{mp-id-2}\begin{split}
&\int_{0}^T {u}(t)e^{ i\lambda_k t} ={x}_{-k}=\overline{x}_{k}=\int_{0}^T\overline{u}(t)e^{ i\lambda_k t}dt,\qquad  \forall
k\in I_1,\\
&\int_{0}^T\Im(u(t))e^{\lambda_k t}dt=\Im({x}_{k})=0 ,\qquad  \forall
k\in I_2,\\
& \int_0^T  t\, \Im(u(t))  dt=\Im(x_0)=0,\ \ \ \ \ 
\int_0^T   \Im(u(t))  dt=\Im(\widehat x_0)=0.
\end{split}\end{equation}
From \eqref{mp-id-2}, we obtain $\jb{\Im(u),f_k}=0$ for every $k\in\Z$. Since $\overline{u}\in X$, it follows that $\Im(u)\in X$, and therefore $\Im(u)=0$, which finally implies $u\in L^2((0,T);\R)$.
\end{proof}

\begin{corollary}\label{PROP:solvability-gen}
Let $I_1,I_2$ be two sets of distinct indices in $\Z^*$ such that  
$$|I_1|=\infty,\ \ \ |I_2|<\infty,\ \ \ \text{and if} \ \ k\in I_1,\ \ \ \text{then}\ \ \ -k\not\in I_1\cup I_2.$$
Assume $\inf_{\substack{j,k\in I_1\cup I_2\\j\neq k}}|\lambda_{j}-\lambda_{k}|>0$, $\inf_{\substack{j,k\in I_1\cup I_2\\j\neq k}}|\lambda_{j}+\lambda_{k}|>0$ and
$$\gamma:=\min\Bigg\{\raisebox{2pt}{$\sup_{\substack{K\subset I_1\\K \text{ finite }}}\inf_{\substack{j,k\in I_1\setminus K\\j\neq k}}|\lambda_{j}-\lambda_{k}|,\ \ \sup_{\substack{K\subset I_1\\K \text{ finite }}}\inf_{\substack{j,k\in I_1\setminus K\\j\neq k}}|\lambda_{j}+\lambda_{k}|$}\Bigg\}.$$ 
%where $K$ runs over the finite subsets of $\Z$.
Let $T>2\pi/\gamma$. We have the following result. For every $x_0,\widehat x_0\in \R,$ and
$$(x_k)_{k\in I_1\cup I_2}\in \left\{(x_k)_{k\in I_1\cup I_2} \in \ell^2(\C)\ :\  x_{k}\in \R,\ \ k\in I_2\right\},$$
there exists $u\in L^2((0,T);\R)$ such that
\begin{equation*}\begin{split}
&x_{k}=\int_{0}^Tu(t)e^{-i\lambda_k t}dt,\qquad  \forall
k\in I_1,\\
&x_{k}=\int_{0}^Tu(t)e^{\lambda_k t}dt,\qquad  \forall
k\in I_2,\\
x_0&= \int_0^T  t\, u(t)  dt,\ \ \ \ \ 
\widehat x_0=\int_0^T   u(t)  dt.\\
\end{split} \end{equation*}
\end{corollary}
\begin{proof}
The result is a direct consequence of Proposition \ref{PROP:solvability-mp-lin-0} by proceeding as follows. 
\begin{itemize}[leftmargin=2em, rightmargin=2em]
    \item We first define the set of indices $\widetilde I_1=I_1\cup \{n\in\Z^*: -n\in I_1\}$.
    \item Second, we consider the sequence $\{\widetilde\lambda_k\}_{k\in \widetilde I_1\cup I_2}$ such that
    $$\widetilde\lambda_k = \lambda_k,\ \ k\in I_1\cup I_2;\qquad \widetilde\lambda_{-k}=-\lambda_{k}\ \ k\in I_1.$$
   % $$=\{\lambda_k\}_{k\in I_1\cup I_2}\cup \{-\lambda_{-k}\}_{k\in I_1}$$
   \item Third, we introduce the numbers $\{\widetilde x_k\}_{k\in \widetilde I_1\cup I_2}$ such that
   $$ \widetilde x_k  =x_k,\ \ k\in I_1\cup I_2;\qquad \widetilde x_{-k}= \overline{x_{k}},\ \ k\in I_1.$$
\end{itemize}
The result is now proved by applying Proposition \ref{PROP:solvability-mp-lin-0} by considering the introduced elements.

\end{proof}


\begin{thebibliography}{100000000}

\bibitem{avdonin}
S.~A. Avdonin, S.~A. Ivanov,
\emph{Families of exponentials: The method of moments in controllability problems for distributed parameter systems},
Cambridge University Press, 1995.

\bibitem{Ball}
J.~M. Ball, J.~E. Marsden, M.~Slemrod,
\emph{Controllability for distributed bilinear systems},
SIAM J. Control Optim. 20 (1982), 575--597.

\bibitem{BV}
V.~Banica, L.~Vega,
\emph{New conservation laws and energy cascade for 1D cubic NLS and the Schrödinger map},
arXiv:2303.08696.

\bibitem{Bea05}
K.~Beauchard,
\emph{Local controllability of a 1D Schrödinger equation},
J. Math. Pures Appl. (9) 84 (2005), no.~7, 851--956.

\bibitem{Karine2}
K.~Beauchard, H.~Lange, H.~Teismann,
\emph{Local exact controllability of a one-dimensional nonlinear Schrödinger equation},
SIAM J. Control Optim. 53 (2015), no.~5, 2781--2818.

\bibitem{BL}
K.~Beauchard, C.~Laurent,
\emph{Local controllability of 1D linear and nonlinear Schrödinger equations with bilinear control},
J. Math. Pures Appl. (9) 94 (2010), no.~5, 520--554.

\bibitem{Karine3}
K.~Beauchard, F.~Marbach, T.~Perrin,
\emph{An obstruction to small-time local controllability for a bilinear Schrödinger equation},
arXiv:2501.03882.

\bibitem{BCP}
K.~Beauchard, R.~Carles, E.~Pozzoli,
\emph{Small-time approximate controllability of the logarithmic Schrödinger equation},
arXiv:2510.14461.

\bibitem{BP04-1}
K.~Beauchard, E.~Pozzoli,
\emph{An example of a small-time globally approximately controllable bilinear Schrödinger equation},
arXiv:2407.05698.

\bibitem{BP04-2}
K.~Beauchard, E.~Pozzoli,
\emph{Small-time approximate controllability of bilinear Schrödinger equations and diffeomorphisms},
arXiv:2410.02383.

\bibitem{BGT_TAMS}
J.~Bernier, B.~Grébert, T.~Robert,
\emph{Dynamics of quintic nonlinear Schrödinger equations in $H^{2/5+}(\T)$}, Trans. Am. Math. Soc. 378, No. 8, 5701-5768 (2025). 

\bibitem{BCMS12}
U.~Boscain, F.~Chittaro, P.~Mason, M.~Sigalotti,
\emph{Adiabatic control of the Schrödinger equation via conical intersections of the eigenvalues},
IEEE Trans. Automat. Control 57 (2012), no.~8, 1970--1983.

\bibitem{BGRS15}
U.~Boscain, J.-P. Gauthier, F.~Rossi, M.~Sigalotti,
\emph{Approximate controllability, exact controllability, and conical eigenvalue intersections for quantum mechanical systems},
Comm. Math. Phys. 333 (2015), no.~3, 1225--1239.

\bibitem{BCMS09}
U.~Boscain, T.~Chambrion, P.~Mason, M.~Sigalotti,
\emph{Controllability of the discrete-spectrum Schrödinger equation driven by an external field},
Ann. Inst. H. Poincaré C Anal. Non Linéaire 26 (2009), no.~1, 329--349.

\bibitem{Bourgain}
J.~Bourgain,
\emph{Fourier transform restriction phenomena for certain lattice subsets and applications to nonlinear evolution equations},
Geom. Funct. Anal. 3 (1993), 107--156.

\bibitem{BBM}
J.~Bourgain, H.~Brezis, P.~Mironescu,
\emph{Lifting in Sobolev spaces},
J. Anal. Math. 80 (2000), 37--86.

\bibitem{BD}
J.~Bourgain, C.~Demeter,
\emph{The proof of the $\ell^2$ decoupling conjecture},
Ann. of Math. 182 (2015), no.~1, 351--389.

\bibitem{UP}
N.~Boussaïd, M.~Caponigro, T.~Chambrion,
\emph{Regular propagators of bilinear quantum systems},
J. Funct. Anal. 278 (2020), no.~6, 108412.

\bibitem{BdCC13}
N.~Boussaïd, M.~Caponigro, T.~Chambrion,
\emph{Weakly coupled systems in quantum control},
IEEE Trans. Automat. Control 58 (2013), no.~9, 2205--2216.

\bibitem{Alessandro}
N.~Boussaid, A.~Duca,
\emph{$H^1$ local exact controllability of some one-dimensional bilinear Schrödinger equations},
arXiv:2504.16641.

\bibitem{BGT}
N.~Burq, P.~Gérard, N.~Tzvetkov,
\emph{Bilinear eigenfunction estimates and the nonlinear Schrödinger equation on surfaces},
Invent. Math. 159 (2005), no.~1, 187--223.

\bibitem{small-time-molecule}
T.~Chambrion, E.~Pozzoli,
\emph{Small-time bilinear control of Schrödinger equations with application to rotating linear molecules},
Automatica 153 (2023), 111028.

\bibitem{coron-2023}
J.-M. Coron, S.~Xiang, P.~Zhang,
\emph{On the global approximate controllability in small time of semiclassical 1D Schrödinger equations},
J. Differential Equations 245 (2023), 1--44.

\bibitem{CHT}
M.~Christ, J.~Holmer, D.~Tataru,
\emph{Low regularity a priori bounds for the modified Korteweg--de Vries equation},
Lib. Math. (N.S.) 32 (2012), no.~1, 51--75.

\bibitem{Guo}
Y. Deng, C. Fan, S. Guo, Z. Guo, Y. Luo,
\emph{Sharp Strichartz estimates on the one-dimensional torus}, arXiv:2609.18670 [math.AP] (2026).

\bibitem{DNY}
Y.~Deng, A.R.~Nahmod,H.~Yue,
\emph{Optimal Local Well-Posedness for the Periodic Derivative Nonlinear Schrödinger Equation}, Commun. Math. Phys. 384 (2021), 1061--1107. 

\bibitem{Duc20b}
A.~Duca,
\emph{Simultaneous global exact controllability in projection of infinite 1D bilinear Schrödinger equations},
Dyn. Partial Differ. Equ. 17 (2020), no.~3, 275--306.

\bibitem{Vahagn2}
A.~Duca, V.~Nersesyan,
\emph{Local exact controllability of the one-dimensional nonlinear Schrödinger equation with Dirichlet boundary conditions},
SIAM J. Control Optim. 63 (2025), no.~1, S20--S36.

\bibitem{komornik}
V.~Komornik, P.~Loreti,
\emph{Fourier series in control theory},
Springer, 2005.

\bibitem{Mor14}
M.~Morancey,
\emph{Simultaneous local exact controllability of 1D bilinear Schrödinger equations},
Ann. Inst. H. Poincaré C Anal. Non Linéaire 31 (2014), no.~3, 501--529.

\bibitem{MN15}
M.~Morancey, V.~Nersesyan,
\emph{Simultaneous global exact controllability of an arbitrary number of 1D bilinear Schrödinger equations},
J. Math. Pures Appl. (9) 103 (2015), no.~1, 228--254.

\bibitem{Vahagn}
A.~Duca, V.~Nersesyan,
\emph{Bilinear control and growth of Sobolev norms for the nonlinear Schrödinger equation},
J. Eur. Math. Soc. 27 (2025), no.~6, 2603--2622.

\bibitem{small-time-momentum}
A.~Duca, E.~Pozzoli,
\emph{Small-time controllability for the nonlinear Schrödinger equation on $\mathbb{R}^N$ via bilinear electromagnetic fields},
SIAM J. Control Optim. 63 (2025), no.~1, S37--S52.

\bibitem{ET}
M.~B. Erdoğan, N.~Tzirakis,
\emph{Dispersive partial differential equations: Well-posedness and applications},
Cambridge Univ. Press, 2016.

\bibitem{Grafakos}
L.~Grafakos,
\emph{Classical Fourier analysis},
3rd ed., Springer, 2014.

\bibitem{GH}
A.~Grünrock, S.~Herr,
\emph{Low regularity local well-posedness of the derivative nonlinear Schrödinger equation},
SIAM J. Math. Anal. 39 (2008), no.~6, 1890--1920.

\bibitem{HHK}
M.~Hadac, S.~Herr, H.~Koch,
\emph{Well-posedness and scattering for the KP-II equation in a critical space},
Ann. Inst. H. Poincaré C Anal. Non Linéaire 26 (2009), no.~3, 917--941.

\bibitem{Haraux}
A.~Haraux,
\emph{Séries lacunaires et contrôle semi-interne des vibrations d'une plaque rectangulaire},
J. Math. Pures Appl. 68 (1989), 457--465.

\bibitem{hardy}
G.~H. Hardy, E.~M. Wright,
\emph{An Introduction to the Theory of Numbers},
6th ed., Oxford University Press, Oxford, 2008.

\bibitem{HGKV}
B.~Harrop-Griffiths, R.~Killip, M.~Visan,
\emph{Scattering for the nonlinear Schrödinger equation with concentrated nonlinearity},
arXiv:2507.14571.

\bibitem{HerrKwak1}
S.~Herr, B.~Kwak,
\emph{Strichartz estimates and global well-posedness of the cubic NLS on $\mathbb{T}^2$},
Forum Math. Pi 12, (2024)

\bibitem{HerrKwak2}
S.~Herr, B.~Kwak,
\emph{Global well-posedness of the cubic nonlinear Schrödinger equation on $\mathbb{T}^2$},
Invent. Math. 245, No. 1, 513-629 (2026). 

\bibitem{HTT}
S.~Herr, D.~Tataru, N.~Tzvetkov,
\emph{Global well-posedness of the energy-critical nonlinear Schrödinger equation},
Duke Math. J. 159 (2011), no.~2, 329--349.

\bibitem{Jaming}
P.~Jaming, K.~Kellay, C.~Saba, Y.~Wang,
\emph{On $L^1$-norms for non-harmonic trigonometric polynomials with sparse frequencies},
J. Fourier Anal. Appl. 31 (2025), no.~2.

\bibitem{KPV}
C.~E. Kenig, G.~Ponce, L.~Vega,
\emph{Quadratic forms for the 1D semilinear Schrödinger equation},
Trans. Amer. Math. Soc. 348 (1996), no.~8, 3323--3353.

\bibitem{KV}
R.~Killip, M.~Visan,
\emph{Scale invariant Strichartz estimates on tori and applications},
Math. Res. Lett. 23 (2016), no.~2, 445--472.

\bibitem{Kishimoto}
N.~Kishimoto,
\emph{Remark on the periodic mass critical nonlinear Schrödinger equation},
Proc. Amer. Math. Soc. 142 (2014), no.~8, 2649--2660.

\bibitem{KT05}
H.~Koch, D.~Tataru,
\emph{Dispersive estimates for principally normal pseudodifferential operators},
Comm. Pure Appl. Math. 58 (2005), no.~2, 217--284.

\bibitem{KT07}
H.~Koch, D.~Tataru,
\emph{A priori bounds for the 1D cubic NLS in negative Sobolev spaces},
IMRN (2007), Art. ID rnm053.

\bibitem{KTV}
H.~Koch, D.~Tataru, M.~Visan,
\emph{Dispersive equations and nonlinear waves},
Birkhäuser, 2014.

\bibitem{Mir09}
M.~Mirrahimi,
\emph{Lyapunov control of a quantum particle in a decaying potential},
Ann. Inst. H. Poincaré C Anal. Non Linéaire 26 (2009), no.~5, 1743--1765.

\bibitem{Mc}
R.~McConnell,
\emph{Nonlinear smoothing for the periodic generalized nonlinear Schrödinger equation},
J. Differential Equations 341 (2022), 353--379.

\bibitem{Ner10}
V.~Nersesyan,
\emph{Global approximate controllability for Schrödinger equation in higher Sobolev norms},
Ann. Inst. H. Poincaré C Anal. Non Linéaire 27 (2010), no.~3, 901--915.

\bibitem{OT}
T.~Ogawa, Y.~Tsutsumi,
\emph{Blow-up of solutions for the nonlinear Schrödinger equation with quartic potential and periodic boundary condition},
Functional-analytic methods for partial differential equations, Proc. Conf. Symp., Tokyo/Jap. 1989, Lect. Notes Math. 1450, 236-251 (1990). 

\bibitem{OW}
T.~Oh, Y.~Wang,
\emph{Normal form approach to the periodic cubic NLS},
J. Anal. Math. 143 (2021), no.~2, 723--762.

\bibitem{Ozawa}
T.~Ozawa,
\emph{Remarks on proofs of conservation laws for nonlinear Schrödinger equations},
Calc. Var. PDE 25 (2006), no.~3, 403--408.

\bibitem{small-time-wave}
E.~Pozzoli,
\emph{Small-time global approximate controllability of bilinear wave equations},
J. Differential Equations 388 (2024), 421--438.

\bibitem{KP1}
T.~Robert,
\emph{Global well-posedness of partially periodic KP-I equation},
Ann. Inst. H. Poincaré C Anal. Non Linéaire 35 (2018), no.~7, 1773--1826.

\bibitem{ModulatedNLS}
T.~Robert,
\emph{Remarks on nonlinear dispersive PDEs with rough dispersion management},
arXiv:2410.23038.

\bibitem{RT}
T.~Robert, L.~Thomann,
\emph{Some topological obstructions to bilinear controllability of nonlinear Schrödinger equations}, in preparation.

\bibitem{Schippa}
R.~Schippa,
\emph{Short-time Fourier transform restriction phenomena and applications to nonlinear dispersive equations},
PhD thesis, 2019.

\bibitem{Serre}
J.-P. Serre,
\emph{Lie algebras and Lie groups},
Springer, 1992.

\bibitem{Tao01}
T.~Tao,
\emph{Multilinear weighted convolution of $L^2$ functions},
Amer. J. Math. 123 (2001), no.~5, 839--908.

\bibitem{Tao06}
T.~Tao,
\emph{Nonlinear dispersive equations: Local and global analysis},
AMS, 2006.

\bibitem{TaoBlog}
T.~Tao,
\emph{The divisor bound},
\url{https://terrytao.wordpress.com/2008/09/23/the-divisor-bound/}.

\bibitem{Triebel3}
H.~Triebel,
\emph{Theory of function spaces III},
Birkhäuser, 2006.

\bibitem{Turinici}
G.~Turinici,
\emph{On the controllability of bilinear quantum systems},
Math. Models Methods Appl. Sci. 74 (2000), 75--92.

\bibitem{Wang}
Y.~Wang,
\emph{Periodic nonlinear Schrödinger equation in critical $H^s(\mathbb{T}^n)$ spaces},
SIAM J. Math. Anal. 45 (2013), no.~3, 1691--1703.

\bibitem{young}
R.~M. Young,
\emph{An Introduction to Non-Harmonic Fourier Series},
Revised edition, Academic Press, 2001.

\end{thebibliography}
\end{document}